\documentclass[a4paper,11pt]{article}
\usepackage{amssymb}
\usepackage{geometry}
\usepackage{amsmath}

\newtheorem{theorem}{Theorem}

\newtheorem{corollary}[theorem]{Corollary}

\newtheorem{definition}[theorem]{Definition}

\newtheorem{lemma}[theorem]{Lemma}

\newtheorem{remark}[theorem]{Remark}

\newtheorem{hypothesis}{Hypothesis}

\newenvironment{proof}[1][Proof]{\textbf{#1.} }{\ \rule{0.5em}{0.5em}}
\newcommand{\RR}{\mathbb{R}}
\newcommand{\loc}{\mathrm{loc}}
\newcommand{\eps}{\varepsilon}

\begin{document}

\title{Well-posedness of the transport equation\\
by stochastic perturbation }
\author{F. Flandoli$^{1}$, M. Gubinelli$^{2}$, E. Priola$^{3}$ \\
{\small {(1) Dipartimento di Matematica Applicata ``U. Dini'', Universit\`a
di Pisa, Italia}}\\
{\small {\ (2) CEREMADE (UMR 7534), Universit\'e Paris Dauphine, France}}\\
{\small {\ (3) Dipartimento di Matematica, Universit\`a di Torino, Italia }}}
\maketitle

\begin{abstract}
We consider the linear transport equation with a globally H\"{o}lder
continuous and bounded vector field, with an integrability condition on the
divergence. While uniqueness may fail for the deterministic PDE, we prove
that a multiplicative stochastic perturbation of Brownian type is enough to
render the equation well-posed. This seems to be the first explicit example
of partial differential equation that become well-posed under the influence
of noise. The key tool is a differentiable stochastic flow constructed and
analyzed by means of a special transformation of the drift of It\^{o}-Tanaka
type.
\end{abstract}

\section{Introduction}

The transport equation in $[0,T]\times \mathbb{R}^{d}$%
\begin{align*}
\partial _{t}u(t,x)+b(t,x)\cdot Du(t,x)& =0,\quad (t,x)\in \lbrack
0,T]\times \mathbb{R}^{d}, \\
u(0,x)& =u_{0}(x),\quad x\in \mathbb{R}^{d}
\end{align*}
driven by a vector field $b:[0,T]\times \mathbb{R}^{d}\rightarrow \mathbb{R}%
^{d}$ which is of class $L^{1}(0,T;W_{loc}^{1,\infty }(\mathbb{R}^{d},%
\mathbb{R}^{d}))$ with a linear growth condition in $x$, can be solved in
several classical ways and in different function spaces. A remarkable
extension has been obtained by R.J. Di Perna and P.L. Lions~\cite
{DiPernaLions} who proved, under the assumption $b\in
L^{1}(0,T;W_{loc}^{1,1}(\mathbb{R}^{d},\mathbb{R}^{d}))$ with a linear
growth condition and ${\mathrm{div}\,}b\in L^{1}(0,T;L^{\infty }(\mathbb{R}%
^{d}))$, that a unique $L^{\infty }([0,T]\times \mathbb{R}^{d})$
solution exists, weak-$\ast $ continuous in time, for any given $u_{0}\in
L^{\infty }(\mathbb{R}^{d})$. Moreover, a generalized notion of flow is
introduced and its existence and uniqueness is proved. L. Ambrosio~\cite
{Ambrosio} proved that uniqueness\ of $L^{\infty }$-solutions is still true
when $W_{loc}^{1,1}(\mathbb{R}^{d},\mathbb{R}^{d})$ is replaced by $BV_{loc}(%
\mathbb{R}^{d},\mathbb{R}^{d})$; furthermore he showed that a sufficient
condition for the uniqueness of the generalized flow is that the negative
part of ${\mathrm{div}\,}b$ is $L^{1}(0,T;L^{\infty }(\mathbb{R}^{d}))$. The
literature following \cite{DiPernaLions} is wide, see a partial review in
\cite{AmbrCrippa}. See also generalizations to transport-diffusion equations
and the associated stochastic differential equations by C. Le Bris and P.L.
Lions \cite{LeBrisLions} and A.~Figalli \cite{Figalli} (in a direction
different from the one of the present paper).

Under weaker conditions on $b$, there are examples of non-uniqueness; see
Section~\ref{sec:positive-example}. The aim of this paper is to show that,
under a suitable random perturbation, $L^{\infty }$-solutions are unique,
even in the case of vector fields $b$ such that the deterministic Cauchy
problem for the transport equation may have multiple solutions. This result
is obtained by introducing a multiplicative white noise in the PDE.
Precisely, we consider the stochastic PDE (SPDE)
\begin{equation}
\begin{aligned} & d_{t}u( t,x) +( b( t,x) \cdot D u( t,x) )
dt+\sum_{i=1}^{d}e_{i}\cdot D u( t,x) \circ dW_{t}^{i}=0,\\ & u( 0,x)
=u_{0}( x),\;\; x \in \mathbb{R}^d, \end{aligned}  \label{SPDE}
\end{equation}
where $e_{1},...,e_{d}$ is the canonical basis of ${\mathbb{R}}^{d}$ and $%
W_{t}=(W_{t}^{1},...,W_{t}^{d})$ is a standard Brownian motion in ${\mathbb{R%
}}^{d}$. The stochastic integration is understood in the Stratonovich sense.

We study existence and uniqueness of $L^{\infty }$-solutions, strong
in the
probabilistic sense, when $b$ is measurable, bounded, globally ${\alpha }$%
-H\"{o}lder continuous in space uniformly in time, for some ${\alpha
}\in (0,1) $ and ${\mathrm{div}\,}b \in L^1_\loc([0,T]\times\RR^d)$.
 In particular, we prove that
  uniqueness always holds  in dimension $d=1$ and
   in any dimension $d \ge 1$
 when $\alpha \in (1/2, 1)$. If $\alpha \in (0, 1/2]$ we still get
uniqueness assuming in addition a global integrability condition on
${\mathrm{div}}\, b$, i.e.,
  \begin{align}\label{fgh}
 {\mathrm{div}\,} b\in L^{p} ( [0,T]\times{\mathbb{R}}^{d} )
 \end{align}
 for some $ p>2$ (recall that a  global integrability
 on ${\mathrm{div}} \, b$  different from \eqref{fgh}
  is  also needed in
the deterministic case to get uniqueness; see
\cite{DiPernaLions,Ambrosio}).

  Moreover, we have existence and uniqueness of $BV_{loc}$%
-solutions, when $u_{0}\in BV_{loc}$, assuming only the H\"{o}lder
condition on $b$, without any assumption on ${\mathrm{div}\,}b$.

Our result gives the first concrete example of partial differential equation
that may lack uniqueness without noise, but is well-posed with a suitable
noise. This phenomenon is well studied for ordinary differential equations
but it is one of the more interesting direction of investigations in the
theory of SPDEs with the ultimate aim of proving the well-posedness of
suitable stochastic perturbations of relevant fluid-dynamics equations. Most
attempts in this direction, until now, focused on additive noise
perturbations, see a discussion in \cite{AFS}. Regularization by
multiplicative noise is, as far as we know, a new phenomenon.

The choice of Stratonovich integral in eq.~\eqref{SPDE} is motivated by two
related facts. On one side (see H. Kunita \cite{K1}), for smooth data and
regular vector field $b$, eq.~\eqref{SPDE} has an explicit solution $%
u(t,x)=u_{0}(\phi _{t}^{-1}(x))$ where $\phi _{t}(x)$ is the flow map giving
the unique strong solution $(X_{t}^{x})_{t\ge 0}$ of the SDE
\begin{equation}
dX_{t}^{x}=b\left( t,X_{t}^{x}\right) dt+dW_{t},\quad t\geq 0,\quad
X_{0}^{x}=x.  \label{SDE intro}
\end{equation}
On the other side, Stratonovich integral is motivated by the Wong-Zakai
principle. Roughly, it state that differential equations driven by regular
random functions usually converge to the Stratonovich version of the
limiting stochastic differential equations provided that these random
functions tend to Brownian motion. In Appendix~\ref{sec:wk} we will prove two
versions of this principle for the stochastic transport equation.

Existence of $L^{\infty }$-solutions to~\eqref{SPDE} does not require all the assumptions stated above on $b$ since a compactness argument require only that $b\in L_{\mathrm{loc}}^{1}(\left[ 0,T\right] \times \mathbb{R}^{d};\mathbb{R}^{d})$
and $\mathrm{div}\,b\in L_{\mathrm{loc}}^{1}(\left[ 0,T\right] \times
\mathbb{R}^{d})$.

As in the deterministic case (see \cite{DiPernaLions,Ambrosio}),
uniqueness of $L^{\infty }$-weak solutions is related to some
form of commutator lemma which allows to perform differential
computations on regularizations of $L^{\infty }$-solutions. In the
deterministic case one has strong convergence to zero of the commutator.
Here we have strong convergence in $L_{loc}^{1}$ only for $d=1$, since~\eqref{fgh} implies that $%
\mathrm{div}\, b = Db\in L_{loc}^{1}\left( \left[ 0,T\right] \times {\mathbb{R}}^{d}\right) $.
For $d>1$ our conditions on $b$ do not guarantee the strong
convergence of the commutator and we are forced to exploit some
non-trivial regularization properties of the stochastic
characteristic equation which come from the non-degeneracy of the
noise. Let us briefly explain this phenomenon.

Observe that formally the mean value $\overline{u}(t,x)=\mathbb{E}[u(t,x)]$
of the solution to eq.~\eqref{SPDE} satisfies
\begin{equation}
\begin{aligned} \frac{d}{dt}\overline{u}( t,x) &+b( t,x) \cdot D \overline
u( t,x) +\frac{1}{2}\Delta \overline u( t,x) dt =0,\\ \overline u( 0,x) &
=u_{0}( x),\;\; x \in \mathbb{R}^d, \end{aligned}  \label{mean-PDE}
\end{equation}
(this may be explained by Lemma \ref{lemma su Ito-Stratonovich}). The
regularizing effect of the viscous term is linked with the regularity
properties of the law of the diffusion $X_{t}^{x}$. At the path-wise level
no regularization of $u(t,x)$ can appear (as witnessed by the characteristics
method). However the non-degeneracy of the diffusion has a remarkable effect
also on integrals of the form
\begin{equation}
\int_{0}^{T}f(s,X_{s}^{x})ds  \label{special integrals}
\end{equation}
where $f:[0,T]\times \mathbb{R}^{d}\rightarrow \mathbb{R}^{d}$ is a \emph{%
deterministic} but possibly time-dependent function. The occupation measure
of a typical trajectory of $X^{x}$ (by occupation measure we mean the
push-forward of Lebesgue measure on $\left[ 0,T\right] $ under the map $%
t\mapsto X_{t}^{x}$) has a density with respect to Lebesgue measure in
dimension one, the local time, see~\cite{RY}. In dimension larger than one
its regularity is less easy, but in any dimension it may be captured by
means of It\^{o} formula and the regularity of solutions of an auxiliary
parabolic equations. Indeed, if we consider a solution $F$ to the parabolic
PDE
\begin{equation*}
\partial _{t}F+\frac{1}{2}\Delta F+b\cdot DF=f
\end{equation*}
on $[0,T]\times \mathbb{R}^{d}$, then an application of the It\^{o} formula
allows us to rewrite the above integral as
\begin{equation*}
\int_{0}^{T}f(s,X_{s}^{x})ds=F(T,X_{T}^{x})-F(0,x)-%
\int_{0}^{T}DF(s,X_{s}^{x})\cdot dW_{s}.
\end{equation*}
The point of this \emph{manouvre} (which we will call the ``It\^{o}--Tanaka
trick'') is to replace the time-average over the diffusion path by a
combination of terms which usually are better behaved than the l.h.s..
Indeed, under appropriate conditions, the non-degeneracy of the diffusion
implies that the solution $F$ of the parabolic PDE is more regular than the
original function $f$. In some sense the It\^{o}--Tanaka trick allows us to
partly transfer the parabolic regularization of the law of the diffusion to
its sample paths.

This basic strategy is our key tool. It gives us the following two
properties of the characteristics equation~\eqref{SDE intro}:

\begin{itemize}
\item[(i)]  Under the assumption that $b(t,\cdot )\in C_{b}^{{\alpha }}({%
\mathbb{R}}^{d};{\mathbb{R}}^{d})$ uniformly in time, equation~%
\eqref{SDE intro} generates a stochastic flow of $C^{1 +
\alpha'}$-diffeomorphisms $\phi _{t}\left( x,\omega \right) $, for
any $0< \alpha' < \alpha$. Moreover this flow is stable under
approximation of the vector field. Related results in $d=1$ have
been proved in \cite{FR}.

\item[(ii)]  Under the integrability assumption
\eqref{fgh} on $\mathrm{div}\,b$, we show that the Jacobian $J\phi
_{t}(\cdot ,\omega )$
of the flow is in $L^{2}(0,T;W_{loc}^{1,2})$ $P$-a.s. (see Theorem \ref{iac}%
). Note that, since the noise is additive, the Jacobian $J\phi
_{t}(x,\omega )$ solves pathwise at least formally the deterministic
ODE
\begin{equation*}
\frac{d}{dt}J\phi _{t}(x,\omega )={\mathrm{div}\,}b(t,\phi _{t}(x,\omega
))J\phi _{t}(x,\omega )
\end{equation*}
and thus
\begin{equation*}
\log J\phi _{t}(x,\omega )=\int_{0}^{t}{\mathrm{div}\,}b(s,\phi
_{s}(x,\omega ))ds.
\end{equation*}
In the deterministic case there is no hope to differentiate in $x$
without further differentiability assumptions on
${\mathrm{div}\,}b$, but in the stochastic case we use again the
improved regularity of integrals of the form (\ref{special
integrals}) thanks to the It\^{o}--Tanaka trick and to some
 classical $L^p$-parabolic regularity results.
  The
required integrability condition on ${\mathrm{div}\,}b$ (see
\eqref{fgh}) turns out to be different from the one imposed in the
deterministic setting \cite {DiPernaLions,Ambrosio}.
\end{itemize}
 When $\alpha \in (1/2, 1)$, we can prove distributional
   convergence of the
 commutator (and so uniqueness for our SPDE)
  by combining    the H\"older regularity  of
   $b$ and   the one of the stochastic flow (so using only (i)).
In the case of $\alpha \in (0, 1/2]$    we can still prove
distributional
 convergence of the
 commutator, but we need  both  (i) and (ii) (and so we
 require \eqref{fgh}).

The differentiability of the flow gives easily existence and
uniqueness of $BV_{\text{loc}}$ solutions (and other more regular
spaces) to the SPDE without any requirement on the divergence (see
Appendix~\ref{sec:bv}).

\medskip A remark on the connections with the work of Y.~Le Jan and
O.~Raimond~\cite{LR} on generalized stochastic flows is important. At
present, no precise comparison can be made between the result of the present
work and those of \cite{LR}, but it is clear that \cite{LR} has been a
source of inspiration for us, like \cite{Ambrosio, DiPernaLions}.

Conceptual similarities between all these works can be seen in the following
results of \cite{LR}. It deals with a stochastic transport-like equation, or
to be more precise, a stochastic continuity equation written in weak form
and a variation-of-constant reformulation of it, see equations (e) and (d)
respectively of Theorem 3.2. An existence and uniqueness result is proved,
in a special class of solutions, under very general assumptions. Finally,
\cite{LR} gives criteria for existence of an associated flow of maps, or on
the contrary for the possibility of coalescence and diffusion. The general
results are applied to examples where the coefficients have a very poor
Sobolev regularity.

What is entirely different between our work and \cite{LR} is that \cite{LR}
deals with stochastic equations which are well posed in the weak sense
(martingale sense) but not necessarily strongly well posed. The relevant
examples, at present, are constructed by means of suitable diffusion
coefficients (and infinite dimensional noise, often, like the isotropic
Brownian motion);\ the drift part does not play a relevant role. The poor
regularity of coefficients mentioned above regards the diffusion
coefficients. On the contrary, our purpose is to deal with a non-regular
drift coefficient (and we choose a trivial but non-degenerate diffusion part
for sake of simplicity), following the philosophy that we randomly perturb a
deterministic transport equation having non-regular drift. Our stochastic
equations are also strongly well posed and they always define a stochastic
flow of diffeomorphisms.

Finding a synthesis of these different approaches to non-regular transport
(or continuity) equations, deterministic and stochastic, would be a very
interesting progress. Let us end this introduction by mentioning a few other
open questions.

The generalization to nonlinear transport equations, where $b$ depends on $u$
itself, would be a major next step for applications to fluid dynamics but it
turns out to be a difficult problem. Specifically there are already some
difficulties in dealing with a vector field $b$ which depends itself on the
random perturbation $W$. There is no obvious extension of the It\^{o}--Tanaka
trick to integrals of the form $\int_{0}^{T}f(\omega ,s,X_{s}^{x}\left(
\omega \right) )ds$ with random $f$. As we will show in Section~\ref
{sec:examples}, it is very easy to produce examples, both for the linear
SPDE (\ref{SPDE}) and for a stochastic version of Euler equation which show
that the particular noise we use does not have any regularizing effect in
this case. Thus new ideas are needed to approach nonlinear problems.

The linear case with deterministic $b$ still contains interesting open
problems. N. Depauw \cite{Dep} gave examples of non-uniqueness in $d=2$ for
divergence free bounded measurable fields $b$ with a condition on the bounded
variation norm. Our results do not cover this case. In particular, the
existence of stochastic flows under $L^{\infty }$ assumptions on $b$ is an
interesting open problem.

Finally, in Section~\ref{sec:examples} we see that the classical
one dimensional example with $b\left( x\right) \simeq \left|
x\right| ^{\gamma }$, $\gamma \in \left( 0,1\right) $, is covered by our
uniqueness result. The differentiablility of the stochastic flow $\phi$ means in
particular that its stretching $J\phi $  around $x=0$ is very
large but finite. However, we can prove $J\phi \in L^{2}(0,T;W_{\mathrm{loc}%
}^{1,2})$ only for $\gamma \in \left( 1/2,1\right) $. It is not clear if $%
\gamma \ge 1/2$ is a natural threshold for the smoothness of the
stretching or it
is just a limitation of our approach.

\paragraph{Acknowledgement.}
The authors would like to thank the anonymous referees for the
careful reading of the first version of this paper and for their
remarks which helped to greatly improve the paper.

\paragraph{Plan.}

The main body of the paper is devoted to the analysis of weak $L^{\infty }$%
-solutions and preliminaries on stochastic flows. In Sect.~\ref{sec:ito} we
prove the existence of a global stochastic flow associated to eq.~(\ref{SDE
intro}) and its differentiability properties. Sect.~\ref{sec:jacobian} is
devoted to prove that under our hypotheses on $\mathrm{div}\,b$ the Jacobian
of the flow is in $L^{2}(0,T;W_{\mathrm{loc}}^{1,2}(\mathbb{R}^{d}))$. In
Sect.~\ref{sec:transp} we prove existence of weak solutions to the SPDE~(\ref
{SPDE}). Sect.~\ref{sec:uniq} is devoted to prove uniqueness of $L^{\infty }$%
-solutions to the SPDE~(\ref{SPDE}). Finally, in
Sect.~\ref{sec:examples} we collect some positive and negative
examples.

Then we present a number of appendixes on related results. Appendix
\ref {sec:bv} is devoted to existence and uniqueness of $BV_{loc}$
solutions, Appendix \ref{sect perturbative} to an equivalent
pathwise formulation of the SPDE,
Appendix \ref{sec:wk} gives Wong-Zakai approximation results
finally Appendix \ref{sec:fractional} gives  uniqueness  results
 by fractional Sobolev spaces  non covered in the main
 text.

\paragraph{Notations.}

\label{sect notations}

Usually we denote by $D_{i}f$ the derivative in the $i$-th coordinate
direction and with $(e_{i})_{i=1,\dots ,d}$ the canonical basis of $\mathbb{R%
}^{d}$ so that $D_{i}f=e_{i}\cdot Df$. For partial derivatives of any order $%
n\ge 1$ we use the notation $D_{i_{1},...,i_{n}}^{n}$. If $\eta :\mathbb{R}%
^{d}\rightarrow \mathbb{R}^{d}$ is a $C^{1}$-diffeomorphism we will denote
by $J\eta (x)=\text{det}[D\eta (x)]$ its Jacobian determinant. For a given
function $f$ depending on $t\in \lbrack 0,T]$ and $x\in {\mathbb{R}}^{d}$,
we will also adopt the notation $f_{t}(x)=f(t,x)$.

Let $T>0$ be fixed. For ${\alpha }\in (0,1)$ define the space $L^{\infty
}\left( 0,T;C_{b}^{\alpha }({\mathbb{R}}^{d})\right) $ as the set of all
bounded Borel functions $f:[0,T]\times {\mathbb{R}}^{d}\rightarrow {\mathbb{R%
}}$ for which
\begin{equation*}
\lbrack f]_{\alpha ,T}=\sup_{t\in \lbrack 0,T]}\sup_{x\neq y\in {\mathbb{R}}%
^{d}}\frac{|f(t,x)-f(t,y)|}{|x-y|^{\alpha }}<\infty
\end{equation*}
($|\cdot |$ denotes the Euclidean norm in ${\mathbb{R}}^{d}$ for every $d$,
if no confusion may arise). This is a Banach space with respect to the usual
norm $\Vert f\Vert _{{\alpha },T}=\Vert f\Vert _{0}+[f]_{\alpha ,T}$ where $%
\Vert f\Vert _{0}=\sup_{(t,x)\in \lbrack 0,T]\times {\mathbb{R}}%
^{d}}|f(t,x)| $.

We write $L^{\infty }\left( 0,T;C_{b}^{\alpha }({\mathbb{R}}^{d};{\mathbb{R}}%
^{d})\right) $ for the space of all vector fields $f:[0,T]\times {\mathbb{R}}%
^{d}\rightarrow {\mathbb{R}}^{d}$ having all components in $L^{\infty
}\left( 0,T;C_{b}^{\alpha }({\mathbb{R}}^{d})\right) $.

Moreover, for $n\geq 1,$ $f\in L^{\infty }\left( 0,T;C_{b}^{n+\alpha }({%
\mathbb{R}}^{d})\right) $ if all spatial partial derivatives $%
D_{i_{1},...,i_{k}}^{k}f\in L^{\infty }\left( 0,T;C_{b}^{\alpha }({\mathbb{R}%
}^{d})\right) $, for all orders $k=0,1,\dots ,n$. Define the corresponding
norm as
\begin{equation*}
\Vert f\Vert _{n+\alpha ,T}=\Vert f\Vert _{0}+\sum_{k=1}^{n}\Vert
D^{k}f\Vert _{0}+[D^{n}f]_{\alpha ,T}
\end{equation*}
where we extend the previous notations $\Vert \cdot \Vert _{0}$ and $[\cdot
]_{\alpha ,T}$ to tensors. The definition of the space $L^{\infty }\left(
0,T;C_{b}^{n+\alpha }({\mathbb{R}}^{d};{\mathbb{R}}^{d})\right) $ is
similar. The previous functions spaces can be defined similarly when $%
T=+\infty $ (i.e., we are considering functions defined on $[0,\infty
)\times \mathbb{R}^{d}$). The spaces $C_{b}^{n+{\alpha }}({\mathbb{R}}^{d})$
and $C_{b}^{n+{\alpha }}({\mathbb{R}}^{d};{\mathbb{R}}^{d})$ are defined as
before but only involve functions $f:{\mathbb{R}}^{d}\rightarrow {\mathbb{R}}%
^{d}$ which do not depend on time. Moreover, we say that $f:\mathbb{R}%
^{d}\rightarrow \mathbb{R}^{d}$ belongs to $C^{n,\alpha }$, $n\in \mathbb{N} $, $\alpha
\in (0,1)$, if $f$ is continuous on $\mathbb{R}^{d}$, $n$-times
differentiable with all continuous derivatives and the derivatives of order $%
n$ are locally $\alpha$-H\"{o}lder continuous. Finally, $C_{0}^{0}(\mathbb{R%
}^{d})$ denotes the space of all real continuous functions defined on $%
\mathbb{R}^{d}$, having compact support and by $C_{0}^{\infty }(\mathbb{R}%
^{d})$ its subspace consisting of infinitely differentiable functions.

For any $r>0$ we denote by $B(r)$ the Euclidean ball centered in 0 of radius
$r$ and by $C_{r}^{\infty }(\mathbb{R}^{d})$ the space of smooth
functions with compact support in $B(r)$; moreover, $\Vert \cdot \Vert
_{L_{r}^{p}}$ and $\Vert \cdot \Vert _{W_{r}^{1,p}}$ stand for,
respectively, the $L^{p}$-norm and the $W^{1,p}$-norm on $B\left( r\right)
$, $p\in \left[ 1,\infty \right] $.
We let also $[f]_{C^\theta_r}=\sup_{x\neq  y \in B(r)}|f(x)-f(y)|/|x-y|^\theta$.

We will often use the standard mollifiers\textit{.} Let $\vartheta :{\mathbb{%
R}}^{d}\rightarrow {\mathbb{R}}$ be a smooth test function such that $0\leq
\vartheta (x)\leq 1$, $x\in {\mathbb{R}}^{d}$, $\vartheta (x)=\vartheta (-x)$%
, $\int_{{\mathbb{R}}^{d}}\vartheta (x)dx=1$, $\mathrm{supp}\,(\vartheta
)\subset $ $B(2)$, $\vartheta (x)=1$ when $x\in B(1)$. For any ${\varepsilon
}>0$, let $\vartheta _{{\varepsilon }}(x)={\varepsilon }^{-d}\vartheta (x/{%
\varepsilon })$ and for any distribution $g:{\mathbb{R}}^{d}\rightarrow {%
\mathbb{R}}^{n}$ we define the mollified approximation $g^{\varepsilon }$ as
\begin{equation}
g^{\varepsilon }(x)=\vartheta _{{\varepsilon }}\ast g(x)=g(\vartheta _{{%
\varepsilon }}(x-\cdot )),\;\;\;x\in {\mathbb{R}}^{d}.  \label{molli}
\end{equation}
If $g$ depends also on time $t$, we consider $g^{\varepsilon
}(t,x)=(\vartheta _{{\varepsilon }}\ast g(t,\cdot ))(x)$, $t\in \lbrack 0,T]$%
, $x\in \mathbb{R}^{d}$.

\medskip Recall that, for any smooth bounded domain $\mathcal{D}$ of $\mathbb{R}^d$,
we have: $f \in W^{\theta, p} (\mathcal  D)$, $\theta \in (0,1)$, $p \ge 1$, if and
only if $f \in L^p (\mathcal D)$ and
\begin{equation*}
[f]_{W^{\theta, p} }^p =\iint_{\mathcal D \times \mathcal D} \frac{|f(x) - f(y)|^p}{%
|x-y|^{\theta p + d}} dx dy < \infty.
\end{equation*}
We have $W^{1, p} (\mathcal D) \subset W^{\theta, p} (\mathcal D)$, $\theta \in (0,1)$.

\bigskip Throughout the paper we will assume a stochastic basis with a $d$%
-dimensional Brownian motion $\left( \Omega,\left( \mathcal{F}{}_{t}\right)
,{}\mathcal{F},P,\left( W_{t}\right) \right) $ to be given. We denote by $%
\mathcal{F}_{s,t}$ the completed $\sigma$-algebra generated by $W_{u}-W_{r}$%
, $s\leq r\leq u\leq t$, for each $0\le s<t$.

\section{Differentiable stochastic flow with $C_{b}^{{\protect\alpha }}$
drift}

\label{sec:ito}

Given $s\in \left[ 0,T\right] $ and $x\in {\mathbb{R}}^{d}$, consider the
stochastic differential equation in ${\mathbb{R}}^{d}$ :
\begin{equation}
dX_{t}=b\left( t,X_{t}\right) dt+dW_{t},\quad t\in \left[ s,T\right] ,\quad
X_{s}=x.  \label{SDE}
\end{equation}
A classical fact is the existence and uniqueness of a weak solution,
obtained for instance by Girsanov transform. Yu.~Veretennikov \cite{V}
proved that boundedness of $b$ (uniformly in $t$) is enough to have
path-wise uniqueness and existence of strong probabilistic solutions. For
related works see also the more recent paper~\cite{Kry-Ro} by N.V.~Krylov
and M.~R\"{o}ckner where strong uniqueness is proved under some
integrability assumption on $b$. These works are based on the technique
introduced by Zvonkin~\cite{Zvonkin} of removing the irregular drift by a
suitable change of coordinates in the SDE. The fact that such a coordinate
change modifies the drift is a consequence of the It\^{o} formula.
Technically, as we will see in the proof of Th.~\ref{th:flow1}, the
It\^{o}--Tanaka trick is similar to the Zvonkin approach. It is worthwhile
to note however that the heuristic behind is not necessarily the same and in
our opinion the It\^{o}--Tanaka point of view has wider range of
applicability (as we demonstrate in the control of the Jacobian of the flow).

A related interesting result has been obtained by A.M.~Davie in~\cite{D}.
Under the assumption that $b$ is measurable and bounded, he proved that the
(deterministic) integral equation
\begin{equation*}
x(t) = x + \int_0^t b(s,x(s)) ds + w(t)
\end{equation*}
has a unique solution $x(\cdot) \in C(0,T;\mathbb{R}^d)$ for all $w\in N^c$
where $N\subset C(0,T;\mathbb{R}^d)$ is a set which has probability zero
according to Wiener measure. This paper contains also the very interesting
key estimate
\begin{equation*}
E \left[\left| \int_0^1 (b(s,x+W_s)-b(s,W_s)) ds\right|^p\right] \le C_{p}
\|b\|_\infty |x|^p, \qquad x\in \mathbb{R}^d
\end{equation*}
where $C_p$ is an absolute constant not depending on $b$. This estimate is
obtained by non-trivial direct computations and show very explicitly the
regularization phenomenon which occurs when considering average values of
functions along the trajectories of diffusions (Brownian motion in this
case).

In all the cited works the analysis of the flows is however missing,
essentially they deal only with (various forms of) path-wise uniqueness of
the SDE. For papers that tackle existence of global flows of homeomorphisms
for SDEs without global Lipschitz coefficients see~\cite{Zhang, FIZ} and the
references therein. However, the assumptions of these works are too strong
for our purposes.

\bigskip

Our key result is the existence of a \emph{differentiable} stochastic flow $%
(x,s,t)\mapsto \phi _{s,t}(x)$ for equation~(\ref{SDE}) under the following
hypothesis:

\begin{hypothesis}
\label{hy1} There exists ${\alpha }\in (0,1)$ such that $b\in L^{\infty
}\left( 0,T;C_{b}^{\alpha }({\mathbb{R}}^{d};{\mathbb{R}}^{d})\right) $.
\end{hypothesis}

Recall the relevant definition from~\cite{K}:

\begin{definition}
A \emph{stochastic flow of diffeomorphisms} (resp. \emph{of class} $%
C^{1,\alpha }$) on $(\Omega ,\left( \mathcal{F}{}_{t}\right) ,$ ${}\mathcal{F%
},P,\left( W_{t}\right) )$ associated to equation (\ref{SDE}) is a map $%
(s,t,x,\omega )\mapsto \phi _{s,t}(x)\left( \omega \right) $, defined for $%
0\leq s\leq t\leq T$, $x\in {\mathbb{R}}^{d}$, $\omega \in \Omega $ with
values in ${\mathbb{R}}^{d}$, such that

\begin{itemize}
\item[(a)]  given any $s\in \left[ 0,T\right] $, $x\in {\mathbb{R}}^{d}$,
the process $X^{s,x}=(X_{t}^{s,x},t\in \lbrack s,T])$ defined as $%
X_{t}^{s,x}=\phi _{s,t}(x)$ is a continuous $\mathcal{F}_{s,t}$-measurable
solution of equation (\ref{SDE}),

\item[(b)]  $P$-a.s., $\phi _{s,t}$ is a diffeomorphism, for all $0\leq
s\leq t\leq T$, and the functions $\phi _{s,t}(x)$, $\phi _{s,t}^{-1}(x)$, $%
D\phi _{s,t}(x)$, $D\phi _{s,t}^{-1}(x)$ are continuous in $(s,t,x)$ (resp.
of class $C^{\alpha }$ in $x$ uniformly in $(s,t)$),

\item[(c)]  $P$-a.s., $\phi _{s,t}(x)=\phi _{u,t}(\phi _{s,u}(x))$ for all $%
0\leq s\leq u\leq t\leq T$ and $x\in {\mathbb{R}}^{d}$ and $\phi _{s,s}(x)=x$%
.
\end{itemize}
\end{definition}

As already mentioned, the main ingredient to obtain the regularity of the
flow is the observation that the time integral $\int_{0}^{t}b\left(
s,X_{s}^{x}\right) ds$ has richer regularity properties than expected only
on the basis of the regularity of $b$. To reveal them we have to use the
regularity theory of parabolic PDEs. Let $\lambda >0$ be fixed. Let us
extend $b$ to the whole $[0,\infty )\times {\mathbb{R}}^{d}$ by setting
\begin{equation}
b(t,x)=b(T,x),\;\;\;t\geq T,\;\;x\in {\mathbb{R}}^{d}.  \label{bb}
\end{equation}
Clearly we have that $b\in L^{\infty }\left( 0,\infty ;C_{b}^{\alpha }({%
\mathbb{R}}^{d};{\mathbb{R}}^{d})\right) $. Given a function $f\in L^{\infty
}\left( 0,\infty ;C_{b}^{\alpha }({\mathbb{R}}^{d};{\mathbb{R}}^{d})\right) $%
, consider the following backward parabolic system (collecting $d$
independent equations):
\begin{equation}
\partial _{t}u_{\lambda }+L^{b}u_{\lambda }-\lambda u_{\lambda
}=f,\;\;\;(t,x)\in \lbrack 0,\infty )\times {\mathbb{R}}^{d},
\label{eq:parabolic}
\end{equation}
where
\begin{equation}
L^{b}u=\frac{1}{2}\Delta u+b\cdot Du  \label{fr}
\end{equation}
and $u:[0,+\infty )\times {\mathbb{R}}^{d}\rightarrow {\mathbb{R}}^{d}$ (%
eq.~\eqref{fr} has to be interpreted componentwise).

Since $b$ and $f$ are only measurable in time instead of continuous, the
notion of solution to \eqref{eq:parabolic} is not standard. We follow~\cite{KP} by prescribing that a function $u:[0,+\infty )\times {\mathbb{R}%
}^{d}\rightarrow {\mathbb{R}}^{d}$ which belongs to $L^{\infty }\left(
0,\infty ;C_{b}^{2+\alpha }({\mathbb{R}}^{d};{\mathbb{R}}^{d})\right) $ is a
solution to \eqref{eq:parabolic} if
\begin{equation*}
u(t,x)-u(s,x)=\int_{s}^{t}[-L^{b}u(r,x)+\lambda u(r,x)+f(r,x)]dr
\end{equation*}
for every $t\geq s\geq 0$, $x\in {\mathbb{R}}^{d}$. From this identity it
follows that $u\left( \cdot ,x\right) $ is Lipschitz continuous for every $%
x\in {\mathbb{R}}^{d}$.
Other regularity properties can be found in \cite{KP}.

\smallskip The next result deals with Schauder estimates and is known even
in a more general form (see \cite{KP} and the references therein). A-priori
estimates of the type \eqref{stima} were first proved in \cite{Br}. We will
only sketch the proof and refer to \cite{KP} for more details. The backward
equation \eqref{eq:parabolic} is not supplemented by the value of the limit $%
u\left( \infty ,x\right) $ and uniqueness is due to the condition of uniform
boundedness of $u$. 

\begin{theorem}
\label{schaud} Let us consider equation \eqref{eq:parabolic} with $b,f\in
L^{\infty }\left( 0,\infty ;C_{b}^{\alpha }({\mathbb{R}}^{d};{\mathbb{R}}%
^{d})\right) .$ Then there exists a unique solution $u=u_{\lambda }$ to
equation \eqref{eq:parabolic} in the space $L^{\infty }\left( 0,\infty
;C_{b}^{2+\alpha }({\mathbb{R}}^{d};{\mathbb{R}}^{d})\right) $. Moreover
there exists $C>0$ (independent on $u$ and $f$) such that
\begin{equation}
\sup_{t\geq 0}\Vert u(t,\cdot )\Vert _{C_{b}^{2+\alpha }({\mathbb{R}}^{d};{%
\mathbb{R}}^{d})}\leq C\sup_{t\geq 0}\Vert f(t,\cdot )\Vert _{C_{b}^{\alpha
}({\mathbb{R}}^{d};{\mathbb{R}}^{d})}  \label{stima}
\end{equation}
\end{theorem}

\begin{proof}
\textbf{Step 1} (uniqueness). Uniqueness follows from the maximum principle $%
\Vert u\Vert _{0}\leq \lambda ^{-1}\Vert f\Vert _{0}$ applied to the
difference of two solutions. The maximum principle under our conditions is
proved in \cite[Theorem 4.1]{KP} (the proof is more delicate than in the
classical case when $b$ and $f$ are continuous in $t$, see \cite[Theorem 8.1.7]{Kr0}). For completeness we give also a self-contained probabilistic proof. From
\cite{V} or \cite{Kry-Ro}, under Hypothesis~\ref{hy1} there exists a unique
strong solution $(X_{t}^{s,x})_{t\ge 0}$ of equation (\ref{SDE}). Let $u\in L^{\infty
}\left( 0,\infty ;C_{b}^{2+\alpha }({\mathbb{R}}^{d};{\mathbb{R}}%
^{d})\right) $ be a solution of equation \eqref{eq:parabolic}. For any given $s$, we may apply
It\^{o} formula to $e^{-\lambda \left( t-s\right) }u\left(
t,X_{t}^{s,x}\right) $ (see Lemma~\ref{lemma formula Ito} below) in the $t$ variable. Taking then expectation we get
\begin{equation*}
e^{-\lambda \left( t-s\right) }E\left[ u\left( t,X_{s,t}^{x}\right) \right]
=u\left( s,x\right) +\int_{s}^{t}e^{-\lambda \left( r-s\right) }E\left[
f\left( r,X_{s,r}^{x}\right) \right] dr.
\end{equation*}
As $t\rightarrow \infty $ (recall that $u$ is bounded) we obtain
\begin{equation*}
u\left( s,x\right) =-\int_{s}^{\infty }e^{-\lambda \left( r-s\right) }E\left[
f\left( r,X_{s,r}^{x}\right) \right] dr.
\end{equation*}
Then $u$ is uniquely determined by $f$ and $b$. This also gives the estimate
$\Vert u\Vert _{0}\leq \lambda ^{-1}\Vert f\Vert _{0}$ mentioned above.

\textbf{Step 2} (existence and estimate~\eqref{stima}). We only recall the
idea of the proof, see \cite{KP} for details. If $b=0$, then the result
follows by using the explicit formula
\begin{equation}
u(t,x)=\int_{t}^{+\infty }e^{-\lambda (r-t)}P_{r-t}f(r,\cdot
)(x)dr=\int_{0}^{+\infty }e^{-\lambda s}P_{s}f(t+s,\cdot )(x)ds  \label{heat}
\end{equation}
(where $(P_{t})$ denotes the forward heat semigroup) and well known
estimates on the spatial derivatives of $P_{t}g$ when $g\in C_{b}^{\alpha
}\left( {\mathbb{R}}^{d}\right) $ and $t>0$.

In the general case, using to the boundedness of $b$ and the maximum
principle, we get easily a-priori estimates for equation \eqref{eq:parabolic}
(assuming that there exists a bounded solution $u$). Then a continuity
method (see~\cite[Lemma 4.3]{KP}) allows to get the existence of the
solution which verifies equation (\ref{eq:parabolic}), along with the estimate~\eqref{stima}.
\end{proof}

In the previous proof we have used It\^{o} formula for solutions of equation %
\eqref{eq:parabolic}, although their regularity in time is not standard for
It\^{o} formula. We give a self-contained proof of the validity of It\^{o} formula in our hypotheses since we have to use it again below in the essential step of the
change of variables from the SDE (\ref{SDE})\ to the SDE (\ref{conjugated
SDE}).

\begin{lemma}
\label{lemma formula Ito}Let $u:[0,+\infty )\times {\mathbb{R}}%
^{d}\rightarrow {\mathbb{R}}$ be a function of class $L^{\infty }\left(
0,\infty ;C_{b}^{2+\alpha }({\mathbb{R}}^{d})\right) $, such that
\begin{equation}
U\left( t,x\right) -U\left( s,x\right) =\int_{s}^{t}V\left( r,x\right) dr
\label{U V}
\end{equation}
for every $t\geq s\geq 0$, $x\in {\mathbb{R}}^{d}$, with $V\in L^{\infty
}\left( 0,\infty ;C_{b}^{\alpha }({\mathbb{R}}^{d})\right) $. Let $\left(
X_{t}\right) _{t\geq 0}$ be a continuous adapted process of the form
\begin{equation*}
X_{t}=X_{0}+\int_{0}^{t}b_{s}ds+\int_{0}^{t}\sigma _{s}dW_{s}
\end{equation*}
where $b$ and $\sigma $ are (resp. ${\mathbb{R}}^{d}$-valued and ${\mathbb{R}%
}^{d\times d}$-valued) progressively measurable processes, $b$ integrable
and $\sigma $ square integrable in $t$ with probability one. Then
\begin{equation*}
U\left( t,X_{t}\right) =U\left( 0,x\right) +\int_{0}^{t}\left( V+b_{s}\cdot
DU+\frac{1}{2}Tr\left( \sigma \sigma ^{T}D^{2}U\right) \right) \left(
s,X_{s}\right) ds
\end{equation*}
\begin{equation}
+\int_{0}^{t}\left\langle DU\left( s,X_{s}\right) ,\sigma
_{s}dW_{s}\right\rangle .  \label{Ito generalized}
\end{equation}
\end{lemma}

\begin{proof}
Set
\begin{equation*}
U_{\varepsilon }\left( t,x\right) =\varepsilon ^{-1}\int_{t}^{t+\varepsilon
}U\left( s,x\right) ds,\quad V_{\varepsilon }\left( t,x\right) =\varepsilon
^{-1}\int_{t}^{t+\varepsilon }V\left( s,x\right) ds.
\end{equation*}
The time derivative (see (\ref{U V}))
\begin{equation*}
\partial _{t}U_{\varepsilon }=\varepsilon ^{-1}\left[ U\left( t+\varepsilon
,x\right) -U\left( t,x\right) \right] =V_{\varepsilon }\left( t,x\right)
\end{equation*}
exists and is continuous. Thus $U_{\varepsilon }$ satisfies the assumptions
of the classical It\^{o} formula. We apply it to $U_{\varepsilon }\left(
t,X_{t}\right) $ and get an identity like (\ref{Ito generalized}) with $%
U_{\varepsilon }$ and $V_{\varepsilon }$ in place of $U$ and $V$. Given $%
t\geq 0$, the r.v. $U_{\varepsilon }\left( t,X_{t}\left( \omega \right)
\right) =\varepsilon ^{-1}\int_{t}^{t+\varepsilon }U\left( s,X_{t}\left(
\omega \right) \right) ds$ converges $P$-a.s. to $U\left( t,X_{t}\left(
\omega \right) \right) $ as $\varepsilon \rightarrow 0$ (we may also use the
fact that, from equation (\ref{U V}), $U$ is globally bounded and continuous
in $(t,x)$ on $[0,+\infty )\times \mathbb{R}^{d}$). Now we use the fact that
$\left| V_{\varepsilon }\right| $, $\left| U_{\varepsilon }\right| $, $%
\left\| DU_{\varepsilon }\right\| $, $\left\| D^{2}U_{\varepsilon }\right\| $
are uniformly bounded in $\left( t,x,\varepsilon \right) $. The $P$-a.s.
convergence of the Lebesgue integral is easy by dominated convergence
theorem, since $\left| V_{\varepsilon }\right| $ is uniformly bounded in $%
\left( t,x,\varepsilon \right) $, $\left| b_{s}\cdot DU_{\varepsilon
}\right| $ is uniformly bounded in $\left( x,\varepsilon \right) $ by a
constant times $\left| b_{s}\right| $ and similarly for $\left| Tr\left(
\sigma \sigma ^{T}D^{2}U_{\varepsilon }\right) \right| $. Finally, the r.v. $%
\int_{0}^{t}\left\langle DU_{\varepsilon }\left( s,X_{s}\right) ,\sigma
_{s}dW_{s}\right\rangle $ converges in probability to $\int_{0}^{t}\left%
\langle DU\left( s,X_{s}\right) ,\sigma _{s}dW_{s}\right\rangle $ because $%
\int_{0}^{t}\left\| \sigma _{s}^{T}\left( DU_{\varepsilon }\left(
s,X_{s}\right) -DU\left( s,X_{s}\right) \right) \right\| ^{2}ds$ converges
in probability to zero (since it converges to zero $P$-a.s., again because
the integrand is bounded by a constant times $\left\| \sigma _{s}\sigma
_{s}^{T}\right\| $). The proof is complete.
\end{proof}

We also need the following simple lemma.

\begin{lemma}
\label{grad} Under the assumptions of Theorem \ref{schaud}, let $u_{\lambda }
$ be the solution to (\ref{eq:parabolic}). Then
\begin{equation*}
\Vert Du_{\lambda }\Vert _{0}\rightarrow 0,\;\;\;as\;\;\lambda \rightarrow
+\infty
\end{equation*}
where the supremum is taken on $[0,\infty )\times {\mathbb{R}}^{d}$. The
choice of $\lambda $ to have, for instance, $\Vert Du_{\lambda }\Vert
_{0}\leq 1/2$, depends only on $\Vert b\Vert _{0}$ and $\Vert f\Vert _{0}$.
\end{lemma}

\begin{proof}
We write $\partial _{t}u_{\lambda }+\frac{1}{2}\Delta u_{\lambda }-\lambda
u_{\lambda }=f-b\cdot Du_{\lambda }$. Using the well known estimate for the
heat semigroup $\sup_{x\in \mathbb{R}^{d}}|DP_{t}g(x)|\leq
Ct^{-1/2}\sup_{x\in \mathbb{R}^{d}}|g(x)|$, $g\in C_{b}({\mathbb{R}}^{d}),$ $%
t>0,$ and, differentiating in formula (\ref{heat}), we get easily, for any $%
\lambda >0$,
\begin{equation*}
\Vert Du_{\lambda }\Vert _{0}\leq \frac{c}{\lambda ^{\frac{1}{2}}}(\Vert
f\Vert _{0}+\Vert b\cdot Du_{\lambda }\Vert _{0})\leq \frac{c}{\lambda ^{%
\frac{1}{2}}}\Vert f\Vert _{0}+\frac{c}{\lambda ^{\frac{1}{2}}}\Vert b\Vert
_{0}\Vert Du_{\lambda }\Vert _{0}.
\end{equation*}
Considering $\lambda >c^{2}\Vert b\Vert _{0}^{2}$, we get
\begin{equation*}
(1-\frac{c}{\lambda ^{\frac{1}{2}}}\Vert b\Vert _{0})\,\Vert Du_{\lambda
}\Vert _{0}\leq \frac{c}{\lambda ^{\frac{1}{2}}}\Vert f\Vert _{0}
\end{equation*}
and the assertion follows. The proof is complete.
\end{proof}

\begin{theorem}
\label{th:flow1} Assume $b\in L^{\infty }\left( 0,T;C_{b}^{\alpha }({\mathbb{%
R}}^{d};{\mathbb{R}}^{d})\right) $. Then we have the following facts:

\begin{itemize}
\item[(i)]  (pathwise uniqueness) For every $s\in \left[ 0,T\right] $, $x\in
{\mathbb{R}}^{d}$, the stochastic equation (\ref{SDE}) has a unique
continuous adapted solution $X^{s,x}=\left( X_{t}^{s,x}\big(\omega \right)
,t\in \left[ s,T\right] ,$ $\omega \in \Omega \big)$.

\item[(ii)]  (differentiable flow) There exists a stochastic flow $\phi
_{s,t}$ of diffeomorphisms for equation (\ref{SDE}). The flow is also of
class $C^{1+{\alpha }^{\prime }}$ for any ${\alpha }^{\prime }<{\alpha }$.

\item[(iii)]  (stability) Let $(b^{n})\subset L^{\infty }\left(
0,T;C_{b}^{\alpha }({\mathbb{R}}^{d};{\mathbb{R}}^{d})\right) $ be a
sequence of vector fields and $\phi ^{n}$ be the corresponding stochastic
flows. If $b^{n}\rightarrow b$ in $L^{\infty }\left( 0,T;C_{b}^{\alpha
^{\prime }}({\mathbb{R}}^{d};{\mathbb{R}}^{d})\right) $ for some $\alpha
^{\prime }>0$, then, for any $p\geq 1$,
\begin{equation}
\lim_{n\rightarrow \infty }\sup_{x\in {\mathbb{R}}^{d}}\sup_{0\leq s\leq
T}E[\sup_{r\in \lbrack s,T]}|\phi _{s,r}^{n}(x)-\phi _{s,r}(x)|^{p}]=0
\label{stability1}
\end{equation}
\begin{equation}
\sup_{n\in \mathbb{N}}\sup_{x\in {\mathbb{R}}^{d}}\sup_{0\leq s\leq
T}E[\sup_{u\in \lbrack s,T]}\Vert D\phi _{s,u}^{n}(x)\Vert ^{p}]<\infty ,
\label{bound}
\end{equation}
\begin{equation}
\lim_{n\rightarrow \infty }\sup_{x\in {\mathbb{R}}^{d}}\sup_{0\leq s\leq
T}E[\sup_{r\in \lbrack s,T]}\Vert D\phi _{s,r}^{n}(x)-D\phi _{s,r}(x)\Vert
^{p}]=0.  \label{stability2}
\end{equation}
\end{itemize}
\end{theorem}


\begin{proof}
\textbf{Step 1.} (auxiliary parabolic systems). For $\lambda >0$ consider
the (vector valued) solution $\psi \in L^{\infty }\left( 0,\infty
;C_{b}^{2+\alpha }({\mathbb{R}}^{d};{\mathbb{R}}^{d})\right) $ to the
parabolic system
\begin{equation}
\partial _{t}\psi _{\lambda }+L^{b}\psi _{\lambda }-\lambda \psi _{\lambda
}=-b,\;\;\;(t,x)\in (0,\infty )\times {\mathbb{R}}^{d}  \label{pde}
\end{equation}
provided by Theorem~\ref{schaud} with $f=-b$. Define
\begin{equation*}
\Psi _{\lambda }(t,x)=x+\psi _{\lambda }(t,x).
\end{equation*}

\begin{lemma}
\label{diff} For $\lambda $ large enough, such that $\sup_{t\geq 0}\Vert
D\psi _{\lambda }(t,\cdot )\Vert _{0}<1$ (see Lemma~\ref{grad}), the
following statements hold:

{\vskip2mm \noindent } (i) Uniformly in $t\in \lbrack 0,+\infty )$, $\Psi
_{\lambda }$ has bounded first and second spatial derivatives and moreover
the second (Fr\'{e}chet) derivative $D_{x}^{2}\Psi _{\lambda }$ is globally $%
\alpha $-H\"{o}lder continuous.

{\vskip2mm \noindent } (ii) \ for any $t\geq 0$, $\Psi _{\lambda }:\mathbb{R}%
^{d}\rightarrow \mathbb{R}^{d}$ is a non-singular diffeomorphism of class $%
C^{2}$.

{\vskip2mm \noindent } (iii) $\Psi _{\lambda }^{-1}$ has bounded first and
second spatial derivatives, uniformly in $t\in \lbrack 0,+\infty )$.
\end{lemma}

\begin{proof}
Assertion (i) follows by Theorem \ref{schaud}.

\vskip2mm (ii) Recall the classical Hadamard theorem (see for instance
\cite[page 330]{protter}): Let $g:\mathbb{R}^{d}\rightarrow \mathbb{R}^{d}$
be of class $C^{2}.$ Suppose that $\lim_{|x|\rightarrow \infty
}|g(x)|=+\infty $ and that the Jacobian matrix $Dg(x)$ is an isomorphism of $%
\mathbb{R}^{d}$ for all $x\in \mathbb{R}^{d}$. Then $g$ is a $C^{2}$%
-diffeomorphism of $\mathbb{R}^{d}$. Applying this result to $\Psi _{\lambda
}$, we get the assertion.

\vskip2mm (iii) We know that $\Psi _{\lambda }^{-1}$ is of class $C^{2}$.
Moreover
\begin{equation*}
D\Psi _{\lambda }^{-1}(t,y)=[D\Psi _{\lambda }(t,\Psi _{\lambda
}^{-1}y)]^{-1}=[I+D\psi _{\lambda }(t,\Psi _{\lambda }^{-1}y)]^{-1}
\end{equation*}
\begin{equation*}
=\sum_{k\geq 0}(-D\psi _{\lambda }(t,\Psi _{\lambda }^{-1}y))^{k},\;\;\;y\in
\mathbb{R}^{d}.
\end{equation*}
It follows that $\sup_{t\geq 0}\Vert D\Psi _{\lambda }^{-1}(t,\cdot )\Vert
_{0}\leq \sum_{k\geq 0}\big (\sup_{t\geq 0}\Vert D\psi _{\lambda }(t,\cdot
)\Vert _{0}\big )^{k}<\infty $. This shows the boundedness of the first
derivative. Arguing in a similar way we get also the boundedness of the
second derivative since
\begin{equation*}
D^{2}\Psi _{\lambda }^{-1}(t,y)=-[D\Psi _{\lambda }(t,\Psi _{\lambda
}^{-1}y)]^{-1}D^{2}\Psi _{\lambda }(t,\Psi _{\lambda }^{-1}(t,y))\left\{
[D\Psi _{\lambda }(t,\Psi _{\lambda }^{-1}y)]^{-1}\right\} ^{\otimes 2}.
\end{equation*}
\end{proof}

\vskip2mm In the sequel we will use a value of $\lambda $ for which Lemma~\ref{diff} holds and simply write $\psi $ and $\Psi $ for $\psi _{\lambda }$
and $\Psi _{\lambda }$.

\smallskip \textbf{Step 2.} (conjugated SDE). Define
\begin{equation*}
\widetilde{b}(t,y)=-\lambda \psi (t,\Psi ^{-1}(t,y)),\quad \widetilde{\sigma
}(t,y)=D\Psi (t,\Psi ^{-1}(t,y))
\end{equation*}
and consider, for every $s\in \left[ 0,T\right] $ and $y\in {\mathbb{R}}^{d}$%
, the SDE
\begin{equation}
Y_{t}=y+\int_{s}^{t}\tilde{\sigma}(u,Y_{u})dW_{u}+\int_{s}^{t}\widetilde{b}%
(u,Y_{u})du,\qquad t\in \lbrack s,T].  \label{conjugated SDE}
\end{equation}
This equation is equivalent to equation (\ref{SDE}), in the following
sense.\ If $X_{t}$ is a solution to \eqref{SDE}, then $Y_{t}=\Psi (t,X_{t})$
verifies equation (\ref{conjugated SDE}) with $y=\Psi (s,x)$: it is
sufficient to apply It\^{o} formula of Lemma~\ref{lemma formula Ito} to $%
\Psi (t,X_{t})$ and use equation (\ref{pde}). It is also possible to show
that given a solution $Y_{t}$ of equation (\ref{conjugated SDE}), then $%
X_{t}=\Psi ^{-1}(t,Y_{t})$ is a solution of \eqref{SDE} with $x=\Psi
^{-1}(s,y)$, but we shall not use this fact.

\smallskip \textbf{Step 3.} (proof of (i) and (ii)). Assertion (i) is known,
see \cite{V}, but we give a proof based on our approach. We have clearly $%
\widetilde{b}\in L^{\infty }\left( 0,T;C_{b}^{1+\alpha }({\mathbb{R}}^{d};{%
\mathbb{R}}^{d})\right) $ and $\widetilde{\sigma }\in L^{\infty }\left(
0,T;C_{b}^{1+\alpha }({\mathbb{R}}^{d};{\mathbb{R}}^{d\times d})\right) $.
By classical results (see \cite[Ch. 2]{K}) this implies existence and
uniqueness of a strong solution $Y$ of equation (\ref{conjugated SDE}) and
even the existence of a $C^{1,\alpha ^{\prime }}$ ($\alpha ^{\prime }<\alpha
$) stochastic flow of diffeomorphisms $\varphi _{s,t}$ associated to
equation (\ref{conjugated SDE}).
Continuity in time is assumed
in~\cite{K}, but it can be easily extended to $L^\infty$ time dependence, as it is done in~\cite{K2} even in greater generality.
The uniqueness of $Y$ implies the path-wise
uniqueness of solutions of the original SDE~\eqref{SDE} since two solutions $%
X,\tilde{X}$ give rise to two processes $Y_{t}=\Psi (t,X_{t})$ and $\tilde{Y}%
_{t}=\Psi (t,\tilde{X}_{t})$ solving~\eqref{conjugated SDE}, then $Y=\tilde{Y%
}$ and then necessarily $X=\tilde{X}$. By the Yamada-Watanabe theorem
path-wise uniqueness together with weak existence (which is a direct
consequence of the Girsanov formula) gives the existence of the (unique)
solution $(X_{t}^{x})_{t\geq s}$ of eq.~\eqref{SDE} starting from $x$ at
time $s$. Moreover setting $\phi _{s,t}=\Psi _{t}^{-1}\circ \varphi
_{s,t}\circ \Psi _{s}$ we realize that $\phi _{s,t}$ is the flow of~%
\eqref{SDE} (in the sense that $X_{t}^{x}=\phi _{s,t}(x))$.

\smallskip \textbf{Step 4.} (proof of (iii)). Let $\psi ^{n}$ be the
solution in $L^{\infty }\left( 0,T;C_{b}^{2+\alpha }({\mathbb{R}}^{d};{%
\mathbb{R}}^{d})\right) $ of the parabolic problem (\ref{pde}) associated to
$b_{n}$. Notice that we can make a choice of $\lambda $ independent of $n$.
Since $b^{n}\rightarrow b$ in $L^{\infty }\left( 0,T;C_{b}^{\alpha }({%
\mathbb{R}}^{d};{\mathbb{R}}^{d})\right) $, by Theorem~\ref{schaud} we have $%
\psi ^{n}\rightarrow \psi $ in $L^{\infty }\left( 0,T;C_{b}^{2+\alpha }({%
\mathbb{R}}^{d};{\mathbb{R}}^{d})\right) $. To prove this last fact one has to write
\begin{equation*}
\partial _{t}\left( \psi ^{n}-\psi \right) +L^{b}\left( \psi ^{n}-\psi
\right) -\lambda \left( \psi ^{n}-\psi \right) =-\left( b^{n}-b\right)
+\left( b^{n}-b\right) \cdot D\psi ^{n}
\end{equation*}
and use the bound
\begin{equation*}
\left\| \left( b^{n}-b\right) \cdot \left( D\psi ^{n}-I\right) \right\|
_{C_{b}^{\alpha }({\mathbb{R}}^{d};{\mathbb{R}}^{d})}\leq C\left\|
b^{n}-b\right\| _{C_{b}^{\alpha }({\mathbb{R}}^{d};{\mathbb{R}}^{d})}
\end{equation*}
which is true since, by Theorem~\ref{schaud}, $D\psi ^{n}$ is uniformly
bounded.

Consider the flows $\varphi _{s,t}^{n}=\Psi _{t}^{n}\circ \phi
_{s,t}^{n}\circ (\Psi _{s}^{n})^{-1}$ which satisfy
\begin{equation}
\varphi _{s,t}^{n}(y)=y+\int_{s}^{t}\widetilde{b}^{n}(u,\varphi
_{s,u}^{n}(y))du+\int_{s}^{t}\widetilde{\sigma }^{n}(u,\varphi
_{s,u}^{n}(y))\cdot dW_{u},
\end{equation}
We have $\widetilde{\sigma }^{n}\rightarrow \widetilde{\sigma }$ and $%
\widetilde{b}^{n}\rightarrow \widetilde{b}$ in $L^{\infty }\left(
0,T;C_{b}^{1+\alpha }({\mathbb{R}}^{d};{\mathbb{R}}^{d\times d})\right) $
and $L^{\infty }\left( 0,T;C_{b}^{1+\alpha }({\mathbb{R}}^{d};{\mathbb{R}}%
^{d})\right) $, respectively. By standard argument using the Gronwall lemma,
the Doob inequality and the Burkholder inequality 
(compare, for instance, with the proof of \cite[Theorem II.2.1]{K}) we
obtain the analog of \eqref{stability1} for the auxiliary flows $\varphi
_{s,t}^{n}$ and $\varphi _{s,t}$. The estimates are standard so we leave them to the reader. We note only that we need to control the difference
$\varphi _{s,t}^{n}(y_{n})-\varphi _{s,t}(y)$ where $y_{n}=\Psi _{s}^{n}( x) $ and $y=\Psi _{s}( x)$. But $| \Psi _{s}^{n}( x) -\Psi
_{s}( x) | $ is  uniformly small in $x\in {\mathbb{R}}^{d}$
so we need to estimate $\varphi _{s,t}^{n}(y+v)-\varphi _{s,t}(y)$,
uniformly in $y$, only with respect to a uniformly small variation $v$. Then it is not
difficult to see that $E[ \sup_{t\in [ s,T] }| \varphi
_{s,t}^{n}(y+v)-\varphi _{s,t}(y)| ^{p}] $ is small in
$y$ uniformly for large $n$ and small $v$. Finally, one has to check that $%
(\Psi _{s}^{n})^{-1}$ converges to $\Psi _{s}^{-1}$ uniformly. This is due
to the fact that $\Psi _{s}^{n}$ converges uniformly to $\Psi _{s}$ with its
derivatives and the Jacobian $J\Psi _{s}$ is uniformly away from zero.

Concerning the derivative of the stochastic flow, first one can prove an
inequality for $D\varphi _{s,t}^{n}(y)$ similar to~\eqref{bound} using the
fact that the equation for $D\varphi _{s,t}^{n}(y)$ has the identity as
initial condition and the coefficients $D\widetilde{b}^{n}\left( \phi
_{s,u}^{n}\right) $ and $D\widetilde{\sigma }^{n}\left( \phi
_{s,u}^{n}\right) $ are uniformly bounded functions (in all variables and $n$%
). Then one has to use the uniform boundedness of the derivatives of $\Psi
_{s}^{n}$ and its inverse, to estimate $D\phi _{s,u}^{n}$.

Arguing as in the proof of \cite[Theorem II.3.1]{K}, we get the following
linear equation for the derivative $D\phi _{s,t}(x)$
\begin{equation}
\begin{split}
\lbrack D\Psi _{t}(\phi _{s,t}(x))]D\phi _{s,t}(x)& =D\Psi
_{s}(x)+\int_{s}^{t}[D^{2}\Psi _{u}(\phi _{s,u}(x))]D\phi _{s,u}(x)dW_{u} \\
& -\lambda \int_{s}^{t}[D\psi _{u}(\phi _{s,u}(x))]D\phi _{s,u}(x)du,
\end{split}
\label{dove}
\end{equation}
$0\leq s\leq t\leq T$, $x\in \mathbb{R}^{d}$. From the convergence $\psi
^{n}\rightarrow \psi $ in $L^{\infty }\left( 0,T;C_{b}^{2+\alpha }({\mathbb{R%
}}^{d};{\mathbb{R}}^{d})\right) $ together with \eqref{bound} and %
\eqref{dove}, we finally obtain, for any $p\geq 1$,
\begin{equation}
\lim_{n\rightarrow \infty }\sup_{x\in {\mathbb{R}}^{d}}\sup_{0\leq s\leq
T}E[\sup_{u\in \lbrack s,T]}\Vert D\phi _{s,u}^{n}(x)-D\phi _{s,u}(x)\Vert
^{p}]=0.  \label{stima1}
\end{equation}
\end{proof}

\vskip2mm We show now that the inverse flow ${\phi }_{s,t}^{-1}$ is directly
related to the solutions of a simple backward stochastic differential
equations, of the same form as the original one (only the drift has opposite
sign).

\begin{lemma}
\label{fra} \bigskip Under Hypothesis~\ref{hy1} the process $\left( \phi
_{s,t}^{-1}\left( y\right) \right) _{s\in \left[ 0,t\right] }$ is the unique
solution of the backward SDE
\begin{equation*}
\phi _{s,t}^{-1}(y)=y-\int_{s}^{t}b(r,\phi _{r,t}^{-1}(y))dr-[W_{t}-W_{s}].
\end{equation*}
and
\begin{equation}
\sup_{x\in {\mathbb{R}}^{d}}\sup_{0\leq u\leq T}E[\sup_{s\in \lbrack
0,u]}\Vert D\phi _{s,u}^{-1}(x)\Vert ^{p}]<\infty   \label{bound-inv}
\end{equation}
for any $p\geq 1$.
\end{lemma}

\begin{proof}
We have $\phi _{s,t}(x)=x+\int_{s}^{t}b(r,\phi _{s,r}(x))dr+W_{t}-W_{s}$ and
then
\begin{equation*}
\phi _{s,t}(\phi _{s,t}^{-1}(y))=\phi _{s,t}^{-1}(y)+\int_{s}^{t}b(r,\phi
_{s,r}(\phi _{s,t}^{-1}(y)))dr+W_{t}-W_{s}.
\end{equation*}
But $\phi _{s,r}\circ \phi _{s,t}^{-1}=\phi _{r,t}^{-1}$ and thus $y=\phi
_{s,t}^{-1}(y)+\int_{s}^{t}b(r,\phi _{r,t}^{-1}(y))dr+W_{t}-W_{s}$. The
proof of the bound~(\ref{bound-inv}) is then similar to that of eq.~(\ref
{bound}) once taken into account the backward character of the equation.
\end{proof}

\begin{remark}
\label{civuole} We shall use the following consequence of the stability
estimates in Th.~\ref{th:flow1}: for any $r>0$ and $p\geq 1$,
\begin{equation*}
\lim_{n\rightarrow \infty }E\left[ \int_{B\left( r\right) }\sup_{t\in
\lbrack 0,T]}|\phi _{t}^{n}(y)-\phi _{t}(y)|^{p}dy\,\right] =0
\end{equation*}
\begin{equation*}
\lim_{n\rightarrow \infty }E\left[ \int_{B\left( r\right) }\sup_{t\in
\lbrack 0,T]}\Vert D\phi _{t}^{n}(y)-D\phi _{t}(y)\Vert ^{p}dy\,\right] =0.
\end{equation*}
Possibly passing to a subsequence, still denoted by $\phi ^{n}$, this
implies that $P$-a.s. the following property holds: for every $t\in \left[
0,T\right] $, $\phi _{t}^{n}$ and $D\phi _{t}^{n}$ converge to $\phi _{t}$
and $D\phi _{t}$ in $L^{p}\left( B\left( r\right) ;{\mathbb{R}}^{d}\right) $
and $L^{p}\left( B\left( r\right) ;{\mathbb{R}}^{d\times d}\right) $
respectively.
\end{remark}

We finish the section with a result of independent interest. It concerns a
special situation when the given vector field $b$ is of zero distributional
divergence and gives rise to a measure-preserving flow.

\begin{lemma}
\label{lemma change} Assume ${\mathrm{div}\,}b(t,\cdot )=0$, $t\in \lbrack
0,T],$ in the sense of distributions. Then the stochastic flow $\phi $ is $P$%
-a.s. measure-preserving, i.e., $J\phi _{s,t}(x)=1$, for all $0\leq s\leq
t\leq T$ and all $x\in {\mathbb{R}}^{d}$, $P$-a.s..
\end{lemma}

\begin{proof}
Let $b_{n}$ be $C_{b,T}^{\infty }$ vector fields that converge to $b$ in $%
C_{b,T}^{{\alpha ^{\prime }}}$ for some $0<\alpha ^{\prime }<\alpha $ and
such that ${\mathrm{div}\,}b_{n}=0$. The functions $b_{n}$ can be
constructed as in \eqref{molli}. Let $\phi ^{n}$ be the associated smooth
diffeomorphism. Applying \cite[Theorem II.3.1]{K} and the well known
Liouville theorem, we get that the diffeomorphism $\phi ^{n}$ preserves the
Lebsegue measure since $b_{n}$ is a divergence-free vector field.
%
%
%
%
%
%
%
Then $J\phi _{s,t}^{n}(x)=1$ for all $0\leq s<t\leq T$ and all $x$, $P$%
-a.s.. Fix $x\in {\mathbb{R}}^{d}$ and $s\in \lbrack 0,T]$. By \eqref{stima1}%
, there exists a subsequence (possibly depending on $x,s$ and still denoted
by $D\phi _{s,u}^{n}$) such that $P$-a.s.
\begin{equation*}
\sup_{s\leq u\leq T}\Vert D\phi _{s,u}^{n}(x)-D\phi _{s,u}(x)\Vert
^{2}\rightarrow 0,\;\;\mbox{as}\;\;n\rightarrow \infty .
\end{equation*}
We find that $J\phi _{s,t}(x)=1$, for any $x\in {\mathbb{R}}^{d}$, and so $%
\phi _{s,t}$ is a measure-preserving diffeomorphism $P$-a.s. for any $%
s<t\leq T$.
\end{proof}

\section{Estimates on the derivative of the Jacobian}

\label{sec:jacobian} The aim of this section is to prove Sobolev
type
estimates on the derivative of the Jacobian of the stochastic flow 
$\phi _{t}=\phi _{0,t}$ associated to equation (\ref{SDE}) starting
at $0$. These estimates will be crucial in the proof of uniqueness
of weak solutions of the SPDE for $d\geq 2$
 (the case $d=1$ will be treated differently).

The basic observation is that the (formal) expression
\begin{equation}
\log J\phi _{t}(x)=\int_{0}^{t}\mathrm{div}\,b(s,\phi _{s}(x))ds
\label{eq:logJ}
\end{equation}
shows the opportunity of exploiting the It\^{o}--Tanaka trick to
regularize
the integrated divergence of $b$. We make the following hypothesis on $%
\mathrm{div}\,b$:

\begin{hypothesis}
\label{hy2} There exists $p \in (2,+\infty )$, such that
\begin{equation}
{\mathrm{div}\,}b\in L^{p} ([0,T]\times \mathbb{R}^{d}) \label{uni}
\end{equation}
(where ${\mathrm{div}\,}b(t,\cdot )$ is understood in distributional
sense).
\end{hypothesis}

 To apply the It\^{o}--Tanaka trick to eq.~\eqref{eq:logJ} the
 relevant PDEs  results are classical   $L^p$-parabolic
 estimates
 (see, for instance, \cite{Kr2})
  which are based on the following function spaces. For
$p \in (1,+\infty )$, we consider the Banach space $H_{p}^{2}(T)$ of
all functions $u\in
 L^{p}(0,T;W^{2,p}(\mathbb{R}%
^{d}))$ such that the distributional derivative $\partial _{t}u\in
L^{p} ([0,T]\times \mathbb{R}^{d})$. The norm is given by
\begin{equation*}
\Vert u\Vert _{{H_{p}^{2}(T)}} \Vert u\Vert _{L^{p}(0,T;W^{2,p}(\mathbb{R}%
^{d}))}+\Vert \partial _{t}u\Vert
_{L^{p}(0,T;L^{p}(\mathbb{R}^{d}))},
\end{equation*}
$u\in {H_{p}^{2}(T)}$.
 The next result is well known (see, for instance,  \cite[Theorem~9 in Section~7.3]{Kr2}).


\begin{theorem}
\label{KR} Consider a Borel and bounded function $l:[0,T]\times
\mathbb{R}^{d}\rightarrow \mathbb{R}^{d}$. For any function $f\in
L^{p} ([0,T]\times \mathbb{R}^{d})$, $p>1$, the Cauchy parabolic
problem
\begin{equation}
\left\{ \begin{aligned} \frac{\partial F}{\partial t}+ \frac{1}{2}
\Delta F + l \cdot D F = f,\;\; \; t \in [0,T[ \\ F(T,x)=0,\;\;\; x
\in \mathbb{R} ^d
\end{aligned}\right.   \label{eq:cauchy-parabolic-problem}
\end{equation}
has a unique solution $F$ in the space $H_{p}^{2}(T)$ and $F\in C( [  0,T]
;W^{1,p}(  \mathbb{R}^{d}) )  $. Moreover,
there exists a positive constant $C=C(p,
 d,T,\Vert l\Vert _{\infty
})$ such that
\begin{equation}
\Vert F\Vert _{{H_{p}^{2}(T)}}
\leq C\Vert
 f\Vert _{L^{p} ([0,T]\times \mathbb{R}^{d})}.  \label{bond}
\end{equation}
If $p\ge 2$ the constant $C$ above can be chosen such that
\begin{equation}
\sup_{t\in\left[
0,T\right]  }\left\|  F(  t,\cdot)  \right\|  _{W^{1,p}\left(
\mathbb{R}^{d}\right)  }
\leq C\Vert
 f\Vert _{L^{p} ([0,T]\times \mathbb{R}^{d})}.  \label{bond-sup}
\end{equation}
\end{theorem}
We prove now a  regularity result for the Jacobian $J \phi$.

\begin{theorem}
\label{iac} Under Hypotheses \ref{hy1} and \ref{hy2}  we
 have
 $J\phi \in $ $
 L^{2}(0,T;W^{1,2}_r)$ P-a.s. for any $r>0$.
\end{theorem}
\begin{proof}
\textbf{Step 1.}  
Recall
the chain rule for Sobolev function: if $f:\mathbb{R}^{d}\rightarrow \mathbb{%
R}$ is a continuous function, of class $W_{loc}^{1,2}(\mathbb{R}^{d})$ and $%
g:\mathbb{R}\rightarrow \mathbb{R}$ is a $C^{\infty }$ function, then $%
g\circ f\in W_{loc}^{1,2}(\mathbb{R}^{d})$ and
\begin{equation*}
 \int_{B(r)}\left| D(g\circ f)(x)\right| ^{2}dx \leq \left( \sup_{x\in B(r)}\left|
g^{\prime }(f(x))\right| \right) ^{2}\int_{B(r)}\left| Df(x)\right|
^{2}dx
\end{equation*}
for every $r>0$. Since $\log J\phi _{t}(x)$ is a  continuous
function, by the previous argument, if we
prove that $\log J\phi _{t}(x)\in W_{loc}^{1,2}(\mathbb{R}^{d})$ for a.e. $%
(\omega ,t)\in \Omega \times \lbrack 0,T]$, then $J\phi _{t}(x)\in
W_{loc}^{1,2}(\mathbb{R}^{d})$ for a.e. $(\omega ,t)\in \Omega
\times \lbrack 0,T]$ and
\begin{equation*}
\int_{B(r)}\left| DJ\phi _{t}(x)\right| ^{2}dx\leq \left( \sup_{x\in
B(r)}\left| J\phi _{t}(x)\right| \right) ^{2}\int_{B(r)}\left| D\log
J\phi _{t}(x)\right| ^{2}dx.
\end{equation*}
Integrating in $t\in \lbrack 0,T]$ and recalling that $\sup_{x\in
B(r),t\in
\lbrack 0,T]}\left| J\phi _{t}(x)\right| $ is finite ($P$-a.s.) because $%
J\phi _{t}(x)$ is continuous in $(t,x)$, we see that in order to
prove the theorem it is sufficient to prove that
\begin{equation*}
\log J\phi _{\cdot }(\cdot )\in L^{2}(0,T,W_{r}^{1,2})
\end{equation*}
for every $r>0$, for a.e. $\omega \in \Omega $. This will be proved
by showing that
\begin{equation} \label{pp}
\log J\phi _{\cdot }(\cdot )\in L^{2}(\Omega \times
(0,T),W_{r}^{1,2}).
\end{equation}
\textbf{Step 2.} Introduce $b^{\varepsilon }(t,x)=(\vartheta
_{\varepsilon }\ast b(t,\cdot ))(x)$, $\varepsilon >0$  and also set
$ b^{0}=b$.
Let $\phi _{t}^{\varepsilon }$ be the flow corresponding to the SDE %
\eqref{SDE} with $b$ replaced by $b^{\varepsilon }$. By well known
results (see \cite{K}), we get, for any $\varepsilon >0$,
\begin{equation*}
\log J\phi _{t}^{\varepsilon
}(x)=\int_{0}^{t}\mathrm{div}\,b^{\varepsilon }(s,\phi
_{s}^{\varepsilon }(x))ds.
\end{equation*}
Since the noise is additive, this can be also proved in an
elementary way by the $\omega $-wise application of the classical
deterministic results to the equation
\begin{equation*}
\rho _{t}^{\varepsilon }(x)=x+\int_{0}^{t}g^{\varepsilon }\left(
s,\rho _{s}^{\varepsilon }(x)\right) ds
\end{equation*}
where $\rho _{t}^{\varepsilon }(x)=\phi _{t}^{\varepsilon }(x)-W_{t}$, $%
g^{\varepsilon }\left( t,y\right) =b^{\varepsilon }(t,y+W_{t})$ (one has $%
J\phi _{t}^{\varepsilon }(x)=J\rho_{t}^{\varepsilon }(x)$ and $\mathrm{div}%
\,g^{\varepsilon }\left( t,y\right) =\mathrm{div}\,b^{\varepsilon
}(s,y+W_{t})$).

Note that, by  Remark~\ref{civuole},
$
J\phi _{t}^{\varepsilon }(x) \to  J\phi _{t}^{}(x)
$
in $L^{2}(\Omega \times (0,T),L_{r}^{2})$ as $\eps \to 0^+$. Define
$$
\psi_{\eps}(t,x)= \int_{0}^{t} \mathrm{div}\,b^{\varepsilon }\left(
s,{\phi }_{s}^{\varepsilon }\left( x\right) \right) ds.
$$
 Possibly passing to a sequence $(\eps_n)_{n\ge 1}$, we have
 that
\begin{equation}\label{fr6}
\psi_{{\eps_n}}(t,x) \to \log J\phi _{t}^{}(x)
 \end{equation} a.e. in $t, x,
\omega$, as $n \to \infty$.  On the other hand,  we have that
 $\left( \psi _{\varepsilon
}\right) _{\varepsilon >0}$ is bounded in $L^{2}(\Omega \times
(0,T);L_{r}^{2})$.
 Indeed, we have, using that $p>2$ and \eqref{bound-inv},
$$
\Big ( E \int_{0}^{T}\int_{B(r)} \left|
 \int_{0}^{t}  \mathrm{div}\,b^{\varepsilon }\left( s,{\phi
}_{s}^{\varepsilon }\left( x\right) \right)ds
   \right| ^{2} dxdt \Big )^{p/2}
$$
$$
\le C_{r, T}  E \int_{0}^{T}\int_{B(r)} \left|
 \int_{0}^{t}  \mathrm{div}\,b^{\varepsilon }\left( s,{\phi
}_{s}^{\varepsilon }\left( x\right) \right)ds
   \right| ^{p} dxdt
$$
$$
\le C_{r, T}'   E \int_{0}^{T} ds\int_{B(r)} \left|
   \mathrm{div}\,b^{\varepsilon }\left( s,{\phi
}_{s}^{\varepsilon }\left( x\right) \right)
   \right| ^{p} dx  \, \le
   C_{r, T}'   \int_{0}^{T} ds\int_{\RR^d} |
   \mathrm{div}\,b^{\varepsilon }\left( s,y
   \right)|^p \, E[ J (\phi_s^{\eps})^{-1}\, (y)]
   \, dy
$$
$$
\le
   C_{r, T}''   \sup_{s \in [0,T] , \, y \in \RR^d}
   E[ J (\phi_s^{\eps})^{-1}\, (y)]
    \int_{0}^{T} ds\int_{\RR^d} |
   \mathrm{div}\,b^{\varepsilon }\left( s,y
   \right)|^p
   \, dy \le C  < \infty.
$$
where $C$ is independent on $\eps >0$. Note that the previous
computation also  shows that  $\psi_{\eps}(t,x)  $ is uniformly
integrable on
  $\Omega \times [0,T] \times B(r)$.

By weak convergence we known that there exists a subsequence
 of $\{\psi_{\eps_n}\}_{n\ge 1}$ (still denoted by $\{\psi_{\eps_n}\}_{n\ge 1}$)
  which  converges weakly in $L^{2}(\Omega \times
(0,T);L_{r}^{2})$ to some function $\psi $.

On the other hand,  almost sure convergence and uniform
integrability together imply that $\psi_{\eps_n}$ converges strongly
 in $L^1 ( \Omega
 \times [0,T] \times B(r))$ to $\log J \phi $. It follows that, for any
 $\eta \in L^{\infty}
( \Omega
 \times [0,T] \times B(r))$, we have
$$
E \int_{0}^{T}\int_{B(r)} \psi_{}(t,x) \eta(t,x) dt dx = E
\int_{0}^{T}\int_{B(r)}  \log J \phi_t(x) \eta(t,x) dt dx.
$$
giving $\psi = \log J \phi $ which means that $\psi_{\eps_n}$
converges weakly in $L^{2}(\Omega \times (0,T);L_{r}^{2})$ to $ \log
J\phi $.

\medskip \textbf{Step 3.} To prove  assertion
\eqref{pp} it is enough to check that the family $\left( \psi
_{\varepsilon }\right) _{\varepsilon >0}$ is bounded in
$L^{2}(\Omega \times (0,T);W_{r}^{1,2})$.
Indeed, once we have proved this fact, we can extract from
 the previous
sequence
 $\psi _{\varepsilon _{n}}$ a subsequence which converges
  weakly in $L^{2}(\Omega
\times (0,T);W_{r}^{1,2})$ to some $\gamma $. This in particular
implies that such subsequence  converges weakly in $L^{2}(\Omega
\times (0,T),L_{r}^{2})$ to  $\gamma $. By the previous step, we
must have that $\gamma = J \phi$.

We introduce the following Cauchy problem, for $\eps \ge 0$,
\begin{equation}
\left\{ \begin{aligned} \frac{\partial F^{\eps}}{\partial t}+
\frac{1}{2} \Delta F^{\eps} + D F^{\eps} \cdot b^{\eps}
={\mathrm{div}\, }b^{\eps}, \;\;
\; t \in [0,T[ \\ F^{\eps}(T,x)=0,\;\;\; x \in \mathbb{R} ^d. \end{aligned}%
\right.   \label{equa}
\end{equation}
Note that by Theorem \ref{KR} and since $p> 2$ we have
\begin{equation}
\Vert F^{\eps}\Vert _{{H_{p}^{2}(T)}}
+\sup_{t\in\left[
0,T\right]  }\left\|  F^\eps\left(  t,\cdot\right)  \right\|  _{W^{1,p}(
\mathbb{R}^{d})  }\leq C\Vert
 \mathrm{div}\, b \Vert _{L^{p}(0,T;L^{p}(\mathbb{
R}^{d}))},  \label{bond4}
\end{equation}
for any $\eps \ge 0$. Using It\^{o} formula we find (remark that
 $F^{\eps}(t, \cdot) \in C^{2}_b(\RR^d)$)
\begin{equation*}
F^{\varepsilon }\left( t,{\phi }_{t}^{\varepsilon }\left( x\right)
\right)
-F^{\varepsilon }\left( 0,x\right) -\int_{0}^{t}DF^{\varepsilon }\left( s,{%
\phi }_{s}^{\varepsilon }\left( x\right) \right) \cdot dW_{s}\int_{0}^{t}%
\mathrm{div}\,b^{\varepsilon }\left( s,{\phi }_{s}^{\varepsilon
}\left( x\right) \right) ds = \psi _{\varepsilon }(t,x ).
\end{equation*}
 By classical results (see \cite{K}) the maps $%
\psi _{\varepsilon }(t,\cdot )$ are differentiable, $P$-a.s. and
\begin{align}\label{ci5}
 D\psi_{\varepsilon }(t,x) D F^{\varepsilon } \left( t,{\phi
}_{t}^{\varepsilon }\left( x\right) \right) D{\phi
}_{t}^{\varepsilon }\left( x\right)
 - DF^{\varepsilon
}\left( 0,x\right) \end{align} $$
 - \int_{0}^{t}D^{2}F^{\varepsilon }(s,{\phi }
 _{s}^{\eps}(x)) D{\phi
}_{s}^{\varepsilon }\left( x\right) dW_{s}
$$
Since we already know that $\left( \psi _{\varepsilon }\right)
_{\varepsilon >0}$ is bounded in $L^{2}(\Omega \times
(0,T),L_{r}^{2})$, to verify that $\left( \psi _{\varepsilon
}\right) _{\varepsilon >0}$ is bounded in $L^{2}(\Omega \times
(0,T);W_{r}^{1,2})$, it is enough to prove that  $\left( D\psi
_{\varepsilon }\right) _{\varepsilon >0}$ is bounded in
$L^{2}(\Omega \times (0,T),L_{r}^{2})$.


To this purpose we only prove the bound for  the critical term
$\int_{0}^{t}D^{2}F^{\varepsilon }(s,{\phi }
 _{s}^{\eps}(x)) D{\phi
}_{s}^{\varepsilon }\left( x\right) dW_{s}$ in \eqref{ci5}.
 The other terms are easier to estimate, we remark only that for the term $DF^{\varepsilon}(  0,x)  $ we use eq.~\eqref{bond4}.

We show  that  there exists  a constant $C>0$ (depending on $r$
 and $T$)
 such that
\begin{equation}
E\left[ \int_{0}^{T}\int_{B(r)}\left|
\int_{0}^{t}D^{2}F^{\varepsilon }\left( s,{\phi }_{s}^{\eps}\left(
x\right) \right) D{\phi }_{s}^{\eps}\left( x\right) dW_{s}\right|
^{2}dxdt\right] \leq C  \label{cia}
\end{equation}
for every $\varepsilon >0$. Note that
\begin{align*}
& \int_{0}^{T}\int_{B(r)}E\left[ \left|
\int_{0}^{t}D^{2}F^{\varepsilon }\left( s,{\phi }_{s}^{\eps}\left(
x\right) \right) D{\phi }_{s}^{\eps}\left( x\right)
dW_{s}\right| ^{2}\right] dxdt \\
&
\leq T\int_{B(r)}E\left[ \int_{0}^{T}\left| D^{2}F^{\varepsilon }\left( s,{%
\phi }_{s}^{\eps}
 \left( x\right) \right) D{\phi }_{s}^{\eps}\left( x\right) \right| ^{2}ds%
\right] dx.
\end{align*}
According to \eqref{bound}, we have
\begin{equation*}
\int_{B(r)}\left( \int_{0}^{T}E\left[ \left| D{\phi
}_{s}^{\eps}\left( x\right) \right| ^{r}\right] ds\right) dx \le C <
\infty ,
\end{equation*}
 for every $r\geq 1$, with $C$ independent on
 $\eps$;
 therefore  by the H\"{o}lder inequality on $%
\Omega \times B(r)\times \lbrack 0,T]$, it is sufficient to prove
that there exists $C >0$ such that, for any $\eps>0$,
\begin{equation} \label{cc}
\int_{B(r)}E\left[ \int_{0}^{T}\left| D^{2}F^{\varepsilon }\left( s,{%
\phi }_{s}^{\eps}\left( x\right) \right) \right| ^{p^{ }}ds  \right]
dx \le C <\infty .
\end{equation}

\textbf{Step 4}. Let us show \eqref{cc}. We have
$$E\left[  \int_{0}^{T} ds\int_{B(r)}\left| D^{2}F^{\varepsilon }\left( s,%
{\phi }_{s}^{\eps}\left( x\right) \right) \right| ^{p^{ }}dx \right]
$$
$$ \le   \int_{0}^{T} ds\int_{\RR^d}\left|
D^{2}F^{\varepsilon }\left( s,{y}\right) \right| ^{p^{ }} E [J(\phi
_{s}^{\eps})^{-1}(y)] \, dy
$$ $$\leq \sup_{s \in [0,T] , \, y \in \RR^d}
   E[ J (\phi_s^{\eps})^{-1}\, (y)]
\, \int_{0}^{T} ds\int_{\RR^d}\left| D^{2}F^{\varepsilon }\left(
s,{y}\right) \right| ^{p^{ }} \,  dy \le C < \infty,
$$
where, using  \eqref{bound-inv} and  \eqref{bond4}, $C$ is
independent on $\eps>0.$ This proves \eqref{cc} and ends the proof.
\end{proof}

\section{Stochastic transport equation. Existence of weak solutions}

\label{sec:transp}

To avoid doubts, let us clarify a convention of language we use in the
sequel. An element $u\in L^{\infty}( \Omega\times[ 0,T] \times\mathbb{R}%
^{d}) $ is an equivalence class. Given $\theta\in L^{1}( \mathbb{R}^{d}) $,
the notation $\int_{\mathbb{R}^{d}}\theta( x) u( t,x) dx$ stands for an
element of $L^{\infty}( \Omega\times[ 0,T] ) $ (so again an equivalence
class) defined by Fubini theorem. When we say that $\int _{\mathbb{R}%
^{d}}\theta( x) u( t,x) dx$ \textit{has a continuous modification} we mean
that there exists an element in the equivalence class that is a continuous
stochastic process (a process with continuous paths, $P$-a.s.). We choose
this language so that it is the same as in the case when $u$ is a measurable
function $u:\Omega\times[ 0,T] \times\mathbb{R}^{d}\rightarrow\mathbb{R}$
instead of an equivalence class (moreover, this way, when the property is
true for the equivalent class it is true for all its representatives).

We still use the name `stochastic process' for the elements of $L^{\infty }(
\Omega\times[ 0,T] \times\mathbb{R}^{d}) $.

In the following definition we perform a Stratonovich integration. This is
well defined when the integrator is a continuous semimartingale adapted to
the filtration of the Brownian motion, see~\cite{K}.

\begin{definition}
\label{def solution} Let $b\in L_{loc}^{1}([0,T]\times \mathbb{R}^{d};%
\mathbb{R}^{d})$, $\mathrm{div}\,b\in L_{loc}^{1}([0,T]\times \mathbb{R}^{d})
$ and $u_{0}\in L^{\infty }(\mathbb{R}^{d})$. A weak $L^{\infty }$-solution
of the Cauchy problem (1) is a stochastic process $u\in L^{\infty }(\Omega
\times \lbrack 0,T]\times \mathbb{R}^{d})$ such that, for every test
function $\theta \in C_{0}^{\infty }(\mathbb{R}^{d})$, the process $\int_{%
\mathbb{R}^{d}}\theta (x)u(t,x)dx$ has a continuous modification which is an
$\mathcal{F}$-semimartingale and
\begin{align*}
\int_{\mathbb{R}^{d}}u(t,x)\theta (x)dx& =\int_{\mathbb{R}%
^{d}}u_{0}(x)\theta (x)dx \\
& +\int_{0}^{t}ds\int_{\mathbb{R}^{d}}u(s,x)\left[ b(s,x)\cdot D\theta (x)+%
\mathrm{div}\,b(s,x)\theta (x)\right] dx \\
& +\sum_{i=1}^{d}\int_{0}^{t}\left( \int_{\mathbb{R}^{d}}u(s,x)D_{i}\theta
(x) dx\right) \circ dW_{s}^{i}.
\end{align*}
\end{definition}

In the previous definition we have used Stratonovich integrals since they
are the natural ones in this framework. However, as usual, one can
reformulate the problem in It\^{o} form and avoid the
semimartingale assumption.

\begin{lemma}
\label{lemma su Ito-Stratonovich}A process $u\in L^{\infty }(\Omega \times
\lbrack 0,T]\times \mathbb{R}^{d})$ is a weak $L^{\infty }$-solution of the
Cauchy problem~\eqref{SPDE} if and only if, for every test function $\theta
\in C_{0}^{\infty }(\mathbb{R}^{d})$, the process $\int_{\mathbb{R}%
^{d}}\theta (x)u(t,x)dx$ has a continuous $\mathcal{F}$-adapted
modification and\
\begin{align*}
\int_{\mathbb{R}^{d}}u(t,x)\theta (x)dx& =\int_{\mathbb{R}%
^{d}}u_{0}(x)\theta (x)dx \\
& +\int_{0}^{t}ds\int_{\mathbb{R}^{d}}u(s,x)\left[ b(s,x)\cdot D\theta (x)+%
\mathrm{div}\,b(s,x)\theta (x)\right] dx
\end{align*}
\begin{equation}
+\sum_{i=1}^{d}\int_{0}^{t}\left( \int_{\mathbb{R}^{d}}u(s,x)D_{i}\theta
(x)dx\right) dW_{s}^{i}+\frac{1}{2}\int_{0}^{t}ds\int_{\mathbb{R}%
^{d}}u(s,x)\Delta \theta (x)dx  \label{Ito form of the SPDE}
\end{equation}
for a.e. $(\omega ,t)\in \Omega \times \lbrack 0,T]$.
\end{lemma}

\begin{proof}
The relation between It\^{o} and Stratonovich integrals is (see \cite{K})
\begin{eqnarray*}
&&\int_{0}^{t} \left( \int_{\mathbb{R}^{d}}u(s,x)D_{i}\theta(x) dx\right) \circ
dW_{s}^{i} \\ && \qquad = \int_{0}^{t}\left( \int_{\mathbb{R}^{d}}u(s,x)D_{i}\theta
(x)dx\right) dW_{s}^{i}
+\frac{1}{2}\left[ \int_{\mathbb{R}^{d}}u(\cdot ,x)D_{i}\theta(x) dx,W_{\cdot
}^{i}\right] _{t}
\end{eqnarray*}
where $\left[ \cdot ,\cdot \right] _{t}$ denotes the joint quadratic
variation. Only the martingale part of $\int_{\mathbb{R}^{d}}u(\cdot
,x)D_{i}\theta(x) dx$ counts in the joint quadratic variation. If we start from
definition \ref{def solution}, the martingale part of $\int_{\mathbb{R}%
^{d}}u(\cdot ,x)D_{i}\theta(x) dx$ is (taking $D_{i}\theta $ in place of $%
\theta $ in the equation)
\begin{equation*}
\sum_{j=1}^{d}\int_{0}^{t}\left( \int_{\mathbb{R}^{d}}u(s,x)D^2_{ij}\theta(x)
dx\right) dW_{s}^{j}.
\end{equation*}
The same is true if, conversely, we start from equation (\ref{Ito form of
the SPDE}). The joint quadratic variation is therefore equal to (see \cite{K}%
)
\begin{equation*}
\int_{0}^{t}\left( \int_{\mathbb{R}^{d}}u(s,x)D^2_{ij}\theta(x) dx\right)
\cdot 1\, ds.
\end{equation*}
Summing over $i$, we get the result. The other details of the equivalence
statement are easy. The proof is complete.
\end{proof}

\begin{remark}
The presence of the Laplacian in the It\^{o} formulation should not suggest
that our SPDE has of parabolic nature. As clarified in the classical
literature on parabolic SPDEs, an SPDE in It\^{o} form is parabolic when a
super-parabolicity condition holds. If the second order operator in the
drift has the form $\sum_{i,j=1}^{d}a_{ij}D_{j}D_{i}u$ and the first order
part of the multiplicative noise has the form $\sum_{i,j=1}^{d}\sigma
_{ij}D_{j}udW_{t}^{i}$, then the operator
\begin{equation*}
\sum_{i,j=1}^{d}\left( a_{ij}-\frac{1}{2}\left( \sigma ^{T}\sigma \right)
_{ij}\right) D_{j}D_{i}u
\end{equation*}
must be strongly elliptic (see for instance \cite{Roz}, \cite{DaPZab}). In
our case this operator is equal to zero.
\end{remark}

We may now prove a very general existence result, similarly to the
deterministic case. The proof is essentially the same of that for SPDEs with
monotone operators, see~\cite{KrylRoz,Pard,Roz}.

\begin{theorem}
Let $b\in L_{\mathrm{loc}}^{1}([0,T]\times \mathbb{R}^{d};\mathbb{R}^{d})$, $%
\mathrm{div}\,b\in L_{\mathrm{loc}}^{1}([0,T]\times \mathbb{R}^{d})$ and $%
u_{0}\in L^{\infty }(\mathbb{R}^{d})$. Then there exists a weak $L^{\infty }$%
-solution of the SPDE \eqref{SPDE}.
\end{theorem}

\begin{proof}
Let $\eta ^{\varepsilon }$ be a $C^{\infty }(\mathbb{R}\times \mathbb{R}^{d})
$ mollifier; let $\zeta \in C_{0}^{\infty }(\mathbb{R}^{d})$ be such that $%
\zeta (x)=1$ for $|x|\leq 1$ and $\zeta (x)=0$ for $|x|\geq 2$. Define $%
b^{\varepsilon }=\eta ^{\varepsilon }\ast \,(\zeta (\varepsilon \,\cdot )\,b)
$ (for $t\not\in \lbrack 0,T]$, we set $b_{t}=0$). Note that each $%
b^{\varepsilon }$ is globally Lipschitz. Let $\phi ^{\varepsilon }$ be the
associated flow and define $u_{t}^{\varepsilon }(x)=u_{0}((\phi
^{\varepsilon })_{t}^{-1}(x))$. It is known (it can be checked by direct
computation, see \cite{K1}) that this is the unique classical solutions of
the associated transport equation, that written in weak It\^{o} form is
equation (\ref{Ito form of the SPDE}) with $u^{\varepsilon }$ and $%
b^{\varepsilon }$ in place of $u$ and $b$. From the representation in terms
of the flow we immediately have $\sup_{x,\omega ,t}|u_{t}^{\varepsilon
}(x)(\omega )|\leq C$ uniformly in $\varepsilon $ so there exists a sequence
$u^{\varepsilon _{n}}$ converging weak-$*$ in $L^{\infty }(\Omega \times
\lbrack 0,T]\times \mathbb{R}^{d})$ and weakly in $L^{2}(\Omega \times
\lbrack 0,T]\times B\left( N\right) )$ for every integer $N>0$, to some $u$,
which belongs to these spaces. For shortness, let us denote: $\varepsilon
_{n}$ by $\varepsilon $; for every $\theta \in L^{1}(\mathbb{R}^{d})$, $%
\int_{\mathbb{R}^{d}}\theta (x)u^{\varepsilon }(t,x)dx$ by $%
u_{t}^{\varepsilon }(\theta )$, including $u_{t}(\theta )$ as the case $%
\varepsilon =0$; $b^{\varepsilon }(s,x)$ by $b_{s}^{\varepsilon }(x)$, again
also for $\varepsilon =0$.

We follow here the arguments of \cite{Pard}, Chapter III. Let $\theta \in
C_{0}^{\infty }(\mathbb{R}^{d})$. The process $u^{\varepsilon }(\theta )$ is
non anticipative (by its definition) and converges weakly in $L^{2}(\Omega
\times \lbrack 0,T])$ (by the weak convergence of $u^{\varepsilon }$ to $u$
in $L^{2}(\Omega \times \lbrack 0,T]\times B\left( N\right) )$ for $N$ such
that $B\left( N\right) $ contains the support of $\theta $). The space of
non anticipative processes is a closed subspace of $L^{2}(\Omega \times
\lbrack 0,T])$, hence weakly closed. Therefore $u(\theta )$ is non
anticipative. Thus It\^{o} integrals of $u(\theta )$ (which is also bounded)
are well defined. Moreover, the mapping $h\mapsto \int_{0}^{\cdot }h\cdot dW$
is linear continuous from the space of non-anticipative $L^{2}(\Omega \times
\lbrack 0,T];\mathbb{R}^{d})$-processes to $L^{2}(\Omega \times \lbrack 0,T])
$;\ then it is also weakly continuous. Therefore $\int_{0}^{\cdot
}u_{s}^{\varepsilon }\left( D\theta \right) \cdot dW_{s}$ converges weakly
in $L^{2}(\Omega \times \lbrack 0,T])$ to $\int_{0}^{\cdot }u_{s}\left(
D\theta \right) \cdot dW_{s}$.

Moreover we have $b^{\varepsilon }\rightarrow b$ and $\mathrm{div}%
\,b^{\varepsilon }\rightarrow \mathrm{div}\,b$ in $L_{\mathrm{loc}%
}^{1}([0,T]\times \mathbb{R}^{d})$. For every $\varepsilon \geq 0$, define
\begin{equation*}
G_{t}^{\varepsilon }(u^{\varepsilon },\theta )=u_{t}^{\varepsilon }(\theta
)-u_{0}\left( \theta \right) -\int_{0}^{t}u_{s}^{\varepsilon }\left(
b_{s}^{\varepsilon }\cdot D\theta +\theta \mathrm{div}\,b_{s}^{\varepsilon
}\right) ds-\frac{1}{2}\int_{0}^{t}u_{s}^{\varepsilon }\left( \Delta \theta
\right) ds.
\end{equation*}
It is not difficult to check that $G^{\varepsilon }(u^{\varepsilon },\theta )
$ converges weakly to $G(u,\theta )$ in $L^{2}(\Omega \times \lbrack 0,T])$.
To this purpose notice that $b^{\varepsilon }\cdot D\theta +\theta \mathrm{%
div}\,b^{\varepsilon }$ converges strongly to $b\cdot D\theta +\theta
\mathrm{div}\,b$ in $L^{1}([0,T]\times \mathbb{R}^{d})$.

Therefore we may pass to the weak $L^{2}(\Omega \times \lbrack 0,T])$-limit
in the equation for $u^{\varepsilon }$ and prove that $u$ satisfies equation
(\ref{Ito form of the SPDE}) for a.e. $\left( \omega ,t\right) \in \Omega
\times \left[ 0,T\right] $.

Finally, the right-hand-side of equation (\ref{Ito form of the SPDE})
defines a continuous stochastic process. Therefore $u(\theta )$ has a
continuous modification. The proof is complete.
\end{proof}

Under more restrictive conditions we may construct a solution related to the
stochastic flow.

\begin{theorem}
\label{them SPDE n1} Assume that hypothesis \ref{hy1} holds and $\mathrm{div}%
\,b\in L_{\mathrm{loc}}^{1}([0,T]\times \mathbb{R}^{d})$. Given $u_{0}\in
L^{\infty }({}{\mathbb{R}}^{d})$, the stochastic process $u(t,x)$ defined as
$u(t,x)=u_{0}(\phi _{t}^{-1}(x))$ is a solution of (\ref{SPDE}).
\end{theorem}

\begin{proof}
\textbf{Step 1}. We first prove the claim when $u_{0}$ has support in some
ball $B\left( R\right) $. Let $b^{\varepsilon }$ be a regularization of $b$
as described in Section~\ref{sect notations}. It converges to $b$ in $%
L^{\infty }\left( 0,T;C_{b}^{\alpha ^{\prime }}({\mathbb{R}}^{d};{\mathbb{R}}%
^{d})\right) $ for every $\alpha ^{\prime }<\alpha $. We apply the argument
of the previous proof to this particular $b^{\varepsilon }$ and the
associated solution $u_{t}^{\varepsilon }(x)=u_{0}((\phi ^{\varepsilon
})_{t}^{-1}(x))$. Let $u$ be one of its weak (and weak-$*$) limits, as
described by the previous proof, and $u^{\varepsilon _{n}}$ the
corresponding sequence. We know that $u$ is a solution of (\ref{SPDE}). We
shall use that $u^{\varepsilon _{n}}$ converges weakly to $u$ in $%
L^{2}(\Omega \times \lbrack 0,T]\times B\left( N\right) )$ for every integer
$N>0$. Thus, for every $\theta \in C_{0}^{\infty }(\mathbb{R}^{d})$, $%
u^{\varepsilon _{n}}\left( \theta \right) $ (we use the notations of the
previous proof) converges to $u\left( \theta \right) $ weakly in $%
L^{2}(\Omega \times \lbrack 0,T])$.

Let us prove that a further subsequence, still denoted by $%
u_{t}^{\varepsilon _{n}}\left( \theta \right) $, converges to the quantity $\int_{{%
\mathbb{R}}^{d}}u_{0}(\phi _{t}^{-1}(x))\theta \left( x\right) dx$ for a.e. $%
\left( \omega ,t\right) \in \Omega \times \lbrack 0,T]$. Since $%
u_{t}^{\varepsilon _{n}}\left( \theta \right) $ is equibounded, the
convergence is also strong in $L^{2}(\Omega \times \lbrack 0,T])$ and thus
weak. Therefore $u_t\left( \theta \right) =\int_{{\mathbb{R}}^{d}}u_{0}(\phi
_{t}^{-1}(x))\theta \left( x\right) dx$. This implies the claim of the
theorem. Let $t\in \left[ 0,T\right] $ be given. We have
\begin{equation*}
u_{t}^{\varepsilon _{n}}\left( \theta \right) =\int_{B\left( R\right)
}u_{0}(y)\theta \left( \phi _{t}^{\varepsilon _{n}}\left( y\right) \right)
J\phi _{t}^{\varepsilon _{n}}(y)dy.
\end{equation*}
By Remark~\ref{civuole}, up to a subsequence still denoted by $\varepsilon
_{n}$, $\theta \left( \phi _{t}^{\varepsilon _{n}}\left( \cdot \right)
\right) J\phi _{t}^{\varepsilon _{n}}(\cdot )$ converges in $L^{1}\left(
B\left( R\right) \right) $ to $\theta \left( \phi _{t}\left( \cdot \right)
\right) J\phi _{t}(\cdot )$, $P$-a.s. and thus, $P$-a.s., $%
u_{t}^{\varepsilon _{n}}\left( \theta \right) $ converges to $\int_{B\left(
R\right) }u_{0}(y)\theta \left( \phi _{t}\left( y\right) \right) J\phi
_{t}(y)dy$. This is what we wanted to prove.

\textbf{Step 2}. Consider now a general $u_{0}\in L^{\infty }({}{\mathbb{R}}%
^{d})$. Let $\zeta \in C_{0}^{\infty }(\mathbb{R}^{d})$ be such that $\zeta
(x)=1$ for $|x|\leq 1$, $\zeta (x)=0$ for $|x|\geq 2$, $\zeta (x)\in \left[
0,1\right] $ for all $x\in \mathbb{R}^{d}$. Define $u_{0}^{n}\left( x\right)
=u_{0}\left( x\right) \zeta (n^{-1}x)$. Let $u^{n}$ be the corresponding
solution given by step 1. We have
\begin{equation*}
u^{n}(t,x)=u_{0}(\phi _{t}^{-1}(x))\zeta (n^{-1}\phi _{t}^{-1}(x)).
\end{equation*}
The function $u^{n}$ converges pointwise to $u(t,x)=u_{0}(\phi _{t}^{-1}(x))
$ and, because of equiboundedness, strongly in $L^{2}(\Omega \times \lbrack
0,T]\times B\left( N\right) )$ for every integer $N>0$, and weak-$*$ in $%
L^{\infty }(\Omega \times \lbrack 0,T]\times \mathbb{R}^{d})$. Moreover, $%
u^{n}$ verifies equation (\ref{Ito form of the SPDE}). It is now easy to
repeat the proof of the previous theorem and check that $u$ is a solution.
The proof is complete.
\end{proof}

\section{Uniqueness of weak solutions}

\label{sec:uniq} In this section, we prove uniqueness for our SPDE.
Our main results are Theorems \ref{thm SPDE n2}
 and  \ref{thm SPDE n3}.

\medskip

Let us recall that in the deterministic case with non-regular vector fields $b$%
, uniqueness of weak $L^{\infty }$-solutions is proved by means of
the concept of renormalized solutions, see
\cite{DiPernaLions,Ambrosio}. The technical tool is a commutator
lemma, where the role of some of the
assumptions on $b$ is immediately visible. Given $g\in L_{loc}^{\infty }(%
\mathbb{R}^{d})$, $v\in L_{loc}^{1}(\mathbb{R}^{d},\mathbb{R}^{d})$ , $%
\mathrm{div}\,v\in L_{loc}^{1}(\mathbb{R}^{d})$, denote by $v\cdot
Dg$ the distribution defined on smooth compactly supported test
functions $\rho $ as
\begin{equation*}
(v\cdot Dg)(\rho )=-\int gv\cdot D\rho \,dx-\int g\rho
\,\mathrm{div}\,v\,dx.
\end{equation*}
Given mollifiers $\vartheta _{\varepsilon }$ as in Section~\ref{sect
notations}, define the \emph{commutator} $\mathcal{R}_{\varepsilon
}[v,g]$ as the smooth function
\begin{equation*}
\mathcal{R}_{\varepsilon }[v,g](x)=\left[ \vartheta _{\varepsilon
}\ast \left( v\cdot Dg\right) -v\cdot D\left( \vartheta
_{\varepsilon }\ast g\right) \right] \left( x\right) .
\end{equation*}
In the case $d=1$ we shall use the following classical form of the
commutator lemma.

\begin{lemma}
\label{lemma:commutator-deterministic} Given $v\in W_{loc}^{1,1}(\mathbb{R}%
^{d},\mathbb{R}^{d})$ , $g\in L_{loc}^{\infty }(\mathbb{R}^{d})$, we
have
\begin{equation*}
\int_{B\left( r\right) }|\mathcal{R}_{\varepsilon }[v,g](x)|dx\leq
C\Vert v\Vert _{W_{r+2}^{1,1}}\Vert g\Vert _{L_{r+2}^{\infty }}
\end{equation*}
for all $r>0$, for some constant $C>0$ independent of $\varepsilon
$, $v$, $g $ and $r$. Moreover, $\lim_{\varepsilon \rightarrow 0}$
$\int_{B\left( r\right) }|\mathcal{R}_{\varepsilon }[v,g](x)|dx=0$.
\end{lemma}

\begin{proof}
In a global form on $\mathbb{R}^{d}$ it is proved for instance in
\cite[Ch. 2]{Lions}, Lemma~2.3. The local form on balls $B\left(
r\right) $ can be easily deduced from the global form by multiplying
$v$ and $g$ by a smooth function which takes values
in $\left[ 0,1\right] $, and it
is equal to 1 on $B\left( r+1\right) $ and to 0 outside $B\left( r+2\right) $%
.
\end{proof}

\medskip The strong $L^{1}$ convergence of commutators requires some
weak form of differentiability of $v$ or $g$. We impose such
differentiability conditions
only in the case $d=1$ where ${\mathrm{div}\,}b=Db$, or in the case of $%
BV_{loc}$ solutions treated in Appendix \ref{sec:bv}. For $d>1$ and $%
L^{\infty }$ solutions, taking advantage of the presence of the
flow, we can prove uniqueness by means the distributional
convergence of commutators \emph{composed with the flow}. This
composition put into play the first derivatives of the Jacobian of
the flow.

The special estimates of Section~\ref{sec:jacobian} becomes the main
tool to prove our first uniqueness result (see Theorem \ref{thm SPDE
n2}).

We start by giving some preliminary easy estimates on the distributional
commutator and on its composition with a $C^{1}$ diffeomorphism.

\begin{lemma} \label{comm}
Given $v\in L_{loc}^{\infty }\left( \mathbb{R}^{d},\mathbb{R}^{d}\right) $, $%
{\mathrm{div}\,}v\in L_{loc}^{1}\left( \mathbb{R}^{d}\right) $,
$g\in L_{loc}^{\infty }\left( \mathbb{R}^{d}\right) $;
for any $\rho \in C_{r}^{\infty }(\mathbb{R}^{d})$ and for
sufficiently small $\varepsilon $ we have
\begin{equation*}
\left| \int \mathcal{R}_{\varepsilon }\left[ v,g\right] \left(
x\right) \rho
\left( x\right) dx\right|  \\
\leq 2\left\| g\right\| _{L_{r+1}^{\infty }}\left[ \left\| {\mathrm{div}}%
v\right\| _{L_{r+1}^{1}}\left\| \rho \right\| _{L_{r}^{\infty
}}+\left\| v\right\| _{L_{r+1}^{\infty }}\left\| D\rho \right\|
_{L_{r}^{1}}\right]
\end{equation*}
and
\begin{equation*}
\lim_{\varepsilon \rightarrow 0}\int \mathcal{R}_{\varepsilon }\left[ v,g%
\right] \left( x\right) \rho \left( x\right) dx=0.
\end{equation*}
\end{lemma}

\begin{proof}
The proof proceeds as the one of
Lemma~\ref{lemma:commutator-deterministic} from \cite[Ch. 2]{Lions}:
we prove the inequality for regular fields, then one can extend to
non-regular ones and prove the convergence again first in the
regular case and then apply an approximation procedure which we will
omit. Let us rewrite the expressions in a suitable way:
\begin{equation*}
\begin{split}
\int \mathcal{R}_{\varepsilon }[v,g](x)\rho (x)dx& =\int [\vartheta
_{\varepsilon }\ast \left( v\cdot Dg\right) ](x^{\prime })\rho
(x^{\prime })dx^{\prime }-\int [v\cdot D\left( \vartheta
_{\varepsilon }\ast g\right)
](x)\rho (x)dx \\
& =\iint [-g(x^{\prime })v(x^{\prime })D_{x^{\prime }}\vartheta
_{\varepsilon }(x-x^{\prime })-g(x^{\prime })\vartheta _{\varepsilon
}(x-x^{\prime })\mathrm{div}\,v(x^{\prime })]\rho (x)\,dxdx^{\prime } \\
& \quad +\iint [g(x^{\prime })v(x)\vartheta _{\varepsilon
}(x-x^{\prime })D_{x}\rho (x)+g(x^{\prime })\rho (x)\vartheta
_{\varepsilon }(x-x^{\prime })\mathrm{div}\,v(x)]\,dxdx^{\prime }.
\end{split}
\end{equation*}
Using $D_{x^{\prime }}\vartheta _{\varepsilon }(x-x^{\prime
})=-D_{x}\vartheta _{\varepsilon }(x-x^{\prime })$ and integrating
by parts in $x$ the first term we get
\begin{equation*}
\begin{split}
& =\iint g(x^{\prime })\vartheta _{\varepsilon }(x-x^{\prime
})D_{x}\rho
(x)[v(x)-v(x^{\prime })]\,dxdx^{\prime } \\
& \qquad +\iint g(x^{\prime })\rho (x)\vartheta _{\varepsilon
}(x-x^{\prime })[\mathrm{div}\,v(x)-\mathrm{div}\,v(x^{\prime
})]\,dxdx^{\prime }.
\end{split}
\end{equation*}
Assume $\varepsilon $ so small that the support of $\vartheta
_{\varepsilon } $ has diameter less than one. We have (using
standard estimates on convolutions)
\begin{equation*}
\left| \int \mathcal{R}_{\varepsilon }\left[ v,g\right] (x)\rho
(x)dx\right| \leq 2\left\| g\right\| _{L_{r+1}^{\infty }}\left\|
v\right\| _{L_{r+1}^{\infty }}\left\| D\rho \right\|
_{L_{r}^{1}}+2\left\| \rho \right\| _{L_{r}^{\infty }}\left\|
g\right\| _{L_{r+1}^{\infty }}\left\| \mathrm{div}\,v\right\|
_{L_{r+1}^{1}}.
\end{equation*}
\end{proof}

\begin{corollary}
\label{corollary commutator}Let $\phi $ be an $C^{1}$ diffeomorphism of $%
\mathbb{R}^{d}$. Assume $v\in L_{\mathrm{loc}}^{\infty }\left( \mathbb{R}%
^{d},\mathbb{R}^{d}\right) $, ${\mathrm{div}\,}v\in L_{\mathrm{loc}%
}^{1}\left( \mathbb{R}^{d}\right) $, $g\in L_{\mathrm{loc}}^{\infty
}\left( \mathbb{R}^{d}\right) $. Moreover, if $d>1$, assume also
$J\phi ^{-1}\in W_{loc}^{1,1}\left( \mathbb{R}^{d}\right) $. Then
for any $\rho \in C_{0}^{\infty }(\mathbb{R}^{d})$ there exists a
constant $C_{\rho }>0$ such that, given any $R>0$ such that
$\mathrm{supp}(\rho \circ \phi ^{-1})\subseteq B(R)$, we have:

\begin{itemize}
\item  for $d>1$,
\begin{equation*}
\begin{split}
& \left| \int \mathcal{R}_{\varepsilon }[v,g](\phi (x))\rho (x)dx\right|  \\
& \qquad \leq C_{\rho }\Vert g\Vert _{L_{R+1}^{\infty }}\left[ \Vert {\mathrm{div}}%
v\Vert _{L_{R+1}^{1}}\Vert J\phi ^{-1}\Vert _{L_{R}^{\infty }}+\Vert
v\Vert _{L_{R+1}^{\infty }}(\Vert D\phi ^{-1}\Vert _{L_{R}^{\infty
}}+\Vert DJ\phi ^{-1}\Vert _{L_{R}^{1}})\right]
\end{split}
\end{equation*}

\item  for $d=1$,
\begin{equation*}
\left| \int \mathcal{R}_{\varepsilon }[v,g](\phi (x))\rho
(x)dx\right| \leq C_{\rho }\Vert J\phi ^{-1}\Vert _{L_{R}^{\infty
}}\Vert v\Vert _{W_{R+2}^{1,1}}\Vert g\Vert _{L_{R+2}^{\infty }}.
\end{equation*}
In both cases we have
\begin{equation*}
\lim_{\varepsilon \rightarrow 0}\int \mathcal{R}_{\varepsilon }\left[ v,g%
\right] \left( \phi \left( x\right) \right) \rho \left( x\right)
dx=0.
\end{equation*}
\end{itemize}
\end{corollary}

\begin{proof}
By a change of variables we have $\int \mathcal{R}_{\varepsilon
}[v,g](\phi (x))\rho (x)dx=\int \mathcal{R}_{\varepsilon
}[v,g](y)\rho _{\phi }(y)dx$ where the function $\rho _{\phi
}(y)=\rho (\phi ^{-1}(y))J\phi ^{-1}(y)$ has the support strictly
contained in the ball of radius $R$.
For $d>1$, by the previous lemma
\begin{equation*}
\left| \int \mathcal{R}_{\varepsilon }\left[ v,g\right] \left( \phi
\left( x\right) \right) \rho \left( x\right) dx\right| \leq 2\left\|
g\right\| _{L_{R+1}^{\infty }}\left[ \left\| {\mathrm{div}\,}v\right\|
_{L_{R+1}^{1}}\left\| \rho _{\phi }\right\| _{L_{R}^{\infty }}+\left\|
v\right\| _{L_{R+1}^{\infty }}\left\| D\rho _{\phi }\right\| _{L_{R}^{1}}%
\right] .
\end{equation*}
To conclude, it is sufficient to note that $\Vert \rho _{\phi
}\Vert _{L_{R}^{\infty }}\leq \Vert \rho \Vert _{L_{r}^{\infty
}}\Vert J\phi ^{-1}\Vert _{L_{R}^{\infty }}$ and, denoting
by $[D\phi ^{-1}(y)]^{\ast }$ the adjoint matrix of $[D\phi
^{-1}(y)]$,
\begin{align*}
\left\| D\rho _{\phi }\right\| _{L_{R}^{1}}& \leq \left\| \lbrack
D\phi ^{-1}(\cdot )]^{\ast }\left( D\rho \circ \phi ^{-1}\right)
(\cdot )\,J\phi ^{-1}\left( \cdot \right) \right\|
_{L_{R}^{1}}+\left\| \left( \rho \circ
\phi ^{-1}\right) \cdot DJ\phi ^{-1}\right\| _{L_{R}^{1}} \\
& \leq \,\left\| D\phi ^{-1}\right\| _{L_{R}^{\infty }}\,\left\|
D\rho \right\| _{L_{R}^{1}}+\left\| \rho \right\| _{L_{R}^{\infty
}}\left\| DJ\phi ^{-1}\right\| _{L_{R}^{1}}.
\end{align*}
Given the bound the convergence follows by approximation.

For $d=1$, we simply have
\begin{eqnarray*}
\left| \int \mathcal{R}_{\varepsilon }[v,g](\phi (x))\rho
(x)dx\right| &=&\left| \int \mathcal{R}_{\varepsilon }[v,g](y)\rho
_{\phi }(y)dx\right|
\\
&\leq &\Vert \rho _{\phi }\Vert _{L_{R}^{\infty }}\int_{B\left(
R\right)
}\left| \mathcal{R}_{\varepsilon }[v,g](y)\right| dx \\
&\leq &C_{\rho }\Vert J\phi ^{-1}\Vert _{L_{R}^{\infty }}\Vert
v\Vert _{W_{R+2}^{1,1}}\Vert g\Vert _{L_{R+2}^{\infty }}
\end{eqnarray*}
where we have used Lemma~\ref{lemma:commutator-deterministic}. The
proof is complete.
\end{proof}

We are now ready to prove our first  uniqueness result of weak
$L^{\infty }$-solutions to the SPDE~\eqref{SPDE}.

\begin{theorem}
\label{thm SPDE n2} Assume that Hypothesis \ref{hy1} holds true.
Moreover,
assume Hypothesis \ref{hy2} for any $d\geq 1$ or simply $Db\in L_{loc}^{1}(%
\left[ 0,T\right] \times \mathbb{R})$ in the case $d=1$. Then, for every $%
u_{0}\in L^{\infty }$, there exists a unique weak $L^{\infty
}$-solution of
the Cauchy problem (\ref{SPDE}) which has the form $u(t,x,\omega )=u_{0}({%
\phi }_{t}^{-1}\left( \omega \right) x)$.
\end{theorem}

\begin{proof}
\textbf{Step 1.} By linearity we have to prove that a weak $L^{\infty }$%
-solution with initial condition $u_{0}=0$ vanishes identically. Let
us
denote by $u$ such a solution. For $y\in {\mathbb{R}}^{d}$ fixed, ${%
\varepsilon }>0$, let us choose the test function $\theta (x)=\vartheta _{{%
\varepsilon }}\left( y-x\right) $ in Definition \ref{def solution}.
Let us define $u^{\varepsilon }(t,\cdot )=\vartheta _{{\varepsilon
}}\ast u(t,\cdot )$. We get
\begin{equation*}
u^{\varepsilon }(t,y)=\int_{0}^{t}A_{\varepsilon }\left( s,y\right)
ds+\sum_{i=1}^{d}\int_{0}^{t}B_{\varepsilon }^{\left( i\right)
}\left( s,y\right) \circ dW_{s}^{i}
\end{equation*}
where
\begin{align*}
A_{\varepsilon }\left( t,y\right) &
=\int_{{}{\mathbb{R}}^{d}}u(t,x)\left\{
b(t,x)\cdot D_{x}[\vartheta _{{\varepsilon }}(y-x)]+{\mathrm{div}\,\,}%
b(t,x)\vartheta _{\varepsilon }(y-x)\right\} dx \\
B_{\varepsilon }^{\left( i\right) }\left( t,y\right) & =\int_{{}{\mathbb{R}}%
^{d}}u(t,x)D_{i}[\vartheta _{{\varepsilon }}\left( y-x\right) ]dx.
\end{align*}
All these functions of $y$, namely $u^{\varepsilon }(t,y)$,
$A_{\varepsilon
}\left( t,y\right) $, $B_{\varepsilon }^{\left( i\right) }\left( t,y\right) $%
, are bounded measurable in $t$, adapted, smooth (of class $C^{3}$
is required for the next computation) in $y$ (from (i) and (ii) of
definition
\ref{def solution}). As a minor remark, we know that $\int_{{}{\mathbb{R}}%
^{d}}u(t,x)\rho \left( x\right) dx$ is adapted when $\rho \in
C_{0}^{\infty }\left( {\mathbb{R}}^{d}\right) $, by definition of
solution, and by
approximation the same property holds for $\rho \in L^{1}\left( {\mathbb{R}}%
^{d}\right) $. From the Stratonovich version of
Kunita-It\^{o}-Wentzel formula (see~\cite[Th. 8.3 page 188]{K}), we
have
\begin{align*}
du^{\varepsilon }(t,{\phi }_{t}(x))& =A_{\varepsilon }\left( t,{\phi }%
_{t}(x)\right) dt+\sum_{i=1}^{d}B_{\varepsilon }^{\left( i\right) }\left( t,{%
\phi }_{t}(x)\right) \circ dW_{t}^{i} \\
& +\left( b\cdot Du^{\varepsilon }\right) (t,{\phi }_{t}(x))dt+%
\sum_{i=1}^{d}\left( D_{i}u^{\varepsilon }\right) (t,{\phi
}_{t}(x))\circ dW_{t}^{i}.
\end{align*}
But
\begin{equation*}
\left( D_{i}u^{\varepsilon }\right) (t,{y})=-{\int_{{\mathbb{R}}^{d}}}%
u(t,x)D_{i}[\vartheta _{{\varepsilon }}\left( y-x\right) ]\,dx
\end{equation*}
hence
\begin{equation*}
du^{\varepsilon }(t,{\phi }_{t}(x))=\left[ A_{\varepsilon }\left( t,{\phi }%
_{t}(x)\right) +\left( b\cdot Du^{\varepsilon }\right) (t,{\phi }_{t}(x))%
\right] dt
\end{equation*}
namely (recall that the initial condition is zero)
\begin{equation*}
u^{\varepsilon }(t,{\phi }_{t}(x))=-\int_{0}^{t}\mathcal{R}_{\varepsilon }%
\left[ b_{s},u_{s}\right] \left( {\phi }_{s}(x)\right) ds
\end{equation*}
where $\mathcal{R}_{\varepsilon }\left[ b_{s},u_{s}\right] $ is the
commutator defined above.

The Stratonovich version of Kunita-It\^{o}-Wentzel formula  as given in \cite[Th. 8.3 page 188]{K} is
optimized only with respect to the
martingale parts of the
processes: the theorems are stated for integrals of the form $%
\int_{0}^{t}f_{s}\left( x\right) dM_{s}$ where $M$ is a continuous
semimartingale, but the assumptions on $f$ are those necessary to
deal with martingales, not simply with processes of bounded
variations. For the bounded variation parts, which in our case take the form
$\int_{0}^{t}f_{s}\left( x\right) ds$, much weaker assumptions are
needed. This is in analogy with Lemma~\ref{lemma formula Ito} proved
above and applies in particular to
the random function $U\left( t,y\right) $ given by the integral $%
\int_{0}^{t}A_{\varepsilon }\left( s,y\right) ds$: $A_{\varepsilon
}$ does not satisfies the conditions of \cite[Th. 8.3 page 188]{K}
but $U\left( t,y\right) $ satisfies those of Lemma~\ref{lemma
formula Ito}. The Stratonovich version of Kunita-It\^{o}-Wentzel
formula extends to this case.

\smallskip \textbf{Step 2.} For every $\rho \in C_{0}^{\infty }({}{\mathbb{R}%
}^{d})$ (see the definition above) we have
\begin{equation*}
\int u^{\varepsilon }(t,{\phi }_{t}(x))\rho \left( x\right) dx=\int
u^{\varepsilon }(t,{y})\rho \left( {\phi }_{t}^{-1}(y)\right) J{\phi }%
_{t}^{-1}(y)dy.
\end{equation*}
Given $t\in \left[ 0,T\right] $, with probability one, $u^{\varepsilon }(t,{%
\cdot })$ converges weak-$\ast $ to $u(t,{\cdot })$ as $\varepsilon
\rightarrow 0$. Moreover, $P$-a.s., the function $y\mapsto \rho ({\phi }%
_{t}^{-1}(y))J{\phi }_{t}^{-1}(y)$ is integrable, since it is
continuous and with compact support. Hence, $P$-a.s.,
\begin{align*}
\lim_{\varepsilon \rightarrow 0}\int u^{\varepsilon }(t,{\phi
}_{t}(x))\rho
\left( x\right) dx& =\int u(t,{y})\rho \left( {\phi }_{t}^{-1}(y)\right) J{%
\phi }_{t}^{-1}(y)dy \\
& =\int u(t,{\phi }_{t}(x))\rho \left( x\right) dx.
\end{align*}
Therefore we have $P$-a.s.
\begin{equation*}
\int u(t,{\phi }_{t}(x))\rho \left( x\right) dx=\lim_{\varepsilon
\rightarrow 0}\int \left( \int_{0}^{t}\mathcal{R}_{\varepsilon
}\left[ b_{s},u_{s}\right] \left( {\phi }_{s}(x)\right) ds\right)
\rho \left( x\right) dx
\end{equation*}
\textit{If we prove that, given $\rho \in C_{0}^{\infty
}({}{\mathbb{R}}^{d}) $ and $t\in \left[ 0,T\right] $, this
$P$-a.s.-limit (which exists) is zero, then we have that $u$ is
identically zero (because ${\phi }_{t}$ is a bijection). }

\smallskip \textbf{Step 3.} Let us check, by means of Corollary \ref
{corollary commutator}, that
\begin{equation*}
s\mapsto \int_{\mathbb{R}^{d}}\mathcal{R}_{\varepsilon }\left[ b_{s},u_{s}%
\right] \left( {\phi }_{s}(x)\right) \rho \left( x\right) dx
\end{equation*}
satisfies the assumptions of Lebesgue dominated convergence theorem ($P$%
-a.s.) on $[0,T]$. 
This will complete the proof.

Assume $\rho \in C_{r}^{\infty }(\mathbb{R}^{d})$ for some $r>0$ and
define the r.v.
$R=\sup_{s\in \lbrack
0,T],x\in B(r)}|\phi _{s}(x)|$ so that the maps $\phi _{s}$ send the
support of $\rho $ strictly in the ball of radius $R$.

Let us give the details of the proof in the case $d>1$, the case
$d=1$ being
similar and easier. By Corollary \ref{corollary commutator}, for every $s\in %
\left[ 0,T\right] $ we have

\begin{equation*}
\begin{split}
& \left| \int \mathcal{R}_{\varepsilon }[b_{s},u_{s}](\phi
_{s}(x))\rho
(x)dx\right| \\ & \qquad \leq C_{\rho }\Vert u_{s}\Vert _{L_{R+1}^{\infty }}\left[ \Vert {%
\mathrm{div}}\,b_{s}\Vert _{L_{R+1}^{1}}\Vert J\phi _{s}^{-1}\Vert
_{L_{R}^{\infty }}+\Vert b_{s}\Vert _{L_{R+1}^{\infty }}(\Vert D\phi
_{s}^{-1}\Vert _{L_{R}^{\infty }}+\Vert DJ\phi _{s}^{-1}\Vert _{L_{R}^{1}})%
\right]
\\ &\qquad  \leq C_{\rho }\left[ \Vert {\mathrm{div}}\,b_{s}\Vert
_{L_{R+1}^{1}}\Vert J\phi _{s}^{-1}\Vert _{L_{R}^{\infty }}+\Vert
D\phi _{s}^{-1}\Vert _{L_{R}^{\infty }}+\Vert DJ\phi _{s}^{-1}\Vert
_{L_{R}^{1}}\right]
\end{split}
\end{equation*}
by the global boundedness of $b$ and $u$. From the properties of the
stochastic flow ${\phi }$ we know that $P(R<\infty )=1$ and that $%
(s,x)\mapsto D{\phi }_{s}^{-1}(x)$ is $P$-a.s. continuous. Hence the term $%
\sup_{s\in \left[ 0,T\right] }\Vert J{\phi }_{s}^{-1}\Vert
_{L_{R}^{\infty }} $ and $\sup_{s\in \left[ 0,T\right] }\Vert D\phi
_{s}^{-1}\Vert
_{L_{R}^{\infty }}$ are $P$-a.s. finite. Moreover, $\int_{0}^{T}\Vert {%
\mathrm{div}\,}b_{s}\Vert _{L_{R+1}^{1}}ds<\infty $. So it remains to
show that
\begin{equation}
\int_{0}^{T}\left\| DJ{\phi }_{s}^{-1}\right\| _{L_{R}^{1}}ds<\infty ,\quad P%
\text{-a.s.}  \label{bound on DJinvX}
\end{equation}
where $R$ is a positive r.v. which is $P$-a.s. finite.
This bound will follow from a similar
bound where $R$ is replaced by an arbitrary positive number.
Moreover, since by Lemma~\ref{fra} the equation for ${\phi
}_{s}^{-1}$ is equal to the equation for ${\phi }_{s} $ (up to a
sign and inversion of time) we can use Theorem~\ref{iac} to
conclude. The proof is complete.
\end{proof}

\medskip Let us formulate our second main result which basically
 only requires Hypothesis \ref{hy1} but with $\alpha > 1/2$  (for any
$d \ge 1$). Here   we do not need the regularity results on the
derivatives of $J \phi$ in Sobolev spaces.

\begin{theorem}
\label{thm SPDE n3} Assume that Hypothesis \ref{hy1} holds true with
$\alpha > 1/2$. Moreover assume
  that  $\mathrm{div}\, b\in
L^1_\loc([0,T]\times\RR^d)$. Then, for  every  $ u_{0}\in L^{\infty
}$, there exists a unique weak $L^{\infty }$-solution of
the Cauchy problem (\ref{SPDE}) which has the form $u(t,x,\omega )=u_{0}({%
\phi }_{t}^{-1}\left( \omega \right) x)$.
\end{theorem}

The proof requires the following lemma, in which we provide a
   special bound
  for the commutator.

\begin{lemma} \label{chissa}
Given $ v    \in L_{loc}^{\infty}\left(
\mathbb{R}^{d},\mathbb{R}^{d}\right) $, ${\mathrm{div}\,}v\in
L_{loc}^{1}\left(  \mathbb{R}^{d}\right)$, $g    \in
L_{loc}^{\infty}\left(  \mathbb{R}^{d}\right)$.

\smallskip (i) For any $\rho \in C_r^{\infty} (  \mathbb{R}^{d})$
and for sufficiently small $\varepsilon$ we have,  for some positive
constant $C_r$,
 \begin{align}\label{dop}
\left|\int  \mathcal{R}_{\varepsilon}\left[  g,v\right] (  x)
\rho(x) dx \right| \le
 C_r \| g\|_{L^{\infty}_{r+1}} \|
\rho\|_{L^{\infty}_{r}}
 \| {\mathrm{div}\,} v
\|_{_{L_{r+1}^{1}}} \\ \nonumber
 + \,  \left|
  \iint g(x') D_x \vartheta_{\eps} ({x-x'}) \,
\big( \rho(x) - \rho(x') \big)\, [v(x)-v(x')] \, dx dx' \right|.
\end{align}
(ii) If  in addition  there exists $\theta \in (0,1)$ such that $
 v \in C_{\loc}^{\theta}(
\mathbb{R}^{d}, \mathbb{R}^{d}),
$ then we have the  uniform bound
$$
\left|\int \mathcal{R}_{\varepsilon}\left[  g,v\right] (  x) \rho(x)
dx \right| \le
 C_r \| g\|_{L^{\infty}_{r+1}} \big( \|
\rho\|_{L^{\infty}_{r}}
 \| {\mathrm{div}\,} v
\|_{_{L_{r+1}^{1}}} +  \, [ v ]_{C^{ \theta}_{r+1}}  \, [ \rho
]_{C^{1 - \theta}_{r+1}} \big).
$$
\end{lemma}
\begin{proof} (i) We start as in the proof of Lemma~\ref{comm}.
  We write
$$
\int\mathcal{R}_{\varepsilon}\left[  g,v\right]  \left( x\right)
\rho\left(  x\right)  dx = J_1(\rho) + J_2 (\rho), $$ where
$$
J_1(\rho) =    \iint g(x') \vartheta_{\varepsilon}(x-x')D_{x}\rho(x)
[v(x)-v(x')]   \, dx dx'
$$
$$
 J_2 (\rho) =  \iint g(x')  \rho(x)
  \vartheta_{\varepsilon}(x-x')[ {\mathrm{div}\,}  v(x)- {\mathrm{div}\,} \,  v(x') ] \, dx dx' .
$$
Let us estimate $J_2$. By changing variables, $x = \eps y + y'$,
 $x' =y'$,
$$
J_2 (\rho) =  \frac{1}{\eps^d}\iint g(x')  \rho(x)
  \vartheta(\frac{x-x'}{\eps})
  [ {\mathrm{div}\,}  v(x)- {\mathrm{div}\,} \,  v(x') ] \, dx dx'
$$
$$
=\iint g(y')  \rho( \eps y + y' )
  \vartheta( y)
  [ {\mathrm{div}\,}  v(\eps y + y')- {\mathrm{div}\,} \,  v(y') ] \, dy dy'.
$$
Hence
$$
|J_2 (\rho)| \le \| g\|_{L^{\infty}_{r+1}} \|
\rho\|_{L^{\infty}_{r}} \iint
 ( |{\mathrm{div}\,}  v(\eps y + y')| +  |{\mathrm{div}\,} \,  v(y') |) \, dy dy'
$$
$$
\le  C_r \| g\|_{L^{\infty}_{r+1}} \| \rho\|_{L^{\infty}_{r}}
 \| {\mathrm{div}\,} v
\|_{_{L_{r+1}^{1}}}
$$
In order to estimate $J_1$, we note that
$$
J_1(\rho) =    \iint g(x') \vartheta_{\varepsilon}(x-x')D_{x} \big(
\rho(x) - \rho(x') \big)\,   [v(x)-v(x')]   \, dx dx'
$$
$$
= -  \iint g(x') D_x \vartheta_{\eps} ({x-x'}) \, \big( \rho(x) -
\rho(x') \big)\, [v(x)-v(x')] \, dx dx'
$$ $$ -
\iint g(x') \vartheta_{\varepsilon}(x-x') \big( \rho(x) - \rho(x')
\big)\,   {\mathrm{div}\,} v(x)   \, dx dx' $$ $$ = J_{11}(\rho) +
J_{12}(\rho).
$$
Let us treat $J_{12}(\rho)$. We find
$$
|J_{12}(\rho)| \le 2 \| g\|_{L^{\infty}_{r+1}} \| \rho
 \|_{L^{\infty}_{r+1}} \int  |{\mathrm{div}\,} v(x)|   \, dx  \int
  \vartheta_{\varepsilon}(x-x')
  \,  dx' \le 2\| g\|_{L^{\infty}_{r+1}} \| \rho
 \|_{L^{\infty}_{r+1}}  \| {\mathrm{div}\,} v
\|_{_{L_{r+1}^{1}}}.
$$
(ii) We only have to estimate  $J_{11}(\rho)$.
 We get
$$
|J_{11}(\rho)| \le \left|
  \frac{1}{\eps^{d +1}}
  \iint g(x') D_x \vartheta (\frac{x-x'}{\eps}) \,
\big( \rho(x) - \rho(x') \big)\, [v(x)-v(x')] \, dx dx' \right|
$$
$$
\le \frac{1}{\eps} [ v ]_{C^{\theta}_{r+1}}  \, [ \rho ]_{C^{1-
\theta}_{r+1}} \,
 \| g \|_{L^{\infty}_{r+1}}\, \frac{1}{\eps^{d }}
  \iint  | D_x \vartheta (\frac{x-x'}{\eps})|  \,
 |x - x'|^{  }   dx dx'
$$
$$
 \le
   \frac{1}{\eps^{d }} \iint  | D_x \vartheta (\frac{x-x'}{\eps})|  \,
    dx dx' \le C \, [ v ]_{C^{\theta}_{r+1}}  \, [ \rho
]_{C^{1- \theta}_{r+1}} \,
 \| g \|_{L^{\infty}_{r+1}},
$$
where $C$ is independent on $\eps$. The proof is complete.
\end{proof}

\medskip The previous result is now   extended to the case in which
commutators are composed with  a flow.

\begin{corollary}
 \label{dim1}Let $\phi$ be a $C^{1}$-diffeomorphism of
 $\mathbb{R}^{d}$ ($J \phi$ denotes its Jacobian). Assume
  that there exists $\theta \in (0,1)$ such that
  $     v \in C_{\loc}^{\theta}(
\mathbb{R}^{d}, \mathbb{R}^{d}) $, ${\mathrm{div}\,}v\in
L_{loc}^{1}\left( \mathbb{R}^{d}\right)$, $g \in
L_{loc}^{\infty}\left( \mathbb{R}^{d}\right)$. Moreover assume that
 $
 \, J \phi \in C_{\loc}^{1- \theta}( \mathbb{R}^{d}).
$  Then, for any $\rho \in C_r^{\infty}(  \mathbb{R}^{d})$ and any
$R>0$ such that $\mathrm{supp}( \rho\circ \phi^{-1}) \subseteq
B(R)$,
 we have the  uniform bound
\begin{equation} \label{f56}
\begin{split}
& \Big | \int\mathcal{R}_{\varepsilon}\left[ g,v\right] \left(
\phi\left( x\right)  \right)  \rho\left( x\right)  dx \Big|
\le
 C_r \| g\|_{L^{\infty}_{r+1}}  \|
\rho\|_{L^{\infty}_{r}}\, \| J\phi^{-1}\| _{L_{R}^{\infty}}\,
 \| {\mathrm{div}\,} v
\|_{_{L_{r+1}^{1}}}
\\ & +  \, C_r \| g\|_{L^{\infty}_{r+1}} [ v ]_{C^{ \theta}_{r+1}} ( \|
D\phi^{-1}\| _{C^{1- \theta }_{R}} \,  \| D\rho\| _{L_{r}^{\infty}}
\,
 \, + \, \| \rho \|_{L_{r}^{\infty}}
 [D\phi^{-1}]_{C^{1- \theta }_{R}}) .
\end{split}
\end{equation}
 In addition,
\[
\lim_{\varepsilon\rightarrow0}\int\mathcal{R}_{\varepsilon}\left[
g,v\right] \left(  \phi\left(  x\right)  \right)  \rho\left(
x\right)  dx=0.
\]
\end{corollary}
\begin{proof}
 By changing variable, we have $ \int\mathcal{R}_{\varepsilon}[
g,v] (  \phi(  x) )  \rho( x) dx=\int\mathcal{R}_{\varepsilon}[ g,v]
(  y)  \rho_{\phi}(  y) dx $ where the function
 $$ \rho_{\phi}(  y)  =\rho(  \phi^{-1}( y) )
 J\phi^{-1}(  y)
 $$
 has the support strictly contained in the ball
of radius $R$. Clearly, $
 \|  \rho_{\phi}\| _{L_{R}^{\infty}}\leq \| \rho\|
_{L_{r}^{\infty}}\| J\phi^{-1}\| _{L_{R}^{\infty}}. $

To prove the result, we  have to check that
 Lemma~\ref{chissa} can be applied with $\rho_{\phi}$
  instead of $\rho$.
  This follows since
$$
[\rho_{\phi} ]_{ C^{1 - \theta }_{R}} \le
 \| J\phi^{-1}\|
_{L_{R}^{\infty}} \, [\rho ({\phi}^{-1} (\cdot) ) ]_{ C^{1- \theta
}_{R}} \, + \, \| \rho \|_{L_{r}^{\infty}}
 [J\phi^{-1}]_{C^{1- \theta }_{R}}
$$
$$
\le  \| D\phi^{-1}\| _{L_{R}^{\infty}} \,  \| D\rho\|
_{L_{r}^{\infty}} \, [ D\phi^{-1}]_{C^{1- \theta }_{R}}
 \, + \, \| \rho \|_{L_{r}^{\infty}}
 [D\phi^{-1}]_{C^{1- \theta }_{R}}.
$$
\end{proof}

\textbf{Proof of Theorem \ref{thm SPDE n3} }. We follow the proof of
Theorem \ref{thm SPDE n2}. The first two steps are just the same.
The only change is in  Step 3.

\smallskip\textbf{Step 3.} We have to check that
$$ s\mapsto\int_{\mathbb{R} ^d}
\mathcal{R}_{\varepsilon}\left[  u_{s} ,b_{s}\right] \left(
{\phi}_{s}(x)\right)  \rho\left(  x\right) dx
$$ satisfies
the assumptions of Lebesgue dominated convergence theorem ($P$-a.s.)
on $[0,T]$.
  Assume $\rho \in C_r^\infty(
\mathbb{R}^{d})$ for some $r > 0$ and define the random variable
 $R = \sup_{s\in[0,T],x\in B(r)}
|\phi_s(x)|$ so that the maps $\phi_s$ send the support of $\rho$
strictly in the ball of radius $R$. Note that
$$
\sup_{s\in[0,T],x\in B(R)} |J\phi_s^{-1}(x)| < \infty.
$$
 Recall that $ D\phi^{-1}_s $ is $P$-a.s. locally $\alpha'$-Holder
continuous, uniformly in $s \in [0,T]$,
 for any $\alpha' \in (0, \alpha)$. Since $\alpha > 1/2$, we infer
  by Corollary  \ref{dim1} with $\theta = 1/2$
 $$
\left|\int \mathcal{R}_{\varepsilon}\left[  u_{s} ,b_{s}\right]
\left( {\phi}_{s}(x)\right)  \rho\left(  x\right) dx \right| \le
 C_{\rho} \| u_s\|_{L^{\infty}_{R}}
 \big(
 \| {\mathrm{div}\,} b_s
\|_{_{L_{R+1}^{1}}} +  \, [ b_s ]_{C^{1/2}_{R+1}}  \, \|
 D\phi^{-1}_s\| _{C^{1/2}_{R}}  \,
 \big).
$$
From the properties of the stochastic flow ${\phi}$ we know that $P(R
< \infty)=1$ and that $(  s,x) \mapsto D{\phi}_{s}^{-1}( x) $ is
$P$-a.s. continuous. Hence the terms $\int_{0}^{T}\| {\mathrm{div}\,
}b_{s}\| _{L_{R}^{1}}ds$ and $\sup_{s\in\left[ 0,T\right]  }([ b_s
]_{C^{1/2}_{R+1}}  \, \|
 D\phi^{-1}_s\| _{C^{1/2}_{R}})$ are $P$-a.s. finite
and so we can apply  the dominated convergence theorem. The proof is
complete.

\section{Examples}

\label{sec:examples}

In this section we give two classes of examples. First, we recall a
classical example of non-uniqueness for the deterministic transport equation
and we observe the improvements obtained by random perturbation (we call it
a `positive example'). Notice however that other relevant examples of
deterministic non-uniqueness, like the one of N. Depauw \cite{Dep}, are not
covered by the results of our work since there $b$ is not H\"{o}lder
continuous in the space variable.

Second, we show by means of two `negative examples' that it is not clear how
to extend the approach of this paper to random fields $b$ and nonlinear SPDEs.

\subsection{A positive example}

\label{sec:positive-example} Without noise, the transport equation with an
H\"{o}lder vector field is not necessarily well-posed. A counter-example can
be easily constructed in $1$d. Consider the function
\begin{equation*}
b(x)=\frac{1}{1-\gamma }\,\mathrm{sign}\,(x)\left( |x|\wedge R\right)
^{\gamma },\quad \gamma \in \left( 0,1\right)
\end{equation*}
(the simplest but less symmetric case $b(x)=\,\frac{1}{1-\gamma }\left(
|x|\wedge R\right) ^{\gamma }$ is similar, and many other variants can be
treated in the same way), where $R>0$ is introduced only to have
boundedness. This function is $C_{b}^{\gamma }\left( \mathbb{R}\right) $ and
$\mathrm{div}b=b^{\prime }\in L^{p}\left( \mathbb{R}\right) $ for all $p\in
\left( 1,1/\left( 1-\gamma \right) \right) $. Hence it satisfies hypothesis
\ref{hy1} for every $\gamma \in \left( 0,1\right) $, hypothesis \ref{hy2}
for $\gamma \in \left( \frac{1}{2},1\right) $ (because of the restriction $%
p>2$) and the condition $Db\in L_{loc}^{1}\left( \left[ 0,T\right] \times
\mathbb{R}\right) $ for every $\gamma \in \left( 0,1\right) $. Hence the
stochastic flow of diffeomorphisms exists and the stochastic transport
equation is well posed in $L^{\infty }$ and in $BV_{\text{loc}}$ (see Appendix~\ref{sec:bv} below) for all $\gamma \in
\left( 0,1\right) $.

On the contrary, the deterministic transport equation is not well posed; let
us recall why. The Cauchy problem
\begin{equation*}
x^{\prime }(t)=b\left( x\left( t\right) \right) ,\quad t\geq 0,\quad x\left(
0\right) =x_{0}
\end{equation*}
has a unique solution for all $x_{0}\neq 0$, denoted by $\phi _{t}\left(
x_{0}\right) $. For $x_{0}=0$ we have two extremal solutions $x_{+}\left(
t\right) $ e $x_{-}\left( t\right) $, $x_{+}\left( t\right) =t^{\frac{1}{%
1-\gamma }}$ and $x_{-}\left( t\right) =-t^{\frac{1}{1-\gamma }}$ for small $%
t$. In addition, for $x_{0}=0$, we have the solution $x\left( t\right)
\equiv 0$, and the solutions $x\left( t\right) =x_{\pm }\left(
t-t_{0}\right) 1_{t\geq t_{0}}$ for every $t_{0}\geq 0$. Given $t>0$ and $%
x\in \left[ x_{-}\left( t\right) ,x_{+}\left( t\right) \right] $, there is a
unique number $t_{0}\left( t,x\right) \geq 0$ such that $x_{sign(x)}\left(
t-t_{0}\left( t,x\right) \right) =x$. The function $\phi _{t}$ maps $\mathbb{%
R}\diagdown \left\{ 0\right\} $ one to one on $\left( -\infty ,x_{-}\left(
t\right) \right) \cup \left( x_{+}\left( t\right) ,\infty \right) $; $\phi
_{t}^{-1}$ will be its inverse, between these sets. With these notations,
given $u_{0}\in L^{\infty }$ and two bounded measurable functions $\gamma
_{+},\gamma _{-}:[0,\infty )\rightarrow \mathbb{R}$, define the function
\begin{equation}
u_{\gamma _{\pm }}\left( t,x\right) =\left\{
\begin{array}{ccc}
u_{0}\left( \phi _{t}^{-1}\left( x\right) \right) & \text{for} &
x>x_{+}\left( t\right) \\
\gamma _{+}\left( t_{0}\left( t,x\right) \right) & \text{for} & 0\leq x\leq
x_{+}\left( t\right) \\
\gamma _{-}\left( t_{0}\left( t,x\right) \right) & \text{for} & x_{-}\left(
t\right) \leq x<0 \\
u_{0}\left( \phi _{t}^{-1}\left( x\right) \right) & \text{for} &
x<x_{-}\left( t\right)
\end{array}
\right. .  \label{solutions}
\end{equation}
These are weak $L^{\infty }$ solutions, for every $\gamma _{+},\gamma _{-}$,
of the deterministic transport equation with the same initial condition $%
u_{0}$. For instance, if $u_{0}=1_{x>0}$ and $\gamma _{+}=\gamma _{-}\equiv
a $ for a constant value $a$, the shape of $u_{\gamma _{\pm }}$ can be
easily worked out. All these functions are solutions both in $L^{\infty }$ and
in $BV_{\text{loc}}$, corresponding to the same $BV_{\text{loc}}$ initial condition $u_{0}$.

\subsection{Negative examples}

It would be interesting to generalize the results of this paper to random
vectorfields $b(t, x,\omega)$, possibly adapted. However our approach faces
a fundamental difficulty: it is very easy to exhibit a counterexample which
shows that in some cases the regularizing effect disappears. Consider in one
dimension the case
\begin{equation*}
b\left(t, x,\omega\right) =\sqrt{\left| x-W_{t}\left( \omega\right) \right| }%
,
\end{equation*}
namely the stochastic differential equation
\begin{equation*}
dX_{t}^{x}=b(t, X_{t}^{x},\cdot ) dt+dW_{t},\quad t\geq0,\quad X_{0}^{x}=x .
\end{equation*}
If $\left( X_{t}^{x}\right) $ is a solution, then $Y_{t}^{x}=X_{t}^{x}-W_{t}$
solves
\begin{equation*}
dY_{t}^{x}=\sqrt{\left| Y_{t}^{x}\right| }dt\quad t\geq0,\quad Y_{0}^{x}=x
\end{equation*}
and viceversa. Hence the non-uniqueness for the latter equation transfer to
the former. In terms of stochastic transport equation, an equation of the
form
\begin{equation*}
\partial_{t}u\left( t,x\right) +\left( b_{0}\left( x-W\left( t\right)
\right) \cdot D u\left( t,x\right) \right) dt+D u\left( t,x\right) \circ
dW\left( t\right) =0
\end{equation*}
may have several pathologies if $b_{0}$ is only H\"{o}lder continuous.

Unfortunately the previous example is not so artificial: something similar
happens in the nonlinear case. Let us argue only formally. Consider the
example in $\mathbb{R}^{2}$%
\begin{equation*}
\partial _{t}\xi \left( t,x\right) +\left( u\left( t,x\right) \cdot D\xi
\left( t,x\right) \right) dt+\sum_{i=1}^{d}D_{i}\xi \left( t,x\right) \circ
dW^{i}\left( t\right) =0
\end{equation*}
where $\xi =\partial _{2}u_{1}-\partial _{1}u_{2}$. This is the vorticity
equation of a 2D ideal fluid described by a stochastic version of Euler
equation. Following \cite{MajBar}, this equation is (formally) equivalent to
the family of stochastic ordinary equations depending on a parameter $a\in
\mathbb{R}^{2}$%
\begin{equation*}
dX_{t}^{a}=\left[ \int_{\mathbb{R}^{2}}K(X_{t}^{a}-X_{t}^{a^{\prime }})\xi
_{0}(X_{t}^{a^{\prime }})da^{\prime }\right] dt+dW_{t}
\end{equation*}
for a suitable kernel $K$, $\xi _{0}$ being the initial condition of the
vorticity equation. This problem is equivalent to
\begin{equation*}
dY_{t}^{a}=\left[ \int_{\mathbb{R}^{2}}K(Y_{t}^{a}-Y_{t}^{a^{\prime }})\xi
_{0}(X_{t}^{a^{\prime }})da^{\prime }\right] dt
\end{equation*}
by the change of variable $Y_{t}^{a}=X_{t}^{a}-W_{t}$, and the equation for $%
\left( Y_{t}^{a}\right) $ corresponds to the classical vorticity equation
\begin{equation*}
\frac{\partial _{t}\xi ^{\prime }\left( t,x\right) }{\partial t}+\left(
u^{\prime }\left( t,x\right) \cdot D\xi ^{\prime }\left( t,x\right) \right)
dt=0\qquad \qquad \xi ^{\prime }=\partial _{2}u_{1}^{\prime }-\partial
_{1}u_{2}^{\prime }
\end{equation*}
with initial condition $\xi _{0}$. This means that the stochastic vorticity
equation is (at least formally) equivalent to the deterministic one. There
is no advantage to introduce that kind of stochastic perturbation.

\appendix

\section{Existence and uniqueness in $BV_{loc}$}

\label{sec:bv} The results proved in Section~\ref{sec:ito} on the stochastic flow allow
one to prove several existence and uniqueness results in spaces more regular
than $L^{\infty }$. We describe here the case $BV_{loc}$, as a less trivial
example. Let us emphasize that only the assumption $b\in L^{\infty }\left(
0,T;C_{b}^{\alpha }({\mathbb{R}}^{d};{\mathbb{R}}^{d})\right) $ is needed;
no condition on $\mathrm{div}\,b$ is imposed. This is again at variance with
the deterministic case.

For more details on the functions of locally bounded variation see \cite{GMS}%
. Let us recall that a function $v\in L_{loc}^{1}(\mathbb{R}^{d})$ is said
to be of locally bounded variation, $v\in BV_{loc}(\mathbb{R}^{d})$ if its
distributional derivatives $D_{i}v$, $i=1,...,d$, are signed Radon measures.
We denote by $Dv$ the vector valued measure with components $D_{i}v$. We
have
\begin{equation*}
\int_{\mathbb{R}^{d}}\vartheta (x)\cdot Dv(dx)=-\int_{\mathbb{R}^{d}}v(x)%
\mathrm{div}\,\vartheta (x)dx
\end{equation*}
for all vector fields $\vartheta \in C_{0}^{1}(\mathbb{R}^{d},\mathbb{R}%
^{d}) $, and also for all $\vartheta \in C_{0}(\mathbb{R}^{d},\mathbb{R}%
^{d}) $ such that the distributional divergence (which has compact support) $%
\mathrm{div}\,\vartheta $ is of class $L^{1}(\mathbb{R}^{d})$.

When $a\in C(\mathbb{R}^{d},\mathbb{R}^{d})$ is a given vector field and $%
v\in BV_{loc}(\mathbb{R}^{d})$, the notation $a\cdot Dv$ stands for the
(scalar) distribution $\theta\mapsto\int_{\mathbb{R}^{d}}\theta(x)a(x)\cdot
Dv(dx)$, $\theta\in C_{0}(\mathbb{R}^{d})$. This is the meaning of the
notation $\int_{\mathbb{R}^{d}}\theta(x)b_{s}(x)\cdot Du_{s}(dx)$ used in
the next definition.

We denote by $L_{\mathcal{F}}^{\infty}(BV_{loc}(\mathbb{R}^{d}))$ the space
of all stochastic processes $u\in L^{\infty}(\Omega\times\lbrack0,T]\times
\mathbb{R}^{d})$ such that $u(\omega,t,\cdot)\in BV_{loc}(\mathbb{R}^{d})$
for a.e. $(\omega,t)\in\Omega\times\lbrack0,T]$, for every $\theta\in
C_{0}^{\infty}(\mathbb{R}^{d})$ and $\vartheta\in L^{\infty}(0,T;C_{0}(%
\mathbb{R}^{d},\mathbb{R}^{d}))$ the processes
\begin{equation*}
\int_{\mathbb{R}^{d}}\theta(x)u(t,x)dx,\quad\int_{\mathbb{R}%
^{d}}\vartheta(t,x)\cdot Du(t,dx)
\end{equation*}
are progressively measurable with respect to $(\mathcal{F}_{t})_{t\in
\lbrack0,T]}$, and $|Du(\omega,t)|$ (the total variation of the measure $%
Du(\omega,t,\cdot)$) has the following property:
\begin{equation}
\int_{0}^{T}|Du(\omega,t)|\left( B\left( r\right) \right) dt<\infty
\label{L1 BV}
\end{equation}
for all $r>0$, for $P$-a.e. $\omega\in\Omega$. We use again the notation $%
u_{t}$ for $u ( t,\cdot) $ below.

\begin{definition}
Let $b\in L^{\infty }\left( 0,T;C_{b}^{\alpha }({\mathbb{R}}^{d};{\mathbb{R}}%
^{d})\right) $ and $u_{0}\in BV_{loc}\left( \mathbb{R}^{d}\right) $. A
stochastic process $u\in L_{\mathcal{F}}^{\infty }\left( BV_{loc}\left(
\mathbb{R}^{d}\right) \right) $ is a $BV_{loc}$-solution of the Cauchy
problem (1) if, for every test function $\theta \in C_{0}^{\infty }\left(
\mathbb{R}^{d}\right) $, the process $\int_{\mathbb{R}^{d}}\theta \left(
x\right) u_{t}\left( x\right) dx$ has a continuous modification which is an $%
\mathcal{F}$-semimartin\-gale and
\begin{align*}
\int_{\mathbb{R}^{d}}u_{t}\left( x\right) \theta \left( x\right) dx& =\int_{%
\mathbb{R}^{d}}u_{0}\left( x\right) \theta \left( x\right) dx \\
& -\int_{0}^{t}\left( \int_{\mathbb{R}^{d}}\theta \left( x\right)
b_{s}\left( x\right) \cdot Du_{s}\left( dx\right) \right) ds \\
& +\sum_{i=1}^{d}\int_{0}^{t}\left( \int_{\mathbb{R}^{d}}u_{s}\left(
x\right) D_{i}\theta \left( x\right) dx\right) \circ dW_{s}^{i}.
\end{align*}
\end{definition}

If $u$ is a $BV_{loc}$-solution and $\mathrm{div}\,b\in L_{loc}^{1}\left( %
\left[ 0,T\right] \times\mathbb{R}^{d}\right) $, then $u$ is also a weak $%
L^{\infty}$-solution. Conversely, if $u$ is a weak $L^{\infty}$-solution, $%
u_{0}\in BV_{loc}\left( \mathbb{R}^{d}\right) $ and $u\in L_{\mathcal{F}%
}^{\infty}\left( BV_{loc}\left( \mathbb{R}^{d}\right) \right) $, then $u$ is
a $BV_{loc}$-solution.

\begin{theorem}
If $b\in L^{\infty }\left( 0,T;C_{b}^{\alpha }({\mathbb{R}}^{d};{\mathbb{R}}%
^{d})\right) $ and $u_{0}\in BV_{loc}\left( \mathbb{R}^{d}\right) $, there
exists a unique $BV_{loc}$-solution $u$, given by $u\left( t,x\right)
=u_{0}\left( \phi _{t}^{-1}\left( x\right) \right) $.
\end{theorem}

\begin{proof}
\textbf{Step 1} (Existence). Let us mention a preliminary known fact. If $%
u_{0}\in BV_{loc}\left( \mathbb{R}^{d}\right) $ and $\varphi $ is a
diffeomorphism of $\mathbb{R}^{d}$ (differentiable in both directions with
continuous derivatives), then $u_{0}\circ \varphi \in BV_{loc}\left( \mathbb{%
R}^{d}\right) $ and the signed Radon measures $D_{i}\left( u_{0}\circ
\varphi \right) \left( dx\right) $, $i=1,...,d$, are defined by
\begin{equation}
\int_{\mathbb{R}^{d}}\theta \left( x\right) D_{i}\left( u_{0}\circ \varphi
\right) \left( dx\right) =\int_{\mathbb{R}^{d}}\theta \left( \varphi
^{-1}\left( y\right) \right) J\varphi ^{-1}\left( y\right) \left(
D_{i}\varphi \right) \left( \varphi ^{-1}\left( y\right) \right) \cdot
Du_{0}\left( dy\right)   \label{identity BV}
\end{equation}
for every $\theta \in C_{0}^{\infty }\left( \mathbb{R}^{d}\right) $. To
prove this claim we use the characterization of $BV_{loc}$ functions $v$ as
the weak $L_{loc}^{1}\left( \mathbb{R}^{d}\right) $ limits of $%
W_{loc}^{1,1}\left( \mathbb{R}^{d}\right) $ (or $C_{0}^{\infty }(\mathbb{R}%
^{d})$) functions $v_{n}$ such that $\sup_{n}\int_{B\left( r\right) }\left|
Dv_{n}\left( x\right) \right| dx<\infty $ for every $r>0$ (see \cite{GMS},
Chapter 4). Take a sequence of functions $u_{0}^{n}\in C_{0}^{\infty }(%
\mathbb{R}^{d})$ such that $u_{0}^{n}\rightharpoonup u_{0}$ in $%
L_{loc}^{1}\left( \mathbb{R}^{d}\right) $, $\sup_{n}\int_{B\left( r\right)
}\left| Du_{0}^{n}\left( x\right) \right| dx<\infty $ for every $r>0$. The
function $u_{0}^{n}\circ \varphi $ is in $C_{0}^{1}(\mathbb{R}^{d})$, $%
u_{0}^{n}\circ \varphi \rightharpoonup u_{0}\circ \varphi $ in $%
L_{loc}^{1}\left( \mathbb{R}^{d}\right) $ and
\begin{align*}
\int_{B\left( r\right) }\left| D_{i}\left( u_{0}^{n}\circ \varphi \right)
\left( x\right) \right| dx& \leq \int_{B\left( r\right) }\left|
Du_{0}^{n}\left( \varphi \left( x\right) \right) \cdot D_{i}\varphi \left(
x\right) \right| dx \\
& =\int_{\varphi \left( B\left( r\right) \right) }\left| Du_{0}^{n}\left(
y\right) \cdot D_{i}\varphi \left( \varphi ^{-1}\left( y\right) \right)
\right| J\varphi ^{-1}\left( y\right) dy
\end{align*}
which implies that $\sup_{n}\int_{B\left( r\right) }\left| D\left(
u_{0}^{n}\circ \varphi \right) \left( x\right) \right| dx<\infty $ for every
$r>0$ and thus $u_{0}\circ \varphi \in BV_{loc}\left( \mathbb{R}^{d}\right) $%
. Similarly
\begin{equation*}
\int_{\mathbb{R}^{d}}\theta \left( x\right) D_{i}\left( u_{0}^{n}\circ
\varphi \right) \left( dx\right) =\int_{\mathbb{R}^{d}}\theta \left( \varphi
^{-1}\left( y\right) \right) J\varphi ^{-1}\left( y\right) D_{i}\varphi
\left( \varphi ^{-1}\left( y\right) \right) \cdot Du_{0}^{n}\left( y\right)
dy.
\end{equation*}
Up to a common subsequence, the measures $D_{i}\left( u_{0}^{n}\circ \varphi
\right) $ and $Du_{0}^{n}$ weakly converge to $D_{i}\left( u_{0}\circ
\varphi \right) $ and $Du_{0}$ respectively, on $B\left( r\right) $ for
every $r>0$ (Proposition 5 of \cite{GMS}, Chapter 4.1.1). We can take the
limit in the previous identity and get (\ref{identity BV}).

Let us prove the existence claim. Let $\phi _{t}$ be the stochastic flow of
diffeomorphisms given under the assumption $b\in L^{\infty }\left(
0,T;C_{b}^{\alpha }({\mathbb{R}}^{d};{\mathbb{R}}^{d})\right) $ and let us
set $u\left( t,x\right) =u_{0}\left( \phi _{t}^{-1}\left( x\right) \right) $%
. For every set $A\subset \mathbb{R}^{d}$, denote by $\phi \left( \omega
,A\right) $ the image of $\left[ 0,T\right] \times A$ under the mapping $%
\left( t,x\right) \mapsto \phi _{t}\left( \omega ,x\right) $. This mapping
is $P$-a.s. continuous, hence $\phi \left( \omega ,B\right) $ is a bounded
set $P$-a.s., for every bounded set $B$.

We have $u\in L^{\infty }\left( \Omega \times \left[ 0,T\right] \times
\mathbb{R}^{d}\right) $ and $\int_{\mathbb{R}^{d}}\theta \left( x\right)
u\left( t,x\right) dx$ is progressively measurable for every $\theta \in
C_{0}^{\infty }\left( \mathbb{R}^{d}\right) $. Since $u_{0}\in
BV_{loc}\left( \mathbb{R}^{d}\right) $, for a.e. $\omega \in \Omega $ we
have that for all $t\in \left[ 0,T\right] $ the function $x\mapsto
u_{0}\left( \phi _{t}^{-1}\left( x\right) \right) $ belongs to $%
BV_{loc}\left( \mathbb{R}^{d}\right) $ and, for $\theta \in L^{\infty
}\left( 0,T;C_{0}\left( \mathbb{R}^{d}\right) \right) $ and $i=1,...,d$,
\begin{equation*}
\int_{\mathbb{R}^{d}}\theta _{t}\left( x\right) D_{i}u_{t}\left( dx\right)
=\int_{\mathbb{R}^{d}}\theta _{t}\left( \phi _{t}^{-1}\left( y\right)
\right) J\phi _{t}^{-1}\left( y\right) \left( D_{i}\phi _{t}\right) \left(
\phi _{t}^{-1}\left( y\right) \right) \cdot Du_{0}\left( dy\right)
\end{equation*}
which is progressively measurable. Moreover, since $|D_{i}u_{t}|\left(
B\left( r\right) \right) $ is the supremum of the quantity $\left| \int_{B\left( r\right)
}\theta \left( x\right) D_{i}u_{t}\left( dx\right) \right| $ over all $%
\theta \in C_{0}^{1}\left( B\left( r\right) \right) $ such that $\left|
\theta \left( x\right) \right| \leq 1$ for all $x\in B\left( r\right) $, and
from the previous identity we get
\begin{align*}
\left| \int_{B\left( r\right) }\theta \left( x\right) D_{i}u_{t}\left(
dx\right) \right| & \leq \int_{\phi _{t}\left( B\left( r\right) \right)
}J\phi _{t}^{-1}\left( y\right) \left| \left( D_{i}\phi _{t}\right) \left(
\phi _{t}^{-1}\left( y\right) \right) \right| \left| Du_{0}\right| \left(
dy\right)  \\
& \leq \left| Du_{0}\right| \left( \phi \left( \omega ,B\left( r\right)
\right) \right) \left\| J\phi _{t}^{-1}\left( \cdot \right) \left| \left(
D_{i}\phi _{t}\right) \left( \phi _{t}^{-1}\left( \cdot \right) \right)
\right| \right\| _{L^{\infty }\left( \phi \left( \omega ,B\left( r\right)
\right) \right) }
\end{align*}
we deduce (\ref{L1 BV}) and even
$
P( \sup_{t\in[ 0,T] }|Du_{t}|( B( r)
) <\infty) =1
$.
Hence $u\in L_{\mathcal{F}}^{\infty }\left( BV_{loc}\left( \mathbb{R}%
^{d}\right) \right) $.

We may now repeat the proof of Theorem \ref{them SPDE n1} and prove that $%
u\left( t,x\right) =u_{0}\left( \phi _{t}^{-1}\left( x\right) \right) $ is a
solution; the difference is that in Section~\ref{sec:transp} we assumed $%
\mathrm{div}\,b\in L_{loc}^{1}$ and $u$ was only $L^{\infty }$, while here $%
\mathrm{div}\,b$ is only a distribution but $u\in BV_{loc}\left( \mathbb{R}%
^{d}\right) $, so we have to write differently the integrals involving $%
b\cdot Du$.

\textbf{Step 2} (Uniqueness). We repeat the first part of the proof of
Theorem 23: by linearity, we treat the case $u_{0}=0$; we mollify a $BV_{loc}
$-solution $u\in L_{\mathcal{F}}^{\infty }\left( BV_{loc}\left( \mathbb{R}%
^{d}\right) \right) $, apply stochastic calculus and prove the same results
as in steps 1 and 2 of the proof of Theorem 23. The commutator $\mathcal{R}%
_{\varepsilon }\left[ u_{s},b_{s}\right] $ is always given by $\mathcal{R}%
_{\varepsilon }\left[ u_{s},b_{s}\right] =\vartheta _{\varepsilon }\ast
\left( b_{s}\cdot Du_{s}\right) -b_{s}\cdot D\left( \vartheta _{\varepsilon
}\ast u_{s}\right) $, but here we stress that $b_{s}\cdot Du_{s}$ is the
distribution having the meaning recalled at the beginning of the section (we
have to use compact support mollifiers). We have to prove that, $P$-a.s.,
\begin{equation*}
\lim_{\varepsilon \rightarrow 0}\int_{0}^{t}\left( \int \mathcal{R}%
_{\varepsilon }\left[ u_{s},b_{s}\right] \left( \phi _{s}\left( x\right)
\right) \rho \left( x\right) dx\right) ds=0
\end{equation*}
for every $\rho \in C_{0}^{\infty }\left( \mathbb{R}^{d}\right) $, namely
step 3 of the proof of Theorem 23. Equivalently, we have to prove that $P$%
-a.s.
\begin{equation*}
\lim_{\varepsilon \rightarrow 0}\int_{0}^{t}\left( \int \mathcal{R}%
_{\varepsilon }\left[ u_{s},b_{s}\right] \left( y\right) \theta _{s}\left(
y\right) dy\right) ds=0
\end{equation*}
where $\theta _{t}\left( x\right) =\rho \left( \phi _{t}^{-1}\left(
x\right) \right) J\phi _{t}^{-1}\left( x\right) $. This is,\ $P$-a.s., a
bounded continuous function of $\left( t,x\right) $, and $\theta _{t}\left(
\omega ,x\right) =0$ for all $x\notin \phi \left( \omega ,K\right) $ and all
$t\in \left[ 0,T\right] $, where $K$ is the support of $\rho $. One has
\begin{equation*}
\mathcal{R}_{\varepsilon }\left[ u_{s},b_{s}\right] \left( y\right) =\int
\vartheta _{\varepsilon }\left( y-z\right) \left[ b_{s}\left( z\right)
-b_{s}\left( y\right) \right] \cdot Du_{s}\left( dz\right)
\end{equation*}
hence (all functions are bounded so we may apply Fubini theorem)
\begin{equation*}
\int_{0}^{t}\left( \int \mathcal{R}_{\varepsilon }\left[ u_{s},b_{s}\right]
\left( y\right) \theta _{s}\left( y\right) dy\right) ds=\int_{0}^{t}\left(
\int f_{\varepsilon }\left( s,z\right) \cdot Du_{s}\left( dz\right) \right)
ds
\end{equation*}
where
\begin{equation*}
f_{\varepsilon }\left( \omega ,s,z\right) =\int \vartheta _{\varepsilon
}\left( y-z\right) \left[ b_{s}\left( z\right) -b_{s}\left( y\right) \right]
\theta _{s}\left( \omega ,y\right) dy.
\end{equation*}
Notice that, for $P$-a.e. $\omega \in \Omega $, $f_{\varepsilon }\left(
\omega ,s,z\right) =0$ for all $z\notin \mathcal{U}_{\varepsilon }\left(
\phi \left( \omega ,K\right) \right) $ and all $s\in \left[ 0,T\right] $.
Here $\mathcal{U}_{\varepsilon }\left( \phi \left( \omega ,K\right) \right) $
is an $\varepsilon $-neighbor of $\phi \left( \omega ,K\right) $ (assuming
that the support of $\vartheta $ is in $B\left( 1/2\right) $). And $%
f_{\varepsilon }\left( \omega ,s,\cdot \right) $ is continuous. Hence $\int
f_{\varepsilon }\left( s,z\right) \cdot Du_{s}\left( dz\right) $ is
meaningful. We have
\begin{equation*}
\left| \int_{0}^{t}\left( \int \mathcal{R}_{\varepsilon }\left[ u_{s},b_{s}%
\right] \left( y\right) \theta _{s}\left( y\right) dy\right) ds\right| \leq
\int_{0}^{t}\left( \int \left| f_{\varepsilon }\left( s,z\right) \right|
\left| Du_{s}\right| \left( dz\right) \right) ds
\end{equation*}
\begin{equation*}
\leq \int_{0}^{t}\left| Du_{s}\right| \left( \mathcal{U}_{1}\left( \phi
\left( \omega ,K\right) \right) \right) \left\| f_{\varepsilon }\left(
\omega ,s,\cdot \right) \right\| _{L^{\infty }\left( \mathcal{U}_{1}\left(
\phi \left( \omega ,K\right) \right) \right) }ds
\end{equation*}
for all $\varepsilon \leq 1$. Moreover, with $C\left( \omega \right)
=\left\| \theta _{\cdot }\left( \omega ,\cdot \right) \right\| _{L^{\infty
}\left( \left[ 0,T\right] \times \phi \left( \omega ,K\right) \right) }$, we
have
\begin{align*}
\left| f_{\varepsilon }\left( \omega ,s,z\right) \right| & \leq C\left(
\omega \right) \int_{\phi \left( \omega ,K\right) \cap B\left( z,\varepsilon
\right) }\vartheta _{\varepsilon }\left( y-z\right) \left| b_{s}\left(
z\right) -b_{s}\left( y\right) \right| dy \\
& \leq C\left( \omega \right) \int_{B\left( 1\right) }\vartheta \left(
h\right) \left| b_{s}\left( z\right) -b_{s}\left( z+\varepsilon h\right)
\right| dh\leq C^{\prime }\left( \omega \right) \varepsilon ^{\alpha }
\end{align*}
uniformly in $\left( s,z\right) $. Hence
\begin{equation*}
\left| \int_{0}^{t}\left( \int \mathcal{R}_{\varepsilon }\left[ u_{s},b_{s}%
\right] \left( y\right) \theta _{s}\left( y\right) dy\right) ds\right| \leq
C^{\prime }\left( \omega \right) \varepsilon ^{\alpha }\int_{0}^{T}\left|
Du_{s}\right| \left( \mathcal{U}_{1}\left( \phi \left( \omega ,K\right)
\right) \right) ds
\end{equation*}
which converges to zero as $\varepsilon \rightarrow 0$, for $P$-a.s. $\omega
\in \Omega $. The proof is complete.
\end{proof}


\section{The perturbative equation}

\label{sect perturbative}

In this section we give a pathwise formulation of the SPDE, that we call
perturbative equation. It does not involve stochastic integrals explicitly.
We use it to prove Theorem \ref{thm WZ 2}: the absence of stochastic
integrals simplify the analysis of some limits. Moreover, we think it may
have other applications.

For streamlining the notations in this section we will let
\begin{equation*}
W_{ts}=W_{t}-W_{s},\;\;\;t\geq s\geq 0,\;\;\;\;\;\;v_{t}(\theta )=\int_{%
\mathbb{R}^{d}}\theta (x)v(t,x)dx,\;\;t\in \lbrack 0,T],
\end{equation*}
for any $v\in L_{\mathrm{loc}}^{1}([0,T]\times \mathbb{R}^{d}),$ \ $\theta
\in C_{0}^{\infty }(\mathbb{R}^{d},\mathbb{R}^{k})$.

\begin{theorem}
\label{lemma:mild} Let $b\in L_{\mathrm{loc}}^{1}([0,T]\times \mathbb{R}^{d};%
\mathbb{R}^{d})$, $\mathrm{div}\,b\in L_{\mathrm{loc}}^{1}([0,T]\times
\mathbb{R}^{d})$ and $u_{0}\in L^{\infty }(\mathbb{R}^{d})$. A process $u\in
L^{\infty }(\Omega \times \lbrack 0,T]\times \mathbb{R}^{d})$ is a weak $%
L^{\infty }$-solution of the SPDE \eqref{SPDE} if and only if, for every $%
\theta \in C_{0}^{\infty }(\mathbb{R}^{d})$, $u_{t}(\theta )$ has a
continuous adapted modification and the perturbative equation
\begin{equation}
u_{t}(\theta )=u_{0}(\theta (\cdot +W_{t}))+\int_{0}^{t}ds\int
[b_{s}(x)\cdot D\theta (x+W_{ts})+\mathrm{div}\,b_{s}(x)\theta
(x+W_{ts})]u_{s}(x)dx  \label{ut}
\end{equation}
holds almost surely in $\omega \in \Omega $, for every $t\in \lbrack 0,T]$.
\end{theorem}

\begin{remark}
\label{remark sulla perturbativa}Given $\theta \in C_{0}^{\infty }\left(
\mathbb{R}^{d}\right) $ and $t\in \left[ 0,T\right] $, the integral on the
right-hand-side of equation (31) is a well defined random variable, in spite
of the only local integrability of $b$ and $\mathrm{div}\,b$. Indeed, for $P$%
-a.e. $\omega \in \Omega $, $W_{\cdot }\left( \omega \right) $ is bounded on
$\left[ 0,T\right] $ and thus $\theta \left( x+W_{ts}\left( \omega \right)
\right) $ and $D\theta \left( x+W_{ts}\left( \omega \right) \right) $ vanish
for $x$ outside a bounded set $B_{\omega }\subset \mathbb{R}^{d}$ ($%
B_{\omega }$ depends on $\omega $ but not on $t$ and $s$ in $\left[ 0,T%
\right] $).
\end{remark}

\begin{proof}
\textbf{Step 1}. Let us prove that a solution $u_{t}\left( \theta \right) $
of the perturbative equation is an $\mathcal{F}$-semimartingale and that
the equation of Definition~\ref{def solution} is verified. We apply It\^{o}
formula to the process $F_{0}\left( W_{t}\right) $ and It\^{o}%
-Wentzell-Kunita formula (see~\cite{K}) to $F_{1}\left( t,W_{t}\right) $
where
\begin{align*}
F_{0}\left( y\right) & =u_{0}\left( \theta \left( \cdot +y\right) \right)  \\
F_{1}\left( t,y\right) & =\int_{0}^{t}ds\int_{\mathbb{R}^{d}}\left[
b_{s}\left( x\right) \cdot D\theta \left( x+y-W_{s}\right) +\mathrm{div}\,%
b_{s}\left( x\right) \theta \left( x+y-W_{s}\right) \right] u_{s}\left(
x\right) dx
\end{align*}
$t\in \left[ 0,T\right] $, $x\in \mathbb{R}^{d}$. Notice that $F_{1}\left(
t,y,\omega \right) $, although being a random field, is of bounded variation
(and more) in $t$, namely it has no martingale part, and
\begin{align*}
f_{1}\left( t,y\right) & =\frac{\partial F_{1}\left( t,y\right) }{\partial t%
} \\
& =\int_{\mathbb{R}^{d}}\left[ b_{t}\left( x\right) \cdot D\theta \left(
x+y-W_{t}\right) +\mathrm{div}\,b_{t}\left( x\right) \theta \left(
x+y-W_{t}\right) \right] u_{t}\left( x\right) dx.
\end{align*}
Moreover $F_{1}\left( t,y,\omega \right) $ and $F_{0}\left( y\right) $ are
smooth in $y$. Thus It\^{o}-Wentzell-Kunita formula reduces to the classical
It\^{o} formula
\begin{align*}
du_{t}\left( \theta \right) & =f_{1}\left( t,W_{t}\right) dt+\left(
DF_{0}\left( W_{t}\right) +DF_{1}\left( t,W_{t}\right) \right) \cdot dW_{t}
\\
& +\frac{1}{2}\left( \Delta F_{0}\left( W_{t}\right) +\Delta F_{1}\left(
t,W_{t}\right) \right) dt.
\end{align*}
Notice that $D_{i}F_{0}\left( W_{t}\right) +D_{i}F_{1}\left( t,W_{t}\right)
=u_{t}\left( D_{i}\theta \right) $ and $\Delta F_{0}\left( W_{t}\right)
+\Delta F_{1}\left( t,W_{t}\right) =u_{t}\left( \Delta \theta \right) $.
A substitution yields
\begin{align}
du_{t}\left( \theta \right) & =\int_{\mathbb{R}^{d}}\left[ b_{t}\left(
x\right) \cdot D\theta \left( x\right) +\mathrm{div}\,b_{t}\left( x\right)
\theta \left( x\right) \right] u_{t}\left( x\right) dx\,dt
\label{Ito-Wentzell-Kunita} \\
& +\sum_{i=1}^{d}u_{t}\left( D_{i}\theta \right) dW_{t}^{i}+\frac{1}{2}%
u_{t}\left( \Delta \theta \right) dt  \notag
\end{align}
which shows that $u_{t}\left( \theta \right) $ is a semimartingale. For the
same reason, also $u_{t}\left( D_{i}\theta \right) $ is a semimartingale,
for each $i=1,...,d$. Hence the Stratonovich integral $\int_{0}^{t}u_{s}%
\left( D_{i}\theta \right) \circ dW_{s}^{i}$ is well defined and is related
to the It\^{o} integral and to the joint quadratic variation by the formula (see~\cite
{K})
\begin{equation*}
\int_{0}^{t}u_{s}\left( D_{i}\theta \right) \circ
dW_{s}^{i}=\int_{0}^{t}u_{s}\left( D_{i}\theta \right) dW_{s}^{i}+\frac{1}{2}%
\left[ M^{i},W^{i}\right] _{t}
\end{equation*}
where $M_{t}^{i}$ is the martingale part of $u_{t}\left( D_{i}\theta \right)
$, which is equal to
\begin{equation*}
\sum_{j=1}^{d}\int_{0}^{t}u_{s}\left( D^2_{ij}\theta \right) dW_{s}^{j}
\end{equation*}
(put $D_{i}\theta $ in place of $\theta $ in (\ref{Ito-Wentzell-Kunita})
above). Thus $\left[ M^{i},W^{i}\right] _{t}=\int_{0}^{t}u_{s}\left(
D_{i}D_{i}\theta \right) ds$. Summarizing,
\begin{equation*}
\sum_{i=1}^{d}\int_{0}^{t}u_{s}\left( D_{i}\theta \right) \circ
dW_{s}^{i}=\sum_{i=1}^{d}\int_{0}^{t}u_{s}\left( D_{i}\theta \right)
dW_{s}^{i}+\frac{1}{2}\int_{0}^{t}u_{s}\left( \Delta \theta \right) ds.
\end{equation*}
Together with equation (\ref{Ito-Wentzell-Kunita}), this proves that $u$
satisfies the equation of Definition~\ref{def solution}.

\textbf{Step 3}. Let us now prove the converse statement. Let $u$ be an $%
L^{\infty }$-solution of the SPDE~\eqref{SPDE}. Given $\theta \in
C_{0}^{\infty }(\mathbb{R}^{d})$ and $y\in \mathbb{R}^{d}$ let us take the
test function $\theta _{y}\left( x\right) =\theta \left( x+y\right) $ in the
weak formulation of the SPDE. We get
\begin{align*}
& d_t\int_{\mathbb{R}^{d}}u\left( t,x\right) \theta \left( x+y\right) dx \\
& =\int_{\mathbb{R}^{d}}\left[ b\left( t,x\right) \cdot D\theta \left(
x+y\right) +\mathrm{div}\,b\left( t,x\right) \theta \left( x+y\right) \right]
u\left( t,x\right) dx\,dt \\
& +\sum_{i=1}^{d}\int_{\mathbb{R}^{d}}u\left( t,x\right) D_{i}\theta \left(
x+y\right) dx\circ dW_{t}^{i}.
\end{align*}
Consider now the random field $\Theta (t,y)=\int_{\mathbb{R}%
^{d}}u(t,x)\theta (x+y)dx$. Given $t_{1}\in (0,T]$, we apply It\^{o}%
-Wentzell-Kunita formula, in Stratonovich form, to $t\mapsto \Theta \left(
t,y-W_{t}\right) $ for $t\in \left[ 0,t_{1}\right] $ (see~\cite{K}). We get
\begin{equation*}
d\Theta \left( t,y-W_{t}\right) =d\Theta \left( t,z\right)
|_{z=y-W_{t}}-D\Theta \left( t,y-W_{t}\right) \circ dW_{t}
\end{equation*}
\begin{align*}
& =\int_{\mathbb{R}^{d}}\left[ b\left( t,x\right) \cdot D\theta \left(
x+y-W_{t}\right) +\mathrm{div}\,b\left( t,x\right) \theta \left(
x+y-W_{t}\right) \right] u\left( t,x\right) dx\,dt \\
& +\sum_{i=1}^{d}\int_{\mathbb{R}^{d}}u\left( t,x\right) D_{i}\theta \left(
x+y-W_{t}\right) dx\circ dW_{t}^{i} \\
& -\sum_{i=1}^{d}\int_{\mathbb{R}^{d}}u\left( t,x\right) D_{i}\theta \left(
x+y-W_{t}\right) dx\circ dW_{t}^{i}.
\end{align*}
The last two terms coincide. Thus, integrating on $\left[ 0,t_{1}\right] $
we get
\begin{align*}
& \Theta \left( t_{1},y-W_{t_{1}}\right) -\Theta \left( 0,y\right)  \\
& =\int_{0}^{t_{1}}\int_{\mathbb{R}^{d}}\left[ b\left( s,x\right) \cdot
D\theta \left( x+y-W_{s}\right) +\mathrm{div}\,b\left( s,x\right) \theta
\left( x+y-W_{s}\right) \right] u\left( s,x\right) dx\,ds.
\end{align*}
All the terms are a.s. smooth functions of $y$, then taking $y=W_{t_{1}}$
and substituting the definition of $\Theta $ we get the perturbative
equation. The proof is complete.
\end{proof}


\section{Wong-Zakai approximation results}

\label{sec:wk} In this Appendix we prove Wong-Zakai results to motivate the
Stratonovich integral in the SPDE (1). Since this is a side result for the
purpose of this work, here we do not aim at full generality.

Wong-Zakai principle states that the solutions to equations where the noise
is approximated by more regular processes converge to the solution of the
stochastic differential equation with Stratonovich integrals. In contrast
with the classical literature on the subject, here we meet a new difficulty:
the approximating equations could miss uniqueness of solutions, since they
are deterministic transport equations depending on a random parameter and
the regularity of $b$ assumed in this work does not suffice for uniqueness.
The most general statement of Wong-Zakai type, thus, would claim that all
possible $L^{\infty}$ solutions of the approximating equations converge to
the unique solution of the SPDE (we have proved uniqueness for the SPDE in
the previous sections). However, for the deterministic transport equation,
under our assumptions on $b$, there is no control on solutions, no
representation in terms of a flow and it is not even clear how to prove bounds
on them (in spite of the fact that, formally speaking, the $L^{\infty}$ norm
should not increase in time).

Because of these difficulties, we restrict ourselves to two manageable
situations. In the first case we regularize not only the noise but also the
field $b$, so that the approximating equations are well posed. In the second
one we consider a sequence of solutions of the approximating equations which
fulfill a uniform bound (such sequences always exist under our hypotheses).

Let us mention that there exist several results of Wong-Zakai type for
stochastic partial differential equations;\ let us quote only \cite{GS},
\cite{TZ} and \cite{BF} and references therein, as examples of results for
parabolic and transport-type equations. Several works are based on
stochastic characteristics and a Wong-Zakai result for them, which is also
one of our strategies below. However, at the level of characteristics, all
works assume sufficient regularity of coefficients to be able to make
estimates of differences of solutions. Under our weak assumptions on $b$, we
use a different approach, based on the compactness method. The strong well
posedness of the limit equation is the key tool. Also the second theorem,
not based on characteristics, is proved by a compactness argument.

Given a $d$-dimensional Brownian motion $W$ on a probability space $\left(
\Omega,F,P\right) $, let $\left( W_{n}\right) _{n\geq1}$, be a sequence of
processes on the same space such that ($T>0$ is given) $W_{n}$ converges in
probability to $W$ in the topology of $C^{0}\left( \left[ 0,T\right] ;%
\mathbb{R}^{d}\right)$.

An example is
\begin{equation}  \label{eq:wz-example}
W_{n}\left( t\right) =\int_{0}^{\infty}n\theta\left( n\left( t-s\right)
\right) W\left( s\right) ds
\end{equation}
where $\theta$ is a smooth non negative function with support in $\left(
-1,1\right) $ and $\int_{-\infty}^{\infty}\theta\left( r\right) dr=1$.

\begin{theorem}
Assume that hypothesis~\ref{hy1} hold and $\mathrm{div}\,b\in L^{1}\left(
0,T;L_{loc}^{1}(\mathbb{R}^{d})\right) $. Assume $u_{0}\in C_{b}^{0}\left(
\mathbb{R}^{d}\right) $. Let $\{b_{n}\}_{n\geq 1}$ be a sequence of
equibounded (in $\left( t,x,n\right) $) measurable fields such that for a.e.
$t\in \left[ 0,T\right] $ we have $b_{n}\left( t,\cdot \right) \in $ $%
C_{b}^{1}(\mathbb{R}^{d})$ and $b_{n}\left( t,\cdot \right) \rightarrow
b\left( t,\cdot \right) $ uniformly on compact sets. Let $u_{n}(t,x)$ be the
unique $L^{\infty }$ solution of the equation
\begin{equation*}
\frac{\partial u_{n}(t,x)}{\partial t}+b_{n}(t,x)\cdot
Du_{n}(t,x)+\sum_{i=1}^{d}D_{i}u_{n}(t,x)\frac{dW_{n}^{i}(t)}{dt}=0,\quad
u_{n}(0,x)=u_{0}(x)
\end{equation*}
and let $u(t,x)=u_{0}\left( \varphi _{t}^{-1}\left( x\right) \right) $ be
the solution of equation \eqref{SPDE} given by Theorem \ref{them SPDE n1}.
Then, for every $t\geq 0$ and $x\in \mathbb{R}^{d}$, $u_{n}\left( t,x\right)
$ converges in probability to $u\left( t,x\right) $.
\end{theorem}

\begin{proof}
We have $u_{n}\left( t,x\right) =u_{0}\left( \varphi _{n,t}^{-1}\left(
x\right) \right) $ where $\varphi _{n,t}$ is the flow associated to the
random equation
\begin{equation*}
\frac{dX_{n}\left( t\right) }{dt}=b_{n}\left( t,X_{n}\left( t\right) \right)
+\frac{dW_{n}\left( t\right) }{dt}.
\end{equation*}
Thus it is sufficient to prove that $\varphi _{n,t}^{-1}\left( x\right)
\rightarrow \varphi _{t}^{-1}\left( x\right) $ in probability, given $t\geq 0
$ and $x\in \mathbb{R}^{d}$.

The equations satisfied by the inverse flows are entirely similar to the
equations for the direct flows. Thus, just for simplicity of notations, let
us prove that $\varphi _{n,t}(x)\rightarrow \varphi _{t}(x)$ in probability,
given $t\geq 0$ and $x\in \mathbb{R}^{d}$. For shortness, denote $\varphi
_{n,t}(x)$ by $X_{n}(t)$ and $\varphi _{t}(x)$ by $X(t)$ (the initial
condition $x$ is given).

Convergence in law would be classical. Thanks to the idea of \cite{GK}, we
can prove also convergence in probability, due to strong uniqueness for the
limit equation. Let us recall some detail, similar to \cite{GK} (but here we
have to deal also with the approximation of the noise).

Recall Lemma~1.1 from~\cite{GK}. To prove the convergence in probability of $%
X_{n}(t)$ to $X(t)$ it is sufficient to prove the following property. Let $%
\{l_{k}\}_{k\geq 1}$ and $\{m_{k}\}_{k\geq 1}$ be two diverging sequences of
natural numbers. We have to prove that there exist $\{k(j)\}_{j\geq 1}$ such
that the pair $(X_{l_{k(j)}}(t),X_{m_{k(j)}}(t))$ converges in law to a
random element supported on the diagonal $\{(x,y)\in \mathbb{R}^{d}\times
\mathbb{R}^{d}:x=y\}$.

We have
\begin{equation*}
X_{n}(t)=x+\int_{0}^{t}b_{n}(s,X_{n}(s))ds+W_{n}(t).
\end{equation*}
Since the sequence $\{b_{n}\}_{n\geq 1}$ is equibounded in all variables,
the processes $\{\int_{0}^{\cdot }b_{n}(s,X_{n}(s))ds\}_{n\geq 1}$ are
equibounded in $C^{1}([0,T];\mathbb{R}^{d})$. Then since $W_{n}\rightarrow W$
in $C([0,T],\mathbb{R}^{d})$ in probability then the laws $\mu _{n}$ \ of $%
X_{n}$ on $C([0,T];\mathbb{R}^{d})$ are tight, hence precompact by Prohorov
theorem. We do not use explicitly this fact; we have described the argument
in this simple case for later reference.

The argument now is similar to the proof of Theorem 2.4 of \cite{GK}. Given
the two subsequences $\{l_{k}\}_{k\geq 1}$ and $\{m_{k}\}_{k\geq 1}$, let us
repeat the previous argument for the process $%
Z_{k}=(X_{l_{k}},X_{m_{k}},W_{l_{k}},W_{m_{k}})$. By Prohorov theorem there
exist $\{k(j)\}_{j\geq 1}$ such that $Z_{k(j)}$ converges in law to a
probability measure $\nu $ on $C([0,T];\mathbb{R}^{4d})$. By Skorokhod
embedding theorem, there exists a new probability space $(\widetilde{\Omega }%
,\widetilde{F},\widetilde{P})$ and random variables $\widetilde{Z}_{k(j)}=(%
\widetilde{X}_{l_{k(j)}},\widetilde{X}_{m_{k(j)}},\widetilde{W}_{l_{k(j)}},%
\widetilde{W}_{m_{k(j)}})$ with the same laws as $Z_{k(j)}$, and a random
variable $\widetilde{Z}=(\widetilde{X}^{(1)},\widetilde{X}^{(2)},\widetilde{W%
}^{(1)},\widetilde{W}^{(2)})$ with law $\nu $, such that $\{\widetilde{Z}%
_{k(j)}\}_{j\geq 1}$ converges $\widetilde{P}$-a.s. to $\widetilde{Z}$ in
the topology of $C([0,T];\mathbb{R}^{4d})$. It is easy to deduce that $%
\widetilde{W}^{(1)}$ and $\widetilde{W}^{(2)}$ are Brownian motions. Using
the bounded continuous function $\varphi (\gamma ^{(1)},\gamma ^{(2)})=\frac{%
\Vert \gamma ^{(1)}-\gamma ^{(2)}\Vert _{0}}{1+\Vert \gamma ^{(1)}-\gamma
^{(2)}\Vert _{0}}$ on $C([0,T];\mathbb{R}^{d})^{2}$, we see that
\begin{align*}
E[\varphi (\widetilde{W}^{(1)},\widetilde{W}^{(2)})]& =\lim_{j\rightarrow
\infty }E[\varphi (\widetilde{W}_{l_{k(j)}},\widetilde{W}_{m_{k(j)}})] \\
& =\lim_{j\rightarrow \infty }E[\varphi (W_{l_{k(j)}},W_{m_{k(j)}})] \\
& =E[\varphi (W,W)]=0
\end{align*}
hence $\widetilde{W}^{(1)}=\widetilde{W}^{(2)}$. With a similar argument,
one can check that $\widetilde{X}_{l_{k(j)}},\widetilde{W}_{l_{k(j)}}$ are
related by the equation
\begin{equation*}
\widetilde{X}_{l_{k(j)}}(t)=x+\int_{0}^{t}b_{l_{k(j)}}(s,\widetilde{X}%
_{l_{k(j)}}(s))ds+\widetilde{W}_{l_{k(j)}}(t)
\end{equation*}
and similarly for the pair $\widetilde{X}_{m_{k(j)}},\widetilde{W}_{m_{k(j)}}
$. From the $\widetilde{P}$-a.s. convergence in $C([0,T];\mathbb{R}^{d})$ of
all processes, the uniform convergence of $b_{l_{k(j)}}(s,\cdot )$ to $%
b(s,\cdot )$ on compact sets (a.s. in $s$), and the equiboundedness of $%
b_{l_{k(j)}}$, by Lebesgue dominated convergence theorem we get
\begin{equation*}
\widetilde{X}^{(i)}(t)=x+\int_{0}^{t}b(s,\widetilde{X}^{(i)}(s))ds+%
\widetilde{W}(t)
\end{equation*}
for $i=1,2$, where $\widetilde{W}=\widetilde{W}^{(1)}=\widetilde{W}^{(2)}$.
By strong uniqueness for this equation we deduce $\widetilde{X}^{(1)}=%
\widetilde{X}^{(2)}$. Hence $(\widetilde{X}_{l_{k(j)}}(t),\widetilde{X}%
_{m_{k(j)}}(t))$ converges in law to a random element supported on the
diagonal $\{(x,y)\in \mathbb{R}^{d}\times \mathbb{R}^{d}:x=y\}$. This
implies the claim. The proof is complete.
\end{proof}

\begin{remark}
\label{remarkWZ} A more difficult form of Wong-Zakai result would be to
prove that any sequence $\{u_{n}\}$ of $L^{\infty }$ solutions of the (a
priori) not well posed equations
\begin{equation*}
\frac{\partial u_{n}}{\partial t}+b\cdot
Du_{n}+\sum_{i=1}^{d}D_{i}u_{n}\cdot \frac{dW_{n}^{i}}{dt}=0,\quad
u^{n}(0,x)=u_{0}(x)
\end{equation*}
converges to the unique solution $u$ of the SPDE (here we do not regularize $%
b$). Under our assumptions on $b$, it is very difficult to deal with $%
L^{\infty }$ solutions of this `deterministic' equation. For instance, if we
want to perform computations (for proving estimates, comparisons, etc.), we
have to regularize the solution and control the behavior of a commutator,
which is an open problem under this regularity of $b$ and $u^{n}$ (this
problem is the same as proving that weak solutions are renormalizable, in
the sense of DiPerna-Lions~\cite{DiPernaLions}, problem solved under other
conditions on $b$).
\end{remark}

As a partial result towards a general convergence statement we propose the
following theorem dealing with convergence of a particular class of
non-unique solutions to the transport equation.

\begin{theorem}
\label{thm WZ 2}Let $b\in L_{\mathrm{loc}}^{1}([0,T]\times \mathbb{R}^{d};%
\mathbb{R}^{d})$, $\mathrm{div}\,b\in L_{\mathrm{loc}}^{1}([0,T]\times
\mathbb{R}^{d})$ and $u_{0}\in L^{\infty }(\mathbb{R}^{d})$. Let $u^{n}$ be
a sequence of $L^{\infty }(\Omega \times \lbrack 0,T]\times \mathbb{R}^{d})$
functions which are weak solutions of the PDEs
\begin{equation}
\partial _{t}u_{t}^{n}+b_{t}\cdot
Du_{t}^{n}+\sum_{i=1}^{d}D_{i}u_{t}^{n}\partial _{t}W_{n}^{i}(t)=0
\label{eq:wz-2}
\end{equation}
with the same initial condition $u_{0}$. Assume the following conditions

\begin{itemize}
\item[(i)]  $\mathcal{F}_{t}^{n}=\sigma (W_{n}(s):s\leq t)$ converge to $%
\mathcal{F}_{t}=\sigma (W_{s}:s\leq t)$ as $n\rightarrow \infty $ for any $t$
in the sense that $\mathbb{E}[F|\mathcal{F}_{t}^{n}]\rightarrow \mathbb{E}[F|%
\mathcal{F}_{t}]$ almost surely for any bounded r.v. $F$.

\item[(ii)]  The family $\{u^{n}\}$ is equibounded in $L^{\infty }(\Omega
\times \lbrack 0,T]\times \mathbb{R}^{d})$ and $u^{n}$ is $\mathcal{F}^{n}$%
-progressively measurable for any $n$.
\end{itemize}

Then up to extraction of a subsequence still denoted by $u^{n}$ we have weak-%
$\ast $ convergence to $u\in L^{\infty }(\Omega \times \lbrack 0,T]\times
\mathbb{R}^{d})$ satisfying the SPDE $P$-almost surely. Moreover under
Hypotheses~\ref{hy1} and~\ref{hy2} the whole sequence converge to the unique
$L^{\infty }$-solution of the SPDE.
\end{theorem}

Note that by regularization and by compactness it is not difficult to show
that the equibounded family of solutions $\{u^n\}_{n\ge 1}$ of the equations~%
\eqref{eq:wz-2} exists (condition (ii)). Moreover the convergence of the
conditional expectations (condition (i)) holds for our example~\ref
{eq:wz-example} since the mollifier is of compact support (in that case $%
\mathcal{F}^n_t \subseteq \mathcal{F}_{t+c/n}$ for some $c>0$ and the
filtration $\mathcal{F}$ is continuous). Unfortunately, as we already
stressed, we are not able to deal with an arbitrary family of solutions to
the approximating problems~\eqref{eq:wz-2}.

\medskip

\begin{proof}
By an adaptation of Theorem~\ref{lemma:mild} the $u^{n}$ are shown to
satisfy the perturbative equation
\begin{equation*}
u_{t}^{n}(\theta )=u_{0}(\theta (\cdot +W_{n}(t)))
\end{equation*}
\begin{equation*}
+\int_{0}^{t}[b(s,x)\cdot D\theta (x+W_{n}(t)-W_{n}(s))+\mathrm{div}%
\,b(s,x)\theta (x+W_{n}(t)-W_{n}(s))]u_{s}^{n}(x)dxds
\end{equation*}
for any $\theta \in C_{0}^{\infty }(\mathbb{R}^{d})$. By equiboundedness in $%
L^{\infty }$ we can pass to a weakly-$\ast $ convergent subsequence (still
called $u^{n}$) and we have that the weak-$\ast $ limit $u$ is in $L^{\infty
}$ and satisfies
\begin{equation*}
u_{t}(\theta )=u_{0}(\theta (\cdot +W_{t}))+\int_{0}^{t}[b(s,x)\cdot D\theta
(x+W_{ts})+\mathrm{div}\,b(s,x)\theta (x+W_{ts})]u_{s}(x)dxds
\end{equation*}
for a.e. $\left( \omega ,t\right) \in \Omega \times \left[ 0,T\right] $. The
right-hand-side of this equation is a well defined random variable, for
every given $\theta \in C_{0}^{\infty }\left( \mathbb{R}^{d}\right) $ and $%
t\in \left[ 0,T\right] $ (see remark \ref{remark sulla perturbativa}). With
little more work one can see that, as a stochastic process in $t$, it is $P$%
-a.s. continuous. Let us prove this in detail for completeness. We give the
details only for the double integral, that we denote by $I_{t}\left( \omega
\right) $. Take a representative $\widetilde{u}$ of the equivalence class of
$u$, that is a bounded measurable function. For all $\omega $ in a full
measure set $\Omega _{1}$, $\widetilde{u}$ is a bounded measurable function
of $\left( t,x\right) $. Let $\Omega _{2}$ be a full measure set where $W$
is continuous on $\left[ 0,T\right] $. For every $\omega $ in the full
measure set $\Omega _{1}\cap \Omega _{2}$ we have two integrals of the form
\begin{equation*}
\int_{0}^{t}\int_{\mathbb{R}^{d}}a\left( s,x\right) \varphi \left(
x+W_{ts}\left( \omega \right) \right) dxds
\end{equation*}
where $a$ is deterministic and integrable, $\varphi \in C_{0}^{\infty
}\left( \mathbb{R}^{d}\right) $ and $W_{ts}\left( \omega \right) $ is
continuous in $t$ and $s$. Then this integral is continuous in $t$. This
proves that the double integral is a continuous process.

By Theorem~\ref{lemma:mild}, $u$ satisfies the SPDE as soon as we can prove
that the continuous process $u_{t}(\theta )$ is adapted. We have that $%
u_{t}^{n}(\theta )$ converges to $u_{t}(\theta )$ weak-$*$ in $L^{\infty
}(\Omega \times \left[ 0,T\right] )$, hence $\mathbb{E}\left[
\int_{0}^{T}u_{t}^{n}(\theta )f_{t}Gdt\right] $ converges to $\mathbb{E}%
\left[ \int_{0}^{T}u_{t}(\theta )f_{t}Gdt\right] $ for every bounded r.v. $G$
and bounded measurable function $f:\left[ 0,T\right] \rightarrow \mathbb{R}$.
Moreover, $u_{t}^{n}(\theta )$ is $\mathcal{F}_{t}^{n}$-measurable. Then for
any bounded r.v. $G$ we have $\mathbb{E}[u_{t}^{n}(\theta )G]=\mathbb{E}%
[u_{t}^{n}(\theta )\mathbb{E}[G|\mathcal{F}_{t}^{n}]]$. Therefore $%
\int_{0}^{T}f_{t}\mathbb{E}[u_{t}^{n}(\theta )\mathbb{E}[G|\mathcal{F}%
_{t}^{n}]]dt$ converges to $\mathbb{E}\left[ \int_{0}^{T}u_{t}(\theta
)f_{t}Gdt\right] $.

Moreover, we can write
\begin{equation*}
\mathbb{E}[u_{t}^{n}(\theta )\mathbb{E}[G|\mathcal{F}_{t}^{n}]]=\mathbb{E}%
[u_{t}^{n}(\theta )\mathbb{E}[G|\mathcal{F}_{t}]]+\mathbb{E}%
[u_{t}^{n}(\theta )(\mathbb{E}[G|\mathcal{F}_{t}^{n}]-\mathbb{E}[G|\mathcal{F%
}_{t}])]
\end{equation*}
and thus, by the a.s. convergence of $\mathbb{E}[G|\mathcal{F}%
_{t}^{n}]\rightarrow \mathbb{E}[G|\mathcal{F}_{t}]$ together with the
dominated convergence theorem we get that $\int_{0}^{T}f_{t}\mathbb{E}%
[u_{t}^{n}(\theta )\mathbb{E}[G|\mathcal{F}_{t}]]dt$ converges to $\mathbb{E}%
[ \int_{0}^{T}u_{t}(\theta )f_{t}Gdt] $. But the quantity $\int_{0}^{T}f_{t}%
\mathbb{E}[u_{t}^{n}(\theta )\mathbb{E}[G|\mathcal{F}_{t}]]dt$ converges
also to $\int_{0}^{T}f_{t}\mathbb{E}[u_{t}(\theta )\mathbb{E}[G|\mathcal{F}%
_{t}]]dt$.

We thus obtain that $\mathbb{E}[u_{t}(\theta )G]=\mathbb{E}[u_{t}(\theta )%
\mathbb{E}[G|\mathcal{F}_{t}]]$ for any bounded r.v. $G$ showing that $%
u_{t}(\theta )$ is $\mathcal{F}_{t}$-measurable. The proof of the first
claim of the theorem is complete.

If Hypotheses~\ref{hy1} and~\ref{hy2} hold, each weak-$\ast$ convergent
subsequence converge to the unique solution of the SPDE so that the
extraction of a subsequence is not necessary.
\end{proof}


\section{Two additional uniqueness results in $L^\infty$}
\label{sec:fractional}

The aim of this section is to prove some complementary uniqueness results for $L^\infty$ weak solutions of the SPDE obtained extending the key estimates in fractional Sobolev spaces. The first result is the following:

\begin{theorem}
\label{thm:aux-1} Let $d \ge 2$. Assume Hypothesis~\ref{hy1}
and also  that
$
{\mathrm{div}\,}b\in L^{q} (0,T;L^{p}(\mathbb{R}^{d}))
$
 for some $q>2 \ge p > \frac{2d}{ d + 2 \alpha}$. Then there exists a unique weak
$L^{\infty}$-solution $u$ of the Cauchy problem~(\ref{SPDE}) and
$
u(t,x)=u_{0}({\phi}_{t}^{-1}( x))
$.
\end{theorem}

The main interest of this result is due to the fact that
   we can consider some $p$ in the critical interval
   $ (1,2]$ not covered by Hypothesis \ref{hy2}.

Another result deals with an additional hypothesis of Sobolev regularity for $b$ (beside the usual H\"older regularity) which allow to relax the hypothesis on $\mathrm{div}\,b$.

\begin{theorem}
\label{thm:aux-2}
 Assume that $\mathrm{div}\, b\in
L^1_\loc([0,T]\times\RR^d)$ and that
$$
  b \in L^1 (0,T; W_{\loc}^{\theta, 1} (\RR^d) ) \cap L^\infty (0,T; C^{\alpha} (\RR^d) )
$$
with $\alpha>0$, $\theta>0$ and $\alpha+\theta>1$. Then
 there exists a unique weak $L^{\infty}$-solution $u$ of the
Cauchy problem (\ref{SPDE}) and
$
u(t,x)=u_{0}({\phi}_{t}^{-1}( x))
$.
\end{theorem}

The proofs of both theorems follow the proof of Theorem~\ref{thm SPDE n3} using the results below on the commutator and on the regularity of the Jacobian of the flow. Since these results are complementary the details of the proofs are left to the reader. The following commutator estimates follows from Lemma~\ref{chissa}.

\begin{corollary} \label{ci4} Assume
  $ v    \in L_{loc}^{\infty}\left(
\mathbb{R}^{d},\mathbb{R}^{d}\right) $, ${\mathrm{div}\,}v\in
L_{loc}^{1}\left(  \mathbb{R}^{d}\right)$, $g    \in
L_{loc}^{\infty}\left(  \mathbb{R}^{d}\right)$.
\begin{itemize}
\item[(i)] If
there exists $\theta  \in (0,1)$ such that
$
  v \in W_{\loc}^{\theta, 1}(
 \mathbb{R}^{d}, \mathbb{R}^{d}),
$
  then 
$$
\left|\int \mathcal{R}_{\varepsilon}\left[  g,v\right] (  x) \rho(x)
dx \right| \le
 C_r \| g\|_{L^{\infty}_{r+1}} \big( \|
\rho\|_{L^{\infty}_{r}}
 \| {\mathrm{div}\,} v
\|_{_{L_{r+1}^{1}}} +  \, [ \rho ]_{C^{1 - \theta}_{r}} \,
 [v]_{W^{\theta,1}_{r+1}} \big).
$$
\item[(ii)] If  there exists $\alpha  \in (0,1)$ such that
$
   v \in C_{\loc}^{\alpha}( \mathbb{R}^{d}, \mathbb{R}^{d}),
$
 then 
$$
\left|\int \mathcal{R}_{\varepsilon}\left[  g,v\right] (  x) \rho(x)
dx \right| \le
 C_r \| g\|_{L^{\infty}_{r+1}} \big( \|
\rho\|_{L^{\infty}_{r}}
 \| {\mathrm{div}\,} v
\|_{_{L_{r+1}^{1}}} +  \, [ v ]_{C^{\alpha}_{r+1}} \,
 [\rho]_{W^{1- \alpha,1}_r} \big).
$$
\end{itemize}
\end{corollary}
 \begin{proof} We have
$$
  \left|
  \iint g(x') D_x \vartheta_{\eps} ({x-x'}) \,
\big( \rho(x) - \rho(x') \big)\, [v(x)-v(x')] \, dx dx' \right|
$$
$$
\le \frac{\eps^{1 - \theta}}{\eps} \, [ \rho ]_{C^{1 -
\theta}_{r}} \,
 \| g \|_{L^{\infty}_{r+1}}\,
  \frac{1}{\eps^{d}}\,
   \iint_{B(r+1)^2} |D_x \vartheta (\frac{x-x'}{\eps})|  \,
  \frac{|v(x)-v(x')|}{|x-x'|^{\theta + d}} \,
  |x-x'|^{ \theta + d}  dx dx'
$$
$$
 \le
  [ \rho ]_{C^{1 - \theta}_{r}}
\,
 \| g \|_{L^{\infty}_{r+1}}
  \| D \theta \|_{\infty}\, [v]_{W^{\theta,1}_{r+1}}
$$
The second statement has a similar proof.
\end{proof}

These results can be extended to the case in which
commutators are composed with  a flow.

 \begin{lemma}
 \label{dim12} Let $\phi$ be a $C^{1}$-diffeomorphism of
 $\mathbb{R}^{d}$ ($J \phi$ denotes its Jacobian). Assume
  $ v    \in L_{loc}^{\infty}\left(
\mathbb{R}^{d},\mathbb{R}^{d}\right) $, ${\mathrm{div}\,}v\in
L_{loc}^{1}\left(  \mathbb{R}^{d}\right)$, $g    \in
L_{loc}^{\infty}\left(  \mathbb{R}^{d}\right)$.

Then, for any $\rho \in C_r^{\infty}(  \mathbb{R}^{d})$ and any
$R>0$ such that $\mathrm{supp}( \rho\circ \phi^{-1}) \subseteq
B(R)$,
 we have a uniform bound of
 $\int\mathcal{R}_{\varepsilon}\left[
g,v\right] \left(  \phi\left(  x\right)  \right)  \rho\left(
x\right)  dx $ under one of the following conditions:
\begin{itemize}
\item[(i)] there exists $\theta   \in (0,1)$ such
that
$
  v \in W_{\loc}^{\theta, 1}(
 \mathbb{R}^{d}, \mathbb{R}^{d})$, $ J\phi \in
C_{\loc}^{1 - \theta}( \mathbb{R}^{d})
$;

\item[(ii)] there exists $\alpha  \in (0,1)$ such that
$
  J \phi \in W_{\loc}^{1 - \alpha, 1}(
 \mathbb{R}^{d} ) $, $ v \in
C_{\loc}^{\alpha}( \mathbb{R}^{d}, \mathbb{R}^{d})
$.
\end{itemize}
Moreover, under one of the previous conditions, we also have
\[
\lim_{\varepsilon\rightarrow0}\int\mathcal{R}_{\varepsilon}\left[
g,v\right] \left(  \phi\left(  x\right)  \right)  \rho\left(
x\right)  dx=0.
\]
\end{lemma}
\begin{proof}
 By a change of variables $ \int\mathcal{R}_{\varepsilon}[
g,v] (  \phi(  x) )  \rho( x) dx=\int\mathcal{R}_{\varepsilon}[ g,v]
(  y)  \rho_{\phi}(  y) dx $ where the function
 $ \rho_{\phi}(  y)  =\rho(  \phi^{-1}( y) )
 J\phi^{-1}(  y)
 $
 has the support strictly contained in the ball
of radius $R$. Clearly, $
 \|  \rho_{\phi}\| _{L_{R}^{\infty}}\leq \| \rho\|
_{L_{r}^{\infty}}\| J\phi^{-1}\| _{L_{R}^{\infty}}. $
To prove the result, we  have to check that Corollary~\ref{ci4}
can be applied with $\rho_{\phi}$ instead of $\rho$.

(i) To apply Corollary~\ref{ci4}~(i), we need to check that
 $ \rho_{\phi} \in C^{1 -\theta}_{loc}$. This follows since
$$
[\rho_{\phi} ]_{ C^{1 -\theta}_{R}} \le
 \| J\phi^{-1}\|
_{L_{R}^{\infty}} \, [\rho ({\phi}^{-1} (\cdot) ) ]_{ C^{1
-\theta}_{R}} \, + \, \| \rho \|_{L_{r}^{\infty}}
 [J\phi^{-1}]_{C^{1 -\theta}_{R}}
$$
$$
\le  \| D\phi^{-1}\| _{L_{R}^{\infty}} \,  \| D\rho\|
_{L_{r}^{\infty}} \, [ D\phi^{-1}]_{C^{1 -\theta}_{R}}
 \, + \, \| \rho \|_{L_{r}^{\infty}}
 [D\phi^{-1}]_{C^{1 -\theta}_{R}}.
$$
and the bound follows.

\medskip (ii) To apply Corollary~\ref{ci4}~(ii), we
  need to check that
 $ \rho_{\phi} \in W^{1- \alpha,1}_{loc}$: first
$$ [\rho_{\phi} ]_{W^{1 - \alpha, 1}_R }
\le  \| J\phi^{-1}\| _{L_{R}^{\infty}} [\rho \circ \phi^{-1}]_{W^{1- \alpha,1}_R}
 \, + \,
 [J \phi^{-1} ]_{W^{1- \alpha,1}_R} \, \| \rho\|
_{L_{r}^{\infty}}
$$
and since
$$
[\rho \circ \phi^{-1}]_{W^{1- \alpha,1}_R} \le \|D (\rho \circ \phi^{-1})\|_{L^1_R} \le \|D\rho\|_{L^1_r} \|D \phi^{-1}\|_{L^\infty_R} \,;
$$
we find
$$
[\rho_{\phi} ]_{W^{1 - \alpha, 1}_R } \, \le \,  C_R \|D\rho\|_{L^1_r} \|D \phi^{-1}\|_{L^\infty_R} \| J\phi^{-1}\| _{L_{R}^{\infty}} + [J \phi^{-1} ]_{W^{1-
\alpha,1}_R} \, \| \rho\| _{L_{r}^{\infty}}
$$
and the bound follows.
\end{proof}

Finally the next theorem extends the analysis of the Jacobian of the flow presented in Section~\ref{sec:jacobian} and links the regularity condition on $J\phi$ required in  Lemma~\ref{dim12} (ii) to the assumption on the divergence of $b$ stated in Theorem~\ref{thm:aux-1}.

\begin{theorem}
\label{iac11} Let $d\ge 2$. Assume  Hypothesis~\ref{hy1}  and  the existence of  $p\in (\frac{2d}{ d + 2 \alpha},2]$ and $q>2$ such that
$
 {\mathrm{div}\,}b\in L^{q} (0,T;L^{p}(\mathbb{R}^{d})) .
$
Then,
 for any
 $r>0$,
 $
  J\phi \in
 L^{p}( 0,T;W_{r}^{1 - \alpha, \, p} )$, $P$-a.s.
\end{theorem}
\begin{proof} The first part of the proof is similar to
 the one of Theorem~\ref{iac}. Indeed Step~1 can be carried on thanks to the  chain rule for fractional Sobolev
 spaces:
 if $f:\mathbb{R}^{d}\rightarrow \mathbb{%
R}$ is a continuous function, of class $W_{loc}^{1-\alpha,p}(\mathbb{R}^{d})$ and $%
g:\mathbb{R^d}\rightarrow \mathbb{R}$ is a $C^{\infty }$ function,
then $ g\circ f\in W_{loc}^{1-\alpha,p}(\mathbb{R}^{d})$  and
\begin{equation*}
 [ (g\circ f)]^p_{W_{r}^{1- \alpha, p}}
 \leq \left( \sup_{x\in B(r)}\left|
g^{\prime }(f(x))\right| \right) ^{p} [  f]^p_{W_{r}^{1- \alpha,
p}},
\end{equation*}
for every $r>0$. The modification of Step~2 does not pose any problem, so we only consider the last steps of the proof.

\medskip \textbf{Step 3.} To prove the assertion it is enough
 to check that the family $\left(
\psi _{\varepsilon }\right) _{\varepsilon >0}$ is bounded in
 $L^{p}(\Omega \times (0,T);W_{r}^{1-\alpha,p})$.

Indeed, once we have proved this fact, we can extract from
 the previous
sequence
 $\psi _{\varepsilon _{n}}$ a subsequence which converges
  weakly in $L^{p}(\Omega
\times (0,T);W_{r}^{1-\alpha,p})$ to some $\gamma $. This in
particular implies that such subsequence  converges weakly in
$L^{p}(\Omega \times (0,T),L_{r}^{p})$ to  $\gamma $ so we must have that $\gamma = J \phi$.

We introduce the following Cauchy problem, for $\eps \ge 0$,
\begin{equation}
\left\{ \begin{aligned} \frac{\partial F^{\eps}}{\partial t}+
\frac{1}{2} \Delta F^{\eps} + D F^{\eps} \cdot b^{\eps}
={\mathrm{div}\, }b^{\eps}, \;\;
\; t \in [0,T[ \\ F^{\eps}(T,x)=0,\;\;\; x \in \mathbb{R} ^d. \end{aligned}%
\right.   \label{equa1}
\end{equation}
  This  problem
 has a unique solution $F^{\eps}$ in the space
 $L^{q}(0,T;W^{2,p}(\mathbb{ R}^{d})$. Moreover,
there exists a positive constant $C=C(p,q,
 d,T,\Vert b\Vert _{\infty
})$ such that
\begin{equation}
\Vert F^{\eps}\Vert _{L^{q}(0,T;W^{2,p}(\mathbb{ R}^{d})) }\leq
C\Vert
 \mathrm{div \,}b \Vert _{L^{q}(0,T;L^{p}(\mathbb{
R}^{d}))},  \label{bond5}
\end{equation}
for any $\eps \ge 0$.  This result  can be proved by using
\cite[Theorem 1.2]{Kr1}  and repeating the
 argument  of the proof
 in \cite[Theorem 10.3]{Kry-Ro}. This argument  works without
difficulties in the present case in which $b$ (and so $b^{\eps}$)
is globally bounded and $\mathrm{div \,}b \in
L^{q}(0,T;L^{p}(\mathbb{%
R}^{d}))$ with  $p, q \in (1, + \infty)$.

From the previous result we can also deduce, since
we are assuming $q>2$, that  $F^{\eps} \in C( [ 0,T] ;W^{1,p}(
\mathbb{R}^{d}) )$, for any $\eps \ge 0$, and moreover there
exists a positive constant $C$ $=C(p,q$ $
 d,T, \Vert b \Vert _{\infty
})$ such that
\begin{equation} \label{df}
\sup_{t\in\left[ 0,T\right] }\left\|   F^{\eps}(  t,\cdot)
\right\| _{W^{1,p}\left( \mathbb{R}^{d}\right)  } \leq C \Vert
 \mathrm{div \,}b \Vert _{L^{q}(0,T;L^{p}(\mathbb{
R}^{d}))}.
\end{equation}
We only give a  sketch of proof of \eqref{df}.
 Define $u^{\eps} (t,x) = F^{\eps} (T-t, x)$; we
have the explicit formula
$$
 u^{\eps}(t,x) = \int_0^t P_{t-s}g^{\eps}(s, \cdot)(x) ds,
$$
where $(P_t) $ is the heat semigroup and $g^{\eps}(t,x)
 = D u^{\eps} (t,x) \cdot b^{\eps} (T-t,x)
 - {\mathrm{div}\, }b^{\eps} (T-t,x)$. We get, since $q>2$
  and $q' =\frac{q}{q-1}<2$,
$$
\| D_x u^{\eps}(t, \cdot)\|_{L^p} \le  c\int_0^t
 \frac{1}{(t-s)^{1/2}} \| g^{\eps}(s, \cdot)\|_{L^p} ds
\le C
 \Big( \int_0^T
 \frac{1}{s^{q'/2}}  ds
  \Big)^{1/q'} \, \Big( \int_0^T
   \| \mathrm{div}\,  b (s, \cdot)\|_{L^p}^{q} ds \Big)^{1/q}
$$
  and so \eqref{df} holds.

  \medskip

  Using It\^{o} formula we find (remark that
 $F^{\eps}(t, \cdot) \in C^{2}_b(\RR^d)$)
\begin{equation} \label{ci51}
F^{\varepsilon }\left( t,{\phi }_{t}^{\varepsilon }\left( x\right)
\right)
-F^{\varepsilon }\left( 0,x\right) -\int_{0}^{t}DF^{\varepsilon }\left( s,{%
\phi }_{s}^{\varepsilon }\left( x\right) \right) \cdot dW_{s}\int_{0}^{t}%
\mathrm{div}b^{\varepsilon }\left( s,{\phi }_{s}^{\varepsilon
}\left( x\right) \right) ds = \psi _{\varepsilon }(t,x ).
\end{equation}
Since we already know that $\left( \psi _{\varepsilon }\right)
_{\varepsilon >0}$ is bounded in $L^{p}(\Omega \times
(0,T),L_{r}^{p})$ and since $p\le 2$, to verify that $\left( \psi _{\varepsilon
}\right) _{\varepsilon >0}$ is bounded in $L^{p}(\Omega \times
(0,T);W_{r}^{1-\alpha,p})$, it is enough to prove that
$
E \int_0^T [\psi _{\varepsilon }(t, \cdot)]_{W_{r}^{1-\alpha,2}}^2
dt \le C,
$
for any $\eps >0$. We give details only for the
most difficult term $\int_{0}^{t}D^{}F^{\varepsilon }(s,{\phi }
 _{s}^{\eps}(x))  dW_{s}$ in \eqref{ci51}. The $F(0,x)$ term can be controlled using~\eqref{df} and the others are of easier estimation.
 We show  that  there exists  a constant $C>0$ (independent on  $\eps$)
 such that
\begin{equation}
E \int_{0}^{T} dt \left [ \int_{0}^{t}D^{}F^{\varepsilon }\left(
s,{\phi }_{s}^{\eps}\left( \cdot\right) \right) dW_{s}\right
]_{W_{r}^{1-\alpha,2}}^{2} \, \leq C \label{cia2}
\end{equation}
 We have
$$
E\left[
\int_{0}^{T}dt\int_{B(r)}\int_{B(r)}\frac{|\int_{0}^{t}(DF^{\varepsilon
}\left( s,{\phi}^{\epsilon}_s \left( x\right) \right)
-DF^{\varepsilon }\left( s, {\phi}^{\epsilon}_{s}\left( x^{\prime
}\right) \right) )dW_{s}|^{2}}{|x-x^{\prime }|^{(1 - \alpha)2 + d
}}dx\,dx^{\prime }\right]
$$
$$
= \int_{0}^{T}    \int_{B(r)}\int_{B(r)} E
\int_{0}^{t}\frac{|DF^{\varepsilon }\left(
s,{{\phi}^{\epsilon}}_{s}\left( x\right) \right)
-DF^{\varepsilon }\left( s,{%
{\phi}^{\epsilon}}_{s}\left( x^{\prime }\right) \right)
|^{2}}{|x-x^{\prime }|^{(1 - \alpha)2 + d } }ds \,  dx\,dx^{\prime }
,
$$
\begin{equation*}
=  E \int_{0}^{T}dt  \int_{0}^{t}ds
\int_{B(r)}\int_{B(r)}\frac{|DF^{\varepsilon }\left(
s,{{\phi}^{\epsilon}}_{s}\left( x\right) \right)
-DF^{\varepsilon }\left( s,{%
{\phi}^{\epsilon}}_{s}\left( x^{\prime }\right) \right)
|^{2}}{|x-x^{\prime }|^{(1 - \alpha)2 + d } }dx\,dx^{\prime }
\end{equation*}
\begin{equation*}
\leq TE\left[ \int_{0}^{T}ds\int_{B(r)}\int_{B(r)}\frac{|DF^{\varepsilon }\left( s,%
{{\phi}^{\epsilon}}_{s}\left( x\right) \right) -DF^{\varepsilon
}\left( s, {\phi}^{\epsilon} _{s}\left( x^{\prime }\right) \right)
|^{2}}{|x-x^{\prime }|^{(1 - \alpha)2 + d} }dx\,dx^{\prime }\right]
,
\end{equation*}
\begin{equation*}
\leq TE\int_{0}^{T} [ DF^{\varepsilon }(s,{{\phi}^{\epsilon}}
_{s}(\cdot ))]_{W_{r}^{1 - \alpha, 2}}^{2} \, ds
\end{equation*}
By the Sobolev embedding the $W^{1 - \alpha,2}_r$-seminorm  can
be controlled by the norm in $W^{1,p}_r$ if
$$
 1 - \frac{d}{p} \ge (1-\alpha) - \frac{d}{2}.
$$
This means if $p \ge \frac{2d}{ d + 2 \alpha}$.
 Then we consider $p_{1} $ such that
$
p>{p_{1}}> \frac{2d}{ d + 2 \alpha}
$
 and show that
\begin{equation} \label{f7}
E\int_{0}^{T}\Vert DF^{\varepsilon }(s,{{\phi}^{\epsilon}}
_{s}(\cdot ))\Vert _{W_{r}^{1,p_{1}}}^{2}ds \le C < \infty,
\end{equation}
where $C$ is independent on $\eps$.

\medskip \noindent \textbf{Step 4.}
 To obtain \eqref{f7} we estimate
\begin{equation*}
E\int_{0}^{T}ds\Big (\int_{B(r)}|D^{2}F^{\varepsilon }\left( s,{{{\phi}^{\epsilon}} }%
_{s}\left( x\right) \right) D{{{\phi}^{\epsilon}} }_{s}\left( x\right) |^{p_{1}}dx\Big )^{%
\frac{2}{p_{1}}}
\end{equation*}
A similar term has been already estimated in the proof of Theorem
\ref{iac}. Since
\begin{equation*}
\int_{B(r)}\left( \int_{0}^{T}E\left[ \left| D{{{\phi}^{\epsilon}}
}_{s}\left( x\right) \right| ^{r}\right] ds\right) ^{\gamma
}dx<\infty ,
\end{equation*}
for every $r,\gamma \geq 1$ (see (\ref{bound})), by the H\"{o}lder
inequality, it is sufficient to prove that
\begin{equation*}
\int_{0}^{T}E\left[ \left( \int_{B(r)}\left| D^{2}F^{\varepsilon }\left( s,{%
{{\phi}^{\epsilon}} }_{s}
\left( x\right) \right) \right| ^{p^{}}dx\right) ^{\frac{2}{%
p^{ }}}\right] dt \le C <\infty.
\end{equation*}
  We have
\begin{align*}
& \int_{0}^{T}E\left[ \left( \int_{B(r)}\left| D^{2}F^{\varepsilon }
\left( s,{%
{{\phi}^{\epsilon}} }_{s}\left( x\right) \right) \right|^{p^{}}dx\right) ^{\frac{2}{%
p^{ }}}\right] dt \\
& =E\left[ \int_{0}^{T}ds\left( \int_{{{\phi}^{\epsilon}}
_{s}(B(r))}\left| D^{2}F^{\varepsilon }\left( s,{y}\right) \right|
^{p^{}}J({{\phi}^{\epsilon}}
_{s})^{-1}(y)dy\right) ^{\frac{2}{p^{ }}}\right]  \\
& \leq \sup_{s \in [0,T] , \, y \in \RR^d}
   E[ J (\phi_s^{\eps})^{-1}\, (y)]^{2/p}
\int_{0}^{T}\Big (\int_{\mathbb{R}^d}\left| D^{2}F^{\varepsilon
}\left( s,{y}\right) \right| ^{p}dy\Big )^{\frac{2}{p}}\, \le C <
\infty,
\end{align*}
where, using  the results of Section~\ref{sec:ito} and the bound \eqref{bond5},
$C$ is independent on $\eps>0.$ The proof  is complete.
\end{proof}

\end{document}